\documentclass[11pt,reqno, a4paper]{amsart}
\usepackage{graphicx} 
\usepackage{mathrsfs}
\usepackage{color}
\usepackage{epsfig}
\usepackage{subcaption}
\usepackage[outdir=./]{epstopdf}
\usepackage{enumerate}
\usepackage{amssymb,amsmath,amsthm,graphicx,amsfonts}
\usepackage{hyperref}
\usepackage[margin=1.0in]{geometry}
\usepackage[ruled,linesnumbered]{algorithm2e}

\usepackage{standalone}
\usepackage{tikz} 
\usepackage{pgfplots}
\usepackage{xstring}
\usepackage{subcaption}
\usepackage{esint}

\newcommand{\projectroot}{.} 
\newcommand{\datapath}{\projectroot/data}

\usepackage{titlesec} 
\titleformat{\section}{\vskip10pt\large\bfseries}{\thesection.}{0.5em}{\centering\vspace{5pt}}
\titleformat{\subsection}{\vskip10pt\normalsize\bfseries}{\thesubsection.}{0.5em}{}

\newtheorem{theorem}{\hspace{1mm}Theorem}[section]
\newtheorem{lemma}[theorem]{\hspace{1mm}Lemma}
\newtheorem{remark}[theorem]{Remark} 
\newtheorem{assumption}[theorem]{\hspace{1mm}Assumption} 
\newtheorem{ciao}[theorem]{\hspace{1mm}Corollary} 

\def\R{\mathcal{R}}

\def\d{\mathrm{d}}
\def\P{\mathcal{P}}

\allowdisplaybreaks

\newcommand{\taufac}{C_{\textup{scal}}}
\newcommand{\tauref}{\tau_{\textup{ref}}}
\newcommand{\hspaceref}{h_{\textup{ref}}}
\newcommand{\numexpheight}{\mathcal{H}}

\newcommand{\massmatrix}{\mathbf{M}_h}
\newcommand{\stiffmatrix}{\mathbf{A}_h}
\newcommand{\solutionvec}{\mathbf{u}_h}
\newcommand{\solutiondtvec}{\mathbf{v}_h}
\newcommand{\loadvectorn}{\mathbf{f}_h^n}

\numberwithin{equation}{section}
\numberwithin{figure}{section}
\numberwithin{table}{section}

\usepackage{xcolor}

\begin{document}

\title[]{Finite element exponential integration for rough solutions of nonlinear wave equations. Part I: Dirichlet boundary conditions on polygonal and polyhedral domains}

\author[]{Jiachuan Cao}
\address{Institute for Applied and 
	Numerical Mathematics, Karlsruhe Institute of Technology, Englerstr.~2, 76131 Karlsruhe, Germany}
\email{\{jiachuan.cao,benjamin.doerich,marlis.hochbruck\}@kit.edu}

\author[]{Benjamin D\"{o}rich}

\author[]{Marlis Hochbruck}

\author[]{Buyang Li}
\address{Department of Applied Mathematics, The Hong Kong Polytechnic University,
	Hung Hom, Hong Kong}
\email{buyang.li@polyu.edu.hk}

\subjclass[2010]{65M12, 65M15, 35Q55}


\keywords{nonlinear wave equation, low regularity, Dirichlet boundary conditions, error estimates}
\thanks{
	{ The first three authors acknowledge funding by the Deutsche Forschungsgemeinschaft (DFG, German Research Foundation) - 
		Project-ID 258734477 - SFB 1173.
		The work of B. Li was partially supported by the National Natural Science Foundation of China (project no. 12525111) and Hong Kong Research Grants Council (project no. 15306123). 
	}
}

\maketitle
\vspace{-20pt}

\begin{abstract}\noindent
{\small 
	We study a fully discrete scheme for nonlinear wave equations on general bounded polygonal/polyhedral domains with initial data $(u^0,v^0)\in H^\gamma(\Omega)\times H^{\gamma-1}(\Omega)$, $0<\gamma\le 1$, subject to the natural compatibility condition associated with the homogeneous Dirichlet boundary condition. The scheme combines an exponential Euler time integrator with a
	finite element spatial discretization. In contrast to existing low-regularity error analyses, which are mostly based on Fourier spectral discretizations, our approach applies to general bounded domains and finite element spatial discretizations. We prove rigorous error estimates for low-regularity solutions. The analysis is based on a frequency decomposition of the underlying elliptic operator, used solely as an analytical regularization device and not in the actual implementation, which allows low-regularity techniques to be extended beyond the Fourier spectral framework. The results also indicate that higher-order finite element methods remain advantageous in spatial approximation even for solutions of limited Sobolev regularity. Numerical experiments on different domains and with different polynomial degrees confirm the predicted convergence behavior.
	
}
\end{abstract}


\setlength\abovedisplayskip{4.5pt}
\setlength\belowdisplayskip{4.5pt}


\section{Introduction}

Nonlinear wave equations of the form
\begin{equation}\label{eq:NLW}
	\left\{
	\begin{array}{ll}
		\partial_{tt}u - \Delta u = f(u), & \text{in } \Omega \times (0, T], \\[1.5mm]
		u = 0, & \text{on } \partial\Omega \times (0,T],\\[1.5mm]
		u|_{t=0} = u^{0}, \quad \partial_{t}u|_{t=0} = v^{0}, & \text{in } \Omega,
	\end{array}
	\right.
\end{equation}
arise in a broad range of applications involving nonlinear wave propagation, including continuum mechanics, nonlinear optics, and materials science. 
Their accurate numerical approximation is therefore of considerable practical and theoretical interest.

For sufficiently smooth solutions, a wide variety of time discretization methods are available for nonlinear wave equations, including trigonometric and exponential integrators, splitting methods, symplectic schemes, and classical finite difference or Runge--Kutta type methods; see, for example, \cite{BFS2022,CHSS2020,HL2021,HL1999,WW2019-IMA}. 
Under appropriate regularity assumptions, these methods often achieve their expected convergence orders and exhibit good long-time behavior. 
However, in many situations of practical and mathematical interest, the solution possesses only limited Sobolev regularity. 
This may be caused by rough initial data, nonsmooth geometries, singular forcing, or boundary effects. 
In such regimes, standard error analyses are no longer applicable, and classical high-order schemes may suffer substantial order reduction or require assumptions incompatible with the actual solution regularity. 
Related regularity-dependent convergence results for time discretizations of nonsmooth problems have also been obtained in more abstract settings; see, for example, \cite{brenner1980rational,lubich1994multistep,WE2016}.

In recent years, substantial progress has been made in low‑regularity integrators for deterministic dispersive equations.  
Important advances have been reported for the Korteweg--de Vries (KdV) equation~\cite{Hofmanova-Schratz-2017,li2025unfiltered}  
and the nonlinear Schr\"{o}dinger equation~\cite{feng2025explicit,ORS20,OS18}, and other important dispersive equations~\cite{ji2025filtered,RS21,schratz2021low}.
For nonlinear wave equations, second‑order convergence in the $H^1\times L^2$ norm has been proved under additional regularity assumptions~\cite{RS21}.  
Exponential‑type low‑regularity integrators based on integration‑by‑parts formulations have also shown efficient time stepping for finite‑energy solutions~\cite{BucDH21}.  
Later, Cao, Li, Lin, and Yao~\cite{CLLY2024} developed high-frequency-recovered exponential integrators based on a frequency-decomposition technique, and established convergence rates of \(2\gamma\) and \(1.5\gamma\) in one and two spatial dimensions, respectively, for initial data
$(u^0,v^0)\in H^{\gamma}\times H^{\gamma-1}$ with $0<\gamma\le1$.  
More recently, Ruff and Schnaubelt~\cite{ruff2025error,RufS25} established convergence results for Lie‑ and Strang‑splitting methods for nonlinear wave equations with finite‑energy solutions in three spatial dimensions.  
However, these analyses rely strongly on Fourier techniques, spectral discretizations, or periodic boundary conditions.

This reliance substantially limits the range of geometries and boundary conditions that can be treated. 
In particular, Fourier-based approaches are not directly suited to general bounded domains, variable coefficients, or nontrivial boundary conditions. 
From the perspective of practical computation, finite element methods are the natural framework in such situations. 
They offer geometric flexibility, compatibility with variational formulations, and a well-developed approximation theory on nonsmooth domains. 
Despite these advantages, a rigorous low-regularity error analysis for fully discrete finite element approximations of nonlinear wave equations remains largely open.

In this paper, we study the nonlinear wave equation \eqref{eq:NLW} on a bounded polygonal/polyhedral domain \(\Omega\subset\mathbb{R}^d\), \(d=2,3\), with low-regularity initial data
\begin{align}\label{initial}
	(u^0,v^0)\in H^\gamma(\Omega)\times H^{\gamma-1}(\Omega),
	\qquad 0<\gamma\le 1,
\end{align}
with the appropriate compatibility condition associated with the homogeneous Dirichlet boundary condition imposed when $\gamma\ge \frac12$.
We develop and analyze a fully discrete finite element method with an exponential Euler time discretization for this low-regularity regime. 
Our objective is to transfer low-regularity techniques, previously developed mainly in Fourier-based settings, to a finite element framework on more general domains. The precise functional setting is formulated in Section~\ref{sec:main_results} in terms of the fractional Hilbert scale generated by the Dirichlet Laplacian.


To overcome the challenges posed by the low regularity of the solution, we combine an exponential-type time discretization with a frequency-localized analytical framework adapted to finite element spaces. 
A central ingredient is a frequency projection operator that decomposes the exact solution into low- and high-frequency parts. 
This device plays the role of an analytical smoothing mechanism, analogous in spirit to ideas used in Fourier-based low-regularity analysis, but it is employed here only as a tool in the proof. 
The numerical method itself remains a standard finite element discretization coupled with exponential time stepping and does not require the computation of eigenfunctions or the implementation of spectral truncations. 

The main contributions of this paper can be summarized as follows.
\begin{itemize}
	\item
	We develop a new analytical framework for deriving error estimates in the weak norm \(L^2(\Omega)\times H^{-1}(\Omega)\) for the fully discrete exponential Euler finite element scheme on bounded polygonal/polyhedral domains with initial data satisfying \eqref{initial}. In abstract form, the scheme reads
	\[
	U_h^{n+1}
	=
	e^{\tau L_h}U_h^n
	+
	\tau\varphi_1(\tau L_h)\P_h F(U_h^n).
	\]
	A key step is a comparison estimate between \(e^{tL}\Pi_r\) and \(e^{tL_h}\P_h\Pi_r\), which links the frequency-localized auxiliary problem to the finite element discretization. Together with low-regularity consistency and stability estimates, this comparison yields a fully discrete convergence result under the low-regularity assumption \eqref{initial}. To the best of our knowledge, this provides the first fully discrete finite element convergence analysis for nonlinear wave equations at this level of regularity.
	
	\item
	The resulting weak-norm estimate exhibits a degree-dependent spatial convergence rate. More precisely, for finite element spaces of polynomial degree \(k\), the fully discrete error satisfies
	\[
	\sup_{0\le n\le T/\tau}
	\|U(t_n)-U_h^n\|_{L^2(\Omega)\times H^{-1}(\Omega)}
	\lesssim
	\tau h^{\frac{\ell_k}{\ell_k+1}(\gamma-1)}
	+
	h^{\frac{\ell_k}{\ell_k+1}\gamma},
	\]
	where
	\[
	\ell_k =
	\begin{cases}
		k+1, & k=1,\\
		k+2, & k\ge2.
	\end{cases}
	\]
	In particular, the spatial contribution is of order 	\(h^{\frac{\ell_k}{\ell_k+1}\gamma}\). Thus the convergence exponent improves from \(2\gamma/3\) for linear elements to \(4\gamma/5\) and \(5\gamma/6\) for quadratic and cubic elements, respectively. This improvement is obtained under the same low-regularity assumption \eqref{initial}, without requiring any additional Sobolev regularity as the polynomial degree is increased. Thus polynomial enrichment improves the asymptotic spatial convergence rate itself, rather than merely the error constant, even for rough solutions of nonlinear wave equations.
	
	\item
	Numerical experiments on different domains and for polynomial degrees \(k=1,2,3\) confirm the predicted degree-dependent rates. They also show a clear improvement in the error versus degrees-of-freedom comparison for higher-order finite elements in the low-regularity regime.
\end{itemize}

The remainder of this paper is organized as follows.
In Section~\ref{sec:main_results}, we introduce the frequency-localized analytical framework for the abstract formulation of the nonlinear wave equation and present the time semidiscrete problem, which serves as a bridge for the analysis of the fully discrete scheme. 
At the end of this section, we state the main convergence results. 
Section~\ref{sec:preliminary} collects several auxiliary results that support the proofs of the main theorems, including error estimates between the continuous and discrete solution operator semigroups. 
In Section~\ref{sec:proof}, we combine these ingredients to establish the main theoretical results. 
Finally, Section~\ref{sec:numerical_experiments} reports a series of numerical experiments for nonlinear wave equations on different domains and under various regularity assumptions. 
These experiments verify the accuracy of the proposed method and confirm the theoretical findings.


\noindent \textbf{Notations.} For the sake of notational simplicity, we write \( A \lesssim B \) (or equivalently \( B \gtrsim A \)) to denote that
\(A \leq C B\) for some constant \( C > 0 \). The constant \( C \) may vary from line to line, but it is always independent of the time step \(\tau\), the mesh size \(h\), the frequency cut-off constant \(r\) and the bounds of the numerical solution. 
Similarly, we use \( A \sim B \) to denote the equivalence \(C^{-1} B \leq A \leq C B\)
for some constant \( C > 0 \). In other words, \( A \sim B \) is equivalent to \( A \lesssim B \lesssim A \).  We also use the notation \(C(\cdot)\) to denote a positive constant whose value may change from line to line but depends only on the quantities specified in the parentheses.

\section{Mathematical setting and the main result}\label{sec:main_results}

In this section, we introduce the abstract functional framework used for the analysis and numerical approximation of nonlinear wave equations on polygonal/polyhedral domains. To handle solutions with low regularity or strong singularities, we combine an exponential Euler integrator in time with a finite element discretization in space. Within this framework, we state the main convergence result under minimal regularity assumptions on the exact solution and introduce a frequency-localized auxiliary scheme that serves as a key analytical bridge in the error analysis.

\subsection{Abstract formulation and functional framework}
We begin with the weak formulation of a general nonlinear wave equation, expressed as
\begin{align}\label{eq:wake_form}
	\langle u_{tt}, v\rangle_{H} + a(u, v) = \langle f(u), v\rangle_{H},
\end{align}
where \(\langle\cdot,\cdot\rangle_{H}\) denotes a bilinear form defined on a Hilbert space \(H\), and \(a(\cdot,\cdot)\) is another bilinear operator on a Hilbert space \(V\), with \(V\) densely embedded in \(H\).

\begin{remark}\upshape
	Throughout this work, we restrict our attention to problems with homogeneous Dirichlet boundary conditions. 
	Within this setting, we have \(H = L^2(\Omega)\), endowed with the standard inner product
	\[
	\langle u,v\rangle_{H} = \int_{\Omega} u v \, dx.
	\]
	Moreover, \(V = H_0^1(\Omega)\),
	the associated bilinear form is given by
	\[
	a(u,v) = \int_{\Omega} \nabla u \cdot \nabla v \, dx.
	\]
\end{remark}

Using the {Riesz} representation theorem for linear operators on Hilbert spaces, we reformulate the nonlinear wave equation into an operator form involving an unbounded operator \(A\):  
\[
A : D(A) \to H, \quad \langle A v, w\rangle_{H} = a(v, w) \quad \text{for all } v \in D(A), \, w \in V,
\]
with
\[
D(A) = \left\{ v \in V \;\middle|\; \exists C = C(v) > 0 \;\text{such that}\; \forall w \in V: |a(v, w)| \leq C \|w\|_H \right\}.
\]
The original nonlinear wave equation is then reduced to the abstract form:
\[
u_{tt} + A u = f(u).
\]

We pose the following assumption on the operator $A$:
\begin{assumption}\label{assump:A}  \

	\begin{itemize}
		\item[(a)] The operator \(A\) is self-adjoint and strictly positive, satisfying the coercivity condition 
		\(\langle Av, v\rangle_{H} \geq C_A \langle v, v\rangle_{H}\), where \(C_A > 0\) for all \(v \in D(A)\).
		\item[(b)] \(A\) has a dense domain \(D(A)\) in \(H\) and a purely point spectrum.
	\end{itemize}
\end{assumption}

We denote by \(\{\lambda_{\alpha}, w_{\alpha}\}_{\alpha \in I}\) the eigenpairs of the operator \(\Lambda = \sqrt{A}\), 
where \(I\) is a countable index set. 
The associated orthonormal eigenvectors \(e_{\alpha} = w_{\alpha} / \|w_{\alpha}\|_{H}\) form a complete orthonormal basis of \(H\).

We introduce a continuous scale of Hilbert spaces
\begin{align}\label{eq:def_Yl}
\mathcal{Y}^{\ell} = D(\Lambda^{\ell}), \qquad \ell \geq 0,
\end{align}
with \(\mathcal{Y}^{0} = H=L^2(\Omega)\). From the definition of \(A\), it is straightforward to verify that \(\mathcal{Y}^{1} = H_0^{1}(\Omega)\). 
For negative indices $\ell$, the spaces are defined by duality,
$
\mathcal{Y}^{\ell} = \big(\mathcal{Y}^{-\ell}\big)^{*}.
$
%
For \(\ell \geq 0\), the inner product on \(\mathcal{Y}^{\ell}\) is given by
\begin{align}\label{eq:Y-ell}
	\langle u_1, u_2 \rangle_{\ell} 
	= \langle u_1, u_2 \rangle_{H} 
	+ \langle \Lambda^{\ell} u_1, \Lambda^{\ell} u_2 \rangle_{H}.
\end{align}
Accordingly, each \(u \in \mathcal{Y}^{\ell} \hookrightarrow H\) admits the expansion
\begin{align}\label{eq:expansion_u}
u = \sum_{\alpha \in I} \widehat{u}_{\alpha} e_{\alpha},
\end{align}
and its norm in \(\mathcal{Y}^{\ell}\) satisfies
 \begin{align}\label{eq:Yl_norm}
		\|u\|_{\mathcal{Y}^{\ell}}
		= \sqrt{\langle u, u \rangle_{\ell}}
		\;\sim\;
		\left( \sum_{\alpha \in I} (1 + \lambda_{\alpha}^{2\ell})
		|\widehat{u}_{\alpha}|^{2} \right)^{\frac{1}{2}}.
	\end{align}
	
	%
	By standard interpolation and duality arguments, one obtains
	\begin{align}\label{eq:relation_Y_H}
		\mathcal Y^\ell =
		\begin{cases}
			H^\ell(\Omega), & 0\le \ell<\frac12,\\
			H_{00}^{1/2}(\Omega), & \ell=\frac12,\\
			H_0^\ell(\Omega), & \frac12<\ell\le1,
		\end{cases}
	\end{align}
	with equivalence of norms. Consequently,
	\begin{align}\label{eq:relation_Yl}
		\mathcal{Y}^{\ell} \hookrightarrow H^{\ell}(\Omega)
		\quad \text{and} \quad
		H^{-\ell}(\Omega) \hookrightarrow \mathcal{Y}^{-\ell},
		\qquad 0\le \ell\le 1.
\end{align}

\subsection{Fully discrete scheme and main results}\label{subsec:main}
To derive a numerical approximation for the second-order nonlinear wave equation, we first reformulate it as a first-order system:
\begin{equation}\label{system}
	\left\{
	\begin{array}{lll}
		\partial_t U - L U = F(U) & \text{for } t \in (0, T], \\[2mm]
		U|_{t=0} = U^0,
	\end{array}
	\right.
\end{equation}
where
\begin{align}\label{system_notation}
	U = \begin{pmatrix}
		u \\ v
	\end{pmatrix}, \quad
	U^0 = \begin{pmatrix}
		u^0 \\ v^0
	\end{pmatrix}, \quad
	F(U) = \begin{pmatrix}
		0 \\ f(u)
	\end{pmatrix}, \quad
	L = \begin{pmatrix}
		0 & 1 \\ -A & 0
	\end{pmatrix}.
\end{align}

For each $U=(u,v)^{\top}$, we define
\[
\|U\|_\ell := \|U\|_{\mathcal{Y}^\ell \times \mathcal{Y}^{\ell-1}} = \left( \|u\|_{\mathcal{Y}^\ell}^2 + \|v\|_{\mathcal{Y}^{\ell-1}}^2 \right)^{\frac{1}{2}}.
\]

As \(A\) is self-adjoint, the associated operator \(L\) is skew-adjoint, ensuring that all eigenvalues of \(L\) are purely imaginary. Consequently, the semigroup \(e^{tL}\) forms a \(C^0\)-semigroup, satisfying the key stability property
\begin{align}\label{semigroup_prop_etL}
	\|e^{tL}\|_{\mathcal{L}(\mathcal{Y}^{\ell}\times\mathcal{Y}^{\ell-1})}\lesssim 1.
\end{align}

For the spatial discretization, let \((X_h^k)_h\) be a family of conforming Lagrange finite element spaces of degree \(k\ge1\), defined on a family of simplicial meshes \((\mathcal{T}_h)_h\), where \(h>0\) denotes the mesh size. We assume that the finite element spaces \(X_h^k\) satisfy the following approximation property for functions in \(\mathcal{Y}^\ell\). 

\begin{assumption}\label{assump:mesh}
	Let \(\mathcal{T}_h\) be a simplicial partition of the domain \(\Omega\), and let \(X_h^k\) (\(k\ge 1\)) be a family of finite-dimensional subspaces satisfying \(\dim(X_h^k)=O(h^{-d})\). Assume that for any \(0\le \ell \le k\) and any \(w\in \mathcal{Y}^{\ell+1}\), the following approximation properties hold:
	\begin{align}
		\inf_{w_h\in X_h^k}\|w-w_h\|_{L^2(\Omega)} &\le C h^{\ell+1}\|w\|_{\mathcal{Y}^{\ell+1}}, \label{eq:L2-approx}\\
		\inf_{w_h\in X_h^k}\|w-w_h\|_{H^1(\Omega)} &\le C h^{\ell}\|w\|_{\mathcal{Y}^{\ell+1}}, \label{eq:H1-approx}
	\end{align}
	where \(C\) is a constant independent of \(h\), but possibly dependent on the domain \(\Omega\).
\end{assumption}

\begin{remark}
	\label{rem:mesh_ass}\upshape
	Typical situations in which Assumption~\ref{assump:mesh} is satisfied are as follows.
	\begin{itemize}
		\item If \(\Omega=[0,1]^d\) with \(d=2,3\), then by the definition of the space \(\mathcal{Y}^{\ell}\) in \eqref{eq:def_Yl}, we have \(\mathcal{Y}^{\ell}\hookrightarrow H^{\ell}(\Omega)\) for all \(\ell\ge 0\). Therefore, Assumption~\ref{assump:mesh} holds for all \(k\ge 1\) on quasi-uniform, shape-regular simplicial partitions, by the standard approximation theory for finite element spaces.
		
		\item If \(\Omega\subset\mathbb{R}^2\) is a polygonal domain, then the eigenfunctions of \(\Lambda=\sqrt{A}\) generally exhibit corner singularities and are not smooth near the vertices. A detailed singularity analysis, based on Kondrat'ev theory, yields
		\begin{align}\label{eq:embedded_Yl}
			\mathcal{Y}^{\ell_{\Omega}} \hookrightarrow H^{\ell_{\Omega}}(\Omega),
			\qquad
			\ell_{\Omega}=1+\frac{\pi}{\theta}-\varepsilon,
		\end{align}
		for arbitrarily small \(\varepsilon>0\), where \(\theta=\max_{j=1,\ldots,m}\theta_j\) is the largest interior angle of \(\Omega\); see \cite[Theorem~9.8 and Remark~9.12]{bourlard1992coefficients}. Consequently, Assumption~\ref{assump:mesh} holds for quasi-uniform, shape-regular simplicial partitions \(\mathcal{T}_h\), provided that the polynomial degree satisfies \(1\le k \le \ell_{\Omega}-1\).
		
		\item More generally, if \(\Omega\) is a polygonal domain, then Kondrat'ev theory shows that the eigenfunctions of \(\Lambda=\sqrt{A}\) belong to appropriate weighted Sobolev (Kondrat'ev) spaces that capture the corner singularities \cite{dauge2006elliptic,kondrat1967boundary}. In this setting, one can construct suitably graded simplicial meshes \(\mathcal{T}_h\) such that Assumption~\ref{assump:mesh} continues to hold; see, for example, \cite{buacuctua2005improving,li2022maximum,li2017W1} for graded-mesh constructions for functions in Kondrat'ev spaces.
	\end{itemize}
	\end{remark}

We define the projection operator \(\P_h : H \to X_h^k\) by
\begin{align}\label{define_Ph}
	\langle \P_h u, v_h \rangle_H = \langle u, v_h \rangle_H \quad \text{for all } v_h \in X^k_h.
\end{align}

Furthermore, let \(A_h : X_h^k \to X_h^k\) denote the finite element discretization of \(A\), 
and define the Ritz projection \(\mathcal{R}_h : V \to X_h^k\) through
\begin{subequations} \label{define_Rh}
	\begin{align}
		a(u_h, v_h ) = 		\langle A_h u_h, v_h \rangle_H
		\quad
		\text{and}
		\quad
		a(\R_h u, v_h ) = a( u, v_h ),
	\end{align}
	It follows that for all \(u \in D(A)\),
	\begin{align}
		\langle A_h (\R_h u), v_h \rangle_H = \langle A u, v_h \rangle_H .
	\end{align}
\end{subequations}
For a time step size \(\tau>0\), let \(t_n=n\tau\). By the variation-of-constants formula applied to \eqref{system}, the exact solution satisfies
\begin{align}\label{mild_form}
	U(t_{n+1})
	= e^{\tau L}U(t_n)
	+ \int_0^\tau e^{(\tau-\sigma)L} F\bigl(U(t_n+\sigma)\bigr)\,\mathrm d\sigma .
\end{align}
Applying the spatial discretization operators \(\P_h\) and \(A_h\), we obtain the spatially semidiscrete system
\begin{align}\label{eq:spacial_semidiscrete}
	U_h(t_{n+1})
	= e^{\tau L_h}U_h(t_n)
	+ \int_0^\tau e^{(\tau-\sigma)L_h}\P_h F\bigl(U_h(t_n+\sigma)\bigr)\,\mathrm d\sigma ,
\end{align}
where, for notational simplicity, we use the same symbol \(\P_h\) for its diagonal extension \(\operatorname{diag}(\P_h,\P_h)\), acting componentwise on vector-valued functions. The initial data are given by \(U_h(0)=\P_h U(0)\), and the discrete operator \(L_h\) is defined by
\begin{align}\label{def_Lh}
	L_h = \begin{pmatrix}
		0 & 1 \\ -A_h & 0
	\end{pmatrix}.
\end{align}

To construct a fully discrete scheme, we next discretize the integral term on the right-hand side of~\eqref{eq:spacial_semidiscrete}. 
Approximating the nonlinear term by
\[\P_h F\big(U_h(t_n+\sigma)\big)\approx \P_h F\big(U_h(t_n)\big)\] we arrive at the fully discrete exponential Euler scheme:

%

\begin{align}\label{eq:fully-discrete}
	U^{n+1}_h = e^{\tau L_h} U^n_h + \tau \varphi_1(\tau L_h) \P_h F(U^n_h), \quad \text{with}\quad \varphi_1(z)=\frac{e^z-1}{z},
\end{align}
where $U^0_h = \P_h U(0)$.

Throughout this article, we impose the following assumptions:
\begin{assumption}\label{assumption3}
	\
	
	\begin{itemize}
		%
		\item[(a)] The projection operator \(\P_h\) is stable in \(\mathcal{Y}^{\ell}\) for \(\ell=-1,0,1\), i.e.,
		\begin{align}\label{eq:projection_stability}
			\|\P_h u\|_{\mathcal{Y}^{\ell}} \leq C_P \|u\|_{\mathcal{Y}^{\ell}},
		\end{align}
		for some constant \(C_P>0\).
		\item[(b)] The linear operator \(L_h\) generates a \(C_0\)-semigroup, and the following bound holds for \(\ell=0,1\):
		\begin{align}\label{assump_etLh}
			\|e^{t L_h} U\|_{\mathcal{Y}^{\ell} \times \mathcal{Y}^{\ell-1}}
			\leq e^{C_L t} \|U\|_{\mathcal{Y}^{\ell} \times \mathcal{Y}^{\ell-1}},
		\end{align}
		for some constant \(C_L > 0\) and for all \(U \in X^k_h \times X^k_h\).
		
		\item[(c)] For $d=2,3$ and $\gamma\in(0,1]$, the nonlinear function \(f\) is differentiable and satisfies the following growth conditions
		\begin{align}\label{eq:growth_condition}
			|f(u)|\leq C_f (1+|u|)^{\xi_{\gamma}},\quad |f^{\prime}(u)|\leq C_f (1+|u|)^{\xi_{\gamma}-1},
		\end{align}
		with 
		\begin{align}\label{eq:xi_gamma}
			\xi_{\gamma}<\left\{\begin{array}{l}
				\infty,\quad \qquad\quad\text{for }d=2,\; \gamma=1,\\[2mm]
				1+\frac{2}{d-2\gamma},\quad\;\text{for } d=2,\ 0<\gamma<1,\ \text{or}\quad d=3,\ 0<\gamma\le 1.
			\end{array}\right.
		\end{align}
	\end{itemize}
\end{assumption}
Assumption~\ref{assumption3}(a), in particular the boundedness in \(\mathcal{Y}^1=H_0^1(\Omega)\), is not restrictive, as it holds for quasi-uniform triangulations and certain adaptively refined meshes; see \cite{bank2014h,diening2021sobolev}.
 By a standard duality argument, this immediately yields the corresponding boundedness on \(\mathcal{Y}^{-1}\).
Assumption~\ref{assumption3}(b) is standard for spatial discretizations for which \(A_h\) is symmetric and positive definite.
The growth condition in Assumption~\ref{assumption3}(c) is imposed in order to derive the stability estimates for the numerical scheme in the subsequent analysis.

\begin{remark}\upshape\label{rem:assumption_corollary}
	By interpolation, Assumptions~\ref{assumption3}(a) and~(b) imply that \(\P_h\) is uniformly bounded on \(\mathcal{Y}^{\theta}\) for all \(\theta\in[-1,1]\), and that \(e^{tL_h}\) is uniformly bounded on \(\mathcal{Y}^{\gamma}\times\mathcal{Y}^{\gamma-1}\) for all \(\gamma\in[0,1]\).
	
	Moreover, the nonlinear mapping \(F\) is locally Lipschitz continuous on \(\mathcal{Y}^{\gamma}\times\mathcal{Y}^{\gamma-1}\) for all $\gamma\in(0,1]$, under Assumptions~\ref{assumption3}(c). Indeed, for any \(U_i=(u_i,v_i)^{\top}\in \mathcal{Y}^{\gamma}\times\mathcal{Y}^{\gamma-1}\), \(i=1,2\), Lemma~\ref{lem:nonlinear_term_estimates}, in particular \eqref{eq:fu-fv_gamma}, yields
	\begin{align*}
		\|F(U_1)-F(U_2)\|_{\gamma}
		= \|f(u_1)-f(u_2)\|_{\mathcal{Y}^{\gamma-1}}
		\lesssim C(\|U_1\|_{\gamma},\|U_2\|_{\gamma})\|U_1-U_2\|_{\gamma}.
	\end{align*}
	Therefore, by the boundedness of \(e^{tL}\) and a standard fixed-point argument, problem~\eqref{system} is locally well posed in \(C([0,T];\mathcal{Y}^{\gamma}\times\mathcal{Y}^{\gamma-1})\) for some \(T>0\).
\end{remark}
%
%

We summarize the main convergence result as follows:
\begin{theorem}[Main result]\label{thm:main}
	Under Assumptions~\ref{assump:A}--\ref{assumption3}, suppose that the exact solution of \eqref{system} satisfies \(U \in C([0,T]; \mathcal{Y}^{\gamma} \times \mathcal{Y}^{\gamma-1})\) for some \(\gamma\in(0,1]\). Then there exists \(h_0\in(0,1)\) such that, for all \(0<h<h_0\) and $\tau>0$ satisfying the step-size condition \(\tau\lesssim h^{\frac{\ell_k}{\ell_k+1}(1-\gamma)+\mu}\) for an arbitrarily small but fixed \(\mu>0\), the fully discrete scheme~\eqref{eq:fully-discrete} satisfies
	\begin{equation}\label{eq:main}
		\sup_{0\le n\le T/\tau}\|U(t_n)-U_h^n\|_0
		\le
		C\bigl(\tau\, h^{\frac{\ell_k}{\ell_k+1}(\gamma-1)} + h^{\frac{\ell_k}{\ell_k+1}\gamma}\bigr),
	\end{equation}
	where
	\begin{align}\label{eq:def_lk}
		\ell_k :=
		\begin{cases}
			k+1, & k=1,\\
			k+2, & k\ge2,
		\end{cases}
	\end{align}
	with \(k\ge1\) denoting the polynomial degree of the finite element space. The constant \(C>0\) depends only on \(T\), \(\sup_{t\in[0,T]}\|U(t)\|_\gamma\), and the constants in Assumptions~\ref{assump:A}--\ref{assumption3}, but is independent of \(h\) and \(\tau\).
	
	Moreover, choosing \(\tau=O(h^{\frac{\ell_k}{\ell_k+1}})\), which balances the two terms on the right-hand side of \eqref{eq:main} and is independent of \(\gamma\), yields the total error estimate
	\begin{equation}\label{eq:main-1}
		\sup_{0\le n\le T/\tau}\|U(t_n)-U_h^n\|_0
		\le
		C\,\tau^\gamma,
	\end{equation}
	for solutions \(U \in C([0,T]; \mathcal{Y}^{\gamma} \times \mathcal{Y}^{\gamma-1})\).
\end{theorem}

\begin{remark}\upshape\label{rem:CFL}
	The first term on the right-hand side of \eqref{eq:main} represents the temporal discretization error, while the second term is the spatial discretization error. Owing to the limited regularity of the exact solution, the exponential integrator does not attain first-order convergence in time under the present assumptions. This loss is reflected in the additional factor \(h^{\frac{\ell_k}{\ell_k+1}(\gamma-1)}\) appearing in the temporal error term. In the borderline case \(\gamma=1\), this factor reduces to \(1\), and the optimal first-order temporal convergence is recovered.
	
	The step-size condition \(\tau\lesssim h^{\frac{\ell_k}{\ell_k+1}(1-\gamma)+\mu}\) is imposed to guarantee the smallness requirement in the bootstrap argument of Lemma~\ref{lem:discrete_bootstrap}, in particular, to ensure that \eqref{eq:tau0h0} holds. Moreover, under this condition, the error bound \eqref{eq:main} implies
	\[
	\sup_{0\le n\le T/\tau}\|U(t_n)-U_h^n\|_0
	\le
	C\bigl(h^{\mu}+h^{\frac{\ell_k}{\ell_k+1}\gamma}\bigr)\to 0
	\qquad\text{as } h,\tau\to0.
	\]
	Finally, the choice \(\tau=O(h^{\frac{\ell_k}{\ell_k+1}})\) balances the temporal and spatial contributions in \eqref{eq:main}, automatically satisfies the above step-size restriction, and yields the overall error estimate \(O(\tau^\gamma)\).
\end{remark}

\begin{remark}\upshape
	We emphasize that the spatial convergence rates differ between the cases \(k=1\) and \(k\ge2\). This distinction is intrinsic to our weak-norm error analysis. More precisely, in the \(\mathcal{Y}^{-1}\) norm the Ritz projection admits a higher-order estimate for \(k\ge2\); see \eqref{eq:estimate_Rh_Hminus1} in Lemma~\ref{lem:PhRh} and also \cite[Section~5.8]{brenner2008mathematical}. By contrast, for \(k=1\) no analogous estimate is available in this form. Consequently, the cases \(k=1\) and \(k\ge2\) must be treated separately in the error analysis; see Lemma~\ref{lem:estimate_etL-etLh}.
	
	Owing to the refined weak-norm estimate established in Lemma~\ref{lem:estimate_etL-etLh}, the spatial convergence rate appearing on the right-hand side of \eqref{eq:main}, although somewhat involved in form, appears to be sharp. It is also worth noting that this rate is only of order \(h^{\frac{\ell_k}{\ell_k+1}\gamma}\), due to the limited regularity of the exact solution, so that the optimal rate for \(k\)-th order finite elements cannot in general be expected. Nevertheless, our estimate shows that higher-order methods still yield a clear improvement in spatial convergence over lower-order ones, even in this low-regularity regime. These observations are confirmed by the numerical experiments in Section~\ref{sec:numerical_experiments}.
\end{remark}
	
%

\subsection{Frequency-localized auxiliary scheme}
Since the classical error analysis of the fully discrete scheme~\eqref{eq:fully-discrete} requires higher regularity of the exact solution, it is not directly applicable in our low-regularity setting. To overcome these difficulties, we introduce a frequency-localized auxiliary scheme, which serves as an analytical bridge in the error analysis.

For \(t_n=n\tau\), the variation-of-constants formula gives
\begin{align}\label{eq:construct_semi_discrete_1}
	U(t_{n+1})=e^{\tau L}U(t_n)+\int_{0}^{\tau}e^{(\tau-s)L}F\big(e^{sL}U(t_n)\big)\,\d s+R_1(t_n),
\end{align}
where
\begin{align}\label{eq:R_1}
	R_1(t_n)
	=\int_{0}^{\tau}e^{(\tau-s)L}\Big(F\big(U(t_n+s)\big)-F\big(e^{sL}U(t_n)\big)\Big)\,\d s.
\end{align}
The remainder term \(R_1(t_n)\) can be controlled by combining the semigroup property of \(e^{tL}\) in \eqref{semigroup_prop_etL} with the variation-of-constants formula; see Lemma~\ref{lem:estimate_R1R2R3R4} below.

We next consider the approximation of the remaining integral term in \eqref{eq:construct_semi_discrete_1}. A standard first-order approximation would replace \(F(e^{sL}U(t_n))\) by \(F(U(t_n))\), as in the classical exponential integrator. However,
\begin{align}\label{eq:EI_error}
	F\big(e^{sL}U(t_n)\big)-F\big(U(t_n)\big)
	=
	\int_0^s \frac{\d}{\d \sigma}F\big(e^{\sigma L}U(t_n)\big)\,\d\sigma,
\end{align}
and the derivative on the right-hand side involves spatial derivatives of \(U(t_n)\). This presents a difficulty under the low-regularity assumption 
\((u,v) \in \mathcal{Y}^{\gamma} \times \mathcal{Y}^{\gamma-1}\) with \(\gamma \in (0,1]\), 
since the term \(\tfrac{\mathrm{d}}{\mathrm{d}\sigma}F(e^{\sigma L}U(t_n))\) requires one additional order of spatial differentiability beyond what is available.

Inspired by \cite{CLLY2024}, 
we introduce a frequency decomposition operator $\Pi_r$, 
defined as the spectral projection of functions onto the low-frequency subspace $H_r \hookrightarrow H$, 
which is spanned by the basis functions
\[
\mathcal{B}_r = \left\{ e_\alpha \in H \;\middle|\; \alpha \in I,\; |\lambda_\alpha| \leq r \right\}.
\]
The projection operator $\Pi_r$ satisfies
\begin{align}\label{eq:projection}
	\Pi_r u \in H_r, 
	\qquad 
	\langle \Pi_r u, v \rangle_H = \langle u, v \rangle_H, 
	\quad \forall\, v \in H_r.
\end{align}
In our subsequent error analysis, we make essential use of the improved regularity of the low-frequency projected component to establish rigorous error bounds.

Thus we approximate the integral term on the right-hand side of \eqref{eq:construct_semi_discrete_1} by
\begin{align}\label{eq:approximate_F}
	\int_{0}^{\tau}e^{(\tau-s)L}F\big(e^{sL}U(t_n)\big)&=\int_{0}^{\tau}e^{(\tau-s)L}F\big(e^{sL}\Pi_r U(t_n)\big)\d s+R_2(t_n)\notag\\
	&=\tau \varphi_1(\tau L)F\big(\Pi_r U(t_n)\big)+R_2(t_n)+R_3(t_n)\notag\\
	&=\tau \varphi_1(\tau L)\Pi_r F(\Pi_r U(t_n))+R_2(t_n)+R_3(t_n)+R_4(t_n),
\end{align}
where \(\varphi_1\) is given by \eqref{eq:fully-discrete}.
The residual terms
\(R_2(t_n)\), \(R_3(t_n)\), and \(R_4(t_n)\) are defined by:
\begin{subequations} \label{eq:R2R3R4}
\begin{align}
	R_2(t_n)=&\int_{0}^{\tau}e^{(\tau-s)L}\left(F(e^{s L}U(t_n))-F(e^{s L}\Pi_r U(t_n))\right)\d s,\\
	R_3(t_n)=&\int_0^\tau e^{(\tau-s) L}\int_0^s \frac{\d}{\d \sigma}F(e^{\sigma L}\Pi_r U(t_n))\d\sigma\d s,\label{eq:R3}\\
	R_4(t_n)=&\tau \varphi_1(\tau L)\left(F\big(\Pi_r U(t_n)\big)-\Pi_r F(\Pi_r U(t_n))\right).
\end{align}
\end{subequations}
The approximation \eqref{eq:approximate_F} together with \eqref{eq:construct_semi_discrete_1} implies
\begin{align}\label{eq:consistency_error}
	U(t_{n+1})=e^{\tau L}U(t_n)+\tau \varphi_1(\tau L)\Pi_r F(\Pi_r U(t_n))+\sum_{i=1}^4 R_i(t_n).
\end{align}
And thus we arrive at the auxiliary scheme
\begin{align}\label{eq:semi-discrete}
	U^{n+1} &= e^{\tau L} U^n + \tau \varphi_1(\tau L) \Pi_r F(\Pi_r U^n),
	\quad n \geq 0.
\end{align}

%
%
\begin{remark}\upshape
	The scheme \eqref{eq:semi-discrete} differs from the classical exponential integrator in that the frequency projection \(\Pi_r\) is applied both to the numerical solution and to the nonlinearity. This provides the additional regularity needed for the consistency analysis under low-regularity initial data. We emphasize, however, that \eqref{eq:semi-discrete} is introduced only as an analytical tool: in general, the spectral projector \(\Pi_r\) is not computationally practical on general domains, and our fully discrete method is therefore not obtained by directly discretizing \eqref{eq:semi-discrete}.
\end{remark}

\section{Preliminary results}\label{sec:preliminary}
In this section, we present the key lemmas that will be used to prove the convergence result. These include the Bernstein estimates for the operator \(\Pi_r\) defined in \eqref{eq:projection} and the projection error estimates for the operator \(\P_h\)  defined in \eqref{define_Ph}.
Finally, we compare the semigroup \(e^{tL}\) of the continuous problem with the semigroup \(e^{tL_h}\) obtained after finite element spatial discretization. This comparison plays a crucial role in establishing the error estimates for the fully discrete scheme of the exponential integrator based on finite element discretization.
\begin{lemma}[Bernstein-type inequality]\label{lem:Bernstein}
	Let $u\in \mathcal{Y}^{\ell}$ for $\ell\in \mathbb{R}$, then there holds
	\begin{align}\label{eq:Bernstein}
		\|\Pi_r u\|_{\mathcal{Y}^{\ell+m}}\lesssim r^{m}\|u\|_{\mathcal{Y}^{\ell}}\quad\text{and}\quad \|\Pi_{>r} u\|_{\mathcal{Y}^{\ell-m}}\lesssim r^{-m}\|u\|_{\mathcal{Y}^{\ell}},
	\end{align}
	for all $m\geq 0$.
\end{lemma}

\begin{proof}
	The result follows directly from the definition of the norms in \(\mathcal Y^\ell\) and the spectral localization properties of \(\Pi_r\) and \(\Pi_{>r}\).
\end{proof}

Next, we use a duality argument to derive an estimate for the $L^2$-projection in $\mathcal{Y}^{-1}$.

\begin{lemma}[Projection error estimate]\label{lem:projection}
	Let \(\P_h\) be the projection operator defined in \eqref{define_Ph}. Then, for every \(\ell\in[0,1]\),
	\begin{align}
		\|(\mathrm{id}-\P_h)u\|_{\mathcal{Y}^{0}}
		&\lesssim h^{\ell} \|u\|_{\mathcal{Y}^{\ell}},
		\qquad\quad u\in \mathcal{Y}^{\ell},\label{eq:lem_projection_1}\\
		\|(\mathrm{id}-\P_h)v\|_{\mathcal{Y}^{-1}}
		&\lesssim h^{\ell} \|v\|_{\mathcal{Y}^{\ell-1}},
		\qquad \,  v\in \mathcal{Y}^{\ell-1}.\label{eq:lem_projection}
	\end{align}
\end{lemma}

\begin{proof}
	We first prove \eqref{eq:lem_projection_1}. Since \(\mathcal Y^0=L^2(\Omega)\) and \(\mathcal Y^1=H_0^1(\Omega)\), it follows from the standard approximation property of the \(L^2\)-projection and Assumption~\ref{assump:mesh} that
	\begin{align*}
		\|(\mathrm{id}-\P_h)u\|_{\mathcal Y^0}=\|(\mathrm{id}-\P_h)u\|_{L^2(\Omega)}\lesssim h\|u\|_{H_0^1(\Omega)}=h\|u\|_{\mathcal Y^1}
	\end{align*}
	for all \(u\in\mathcal Y^1\). On the other hand, \(\P_h\) is bounded on \(\mathcal Y^0=L^2(\Omega)\), and hence \(\|(\mathrm{id}-\P_h)u\|_{\mathcal Y^0}\lesssim \|u\|_{\mathcal Y^0}\) for all \(u\in\mathcal Y^0\). Interpolating between the above two estimates, we obtain \eqref{eq:lem_projection_1} for all \(u\in\mathcal Y^\ell\), \(\ell\in[0,1]\).
	
	Next we prove \eqref{eq:lem_projection}. By duality, for all \(v\in\mathcal Y^0\),
	\begin{align*}
		\|(\mathrm{id}-\P_h)v\|_{\mathcal Y^{-1}}
		&=\sup_{\|w\|_{\mathcal Y^1}=1}\langle (\mathrm{id}-\P_h)v,w\rangle_H
		=\sup_{\|w\|_{\mathcal Y^1}=1}\langle v,(\mathrm{id}-\P_h)w\rangle_H\\
		&\leq \|v\|_{\mathcal Y^0}\sup_{\|w\|_{\mathcal Y^1}=1}\|(\mathrm{id}-\P_h)w\|_{\mathcal Y^0}
		\lesssim h\|v\|_{\mathcal Y^0},
	\end{align*}
	where we used \eqref{eq:lem_projection_1} with \(\ell=1\) in the last step. Moreover, by Assumption~\ref{assumption3}(a), \(\P_h\) is bounded on \(\mathcal Y^{-1}\), and hence \(\|(\mathrm{id}-\P_h)v\|_{\mathcal Y^{-1}}\lesssim \|v\|_{\mathcal Y^{-1}}\) for all \(v\in\mathcal Y^{-1}\). Interpolating between these two estimates, we obtain \eqref{eq:lem_projection} for all \(v\in\mathcal Y^{\ell-1}\), \(\ell\in[0,1]\).
\end{proof}

Furthermore, under Assumption~\ref{assump:A}, by combining the duality argument with the Aubin--Nitsche technique, we obtain the following high-order approximation estimates for the projection operators \(\P_h\) and \(\R_h\).
\begin{lemma}[Error estimates of projection operators $\P_h$ and $\R_h$]\label{lem:PhRh}
	Under Assumption~\ref{assump:mesh}, let $\P_h$ and $\R_h$ be the projection operators onto the finite element space $X_h^k$, defined by \eqref{define_Ph} and \eqref{define_Rh}. Then for any $u\in \mathcal{Y}^{k+1}$, the following estimates hold for $k\geq 1$:
	\begin{align}\label{eq:estimate_PhRh_L2}
		\|u-\P_h u\|_{\mathcal{Y}^0}+\|u-\R_h u\|_{\mathcal{Y}^0}\lesssim h^{k+1}\|u\|_{\mathcal Y^{k+1}},
	\end{align}
	and
	\begin{align}\label{eq:estimate_Ph_Hminus1}
		\|u-\P_h u\|_{\mathcal{Y}^{-1}}\lesssim h^{k+2}\|u\|_{\mathcal Y^{k+1}}.
	\end{align}
	Moreover, if $k\geq 2$, there holds
	\begin{align}\label{eq:estimate_Rh_Hminus1}
		\|u-\R_h u\|_{\mathcal{Y}^{-1}}\lesssim h^{k+2}\|u\|_{\mathcal Y^{k+1}}.
	\end{align}
\end{lemma}

\begin{proof}
	Since $\mathcal{Y}^0=L^2(\Omega)$, the estimate \eqref{eq:estimate_PhRh_L2} follows directly from the approximation property in Assumption~\ref{assump:mesh}, the best approximation property of $\P_h$ in $L^2(\Omega)$, and the best approximation property of $\R_h$ in $H_0^1(\Omega)$ together with the standard Aubin--Nitsche argument.
	
	The estimate \eqref{eq:estimate_Ph_Hminus1} is easily verified by duality and by using the error estimate \eqref{eq:estimate_PhRh_L2}:
	\begin{align*}
		\|u-\P_h u\|_{\mathcal{Y}^{-1}}
		&=\sup_{\|w\|_{\mathcal{Y}^1}=1}\langle u-\P_h u, w\rangle_H
		=\sup_{\|w\|_{\mathcal{Y}^1}=1}\langle u-\P_h u, w-\P_h w\rangle_H\\
		&\leq \sup_{\|w\|_{\mathcal{Y}^1}=1}\|u-\P_h u\|_{\mathcal{Y}^0}\|w-\P_h w\|_{\mathcal{Y}^0}
		\lesssim h^{k+2}\|u\|_{\mathcal Y^{k+1}}.
	\end{align*}
	
	For \eqref{eq:estimate_Rh_Hminus1}, let $\phi=A^{-1}w$. Then
	\begin{align}\label{eq:estimate_u-Rhu}
		\|u-\R_h u\|_{\mathcal{Y}^{-1}}
		&=\sup_{\|w\|_{\mathcal{Y}^1}=1}\langle u-\R_h u, w\rangle_H
		=\sup_{\|w\|_{\mathcal{Y}^1}=1}\langle u-\R_h u, A\phi\rangle_H\notag\\
		&=\sup_{\|w\|_{\mathcal{Y}^1}=1}a \bigl(u-\R_h u, \phi\bigr)
		=\sup_{\|w\|_{\mathcal{Y}^1}=1}a \bigl(u-\R_h u, \phi-\R_h \phi\bigr),
	\end{align}
	where we have used the orthogonality
	\begin{align*}
		a \bigl(u-\R_h u, \R_h \phi\bigr)=0.
	\end{align*}
	By using the elliptic estimate $\|\phi\|_{\mathcal{Y}^3}\lesssim\|w\|_{\mathcal{Y}^1}$, and the best approximation property of $\R_h$ for $u\in \mathcal{Y}^{k+1}$ and $\phi\in \mathcal{Y}^{3}$ in the $\mathcal{Y}^1$ norm, together with Assumption~\ref{assump:mesh}, we have
	\begin{align*}
		a \bigl(u-\R_h u, \phi-\R_h \phi\bigr)
		\lesssim \|u-\R_h u\|_{\mathcal{Y}^1}\|\phi-\R_h \phi\|_{\mathcal{Y}^1}
		\lesssim h^{k}\|u\|_{\mathcal{Y}^{k+1}}h^2\|\phi\|_{\mathcal{Y}^3}
		\lesssim h^{k+2}\|u\|_{\mathcal{Y}^{k+1}}\|w\|_{\mathcal{Y}^1}.
	\end{align*}
	By combining this with \eqref{eq:estimate_u-Rhu} we conclude \eqref{eq:estimate_Rh_Hminus1}.
\end{proof}

We are now in a position to estimate the difference between the two semigroups.
\begin{lemma}\label{lem:estimate_etL-etLh}
	Let \(V(0) \in \mathcal{Y}^\gamma \times \mathcal{Y}^{\gamma-1}\), and let the operators \(L\) and \(L_h\) be as defined in \eqref{system_notation} and \eqref{def_Lh}, respectively. Then, under Assumptions~\ref{assump:A}, \ref{assump:mesh}, and~\ref{assumption3}, the following estimate holds for any fixed \(T_*>0\), \(1\le r\le h^{-1}\), and \(t\in[0,T_*]\):
	\begin{align}\label{eq:improved_etL-etLh_Lipschitz}
		\left\|e^{t L} \Pi_r V(0) - e^{t L_h} \P_h \Pi_r V(0)\right\|_0 
		\lesssim \left\{\begin{array}{l}
			h^{2} r^{3-\gamma} \,\|V(0)\|_{\mathcal{Y}^\gamma \times \mathcal{Y}^{\gamma-1}},\quad\quad\qquad\qquad \text{for}\quad k=1,\\[2mm]
			\big(h^{k+2} r^{k+3-\gamma} + h^{\gamma}\big)\|V(0)\|_{\mathcal{Y}^\gamma \times \mathcal{Y}^{\gamma-1}},\quad \text{for}\quad k\geq 2,
		\end{array}\right.
	\end{align}
	where $k$ is the polynomial degree of the finite element basis functions of \(X^k_h\). The hidden constant may depend on \(T_*\), but is independent of \(h\) and \(r\).
\end{lemma}

\begin{proof}[Proof (sketch)]
	This result is equivalent to estimating the error of the finite element semidiscretization for the linear wave equation with initial data \(\Pi_r V(0)\) in the weak norm \(\mathcal{Y}^0 \times \mathcal{Y}^{-1}\). Owing to the presence of the frequency localization operator \(\Pi_r\), the initial data are effectively regularized. In particular, by the Bernstein-type inequality in Lemma~\ref{lem:Bernstein}, arbitrary powers of the differential operator \(A\) acting on the low-frequency part remain bounded, with bounds depending only on suitable positive powers of \(r\).
	
	This observation allows us to combine the regularizing property of \(\Pi_r\) with the approximation estimates in Lemma~\ref{lem:PhRh}. We note, however, that the proof differs essentially between the cases \(k\ge2\) and \(k=1\). For \(k\ge2\), it relies on the sharp \(\mathcal{Y}^{-1}\)-error estimate for the Ritz projection. For \(k=1\), since such an estimate is unavailable, we instead use the sharp \(\mathcal{Y}^0\)-error estimate for the Ritz projection and a different energy-based argument.
	
	The detailed proof is deferred to Appendix~\ref{sec:proof_or_etL-etLh}.
\end{proof}

We next collect several estimates for the nonlinear term in \eqref{eq:NLW}:
\begin{lemma}\label{lem:nonlinear_term_estimates}
	Suppose that the nonlinear function $f$ satisfies the growth conditions \eqref{eq:growth_condition}. Then, the following estimates hold:
	\begin{align}
		\|f(u)\|_{\mathcal{Y}^{\gamma-1}}&\leq C\big(\|u\|_{\mathcal{Y}^{\gamma}}\big),\label{eq:f_gamma-1}\\[2mm]
		\|f(u)-f(v)\|_{\mathcal{Y}^{\gamma-1}}& \leq C\big(\|u\|_{\mathcal{Y}^{\gamma}},\|v\|_{\mathcal{Y}^{\gamma}}\big)\|u-v\|_{\mathcal{Y}^{\gamma}},\label{eq:fu-fv_gamma}\\[2mm]
		\|f(u)\|_{\mathcal{Y}^{-1}}&\leq C\big(\|u\|_{\mathcal{Y}^{\gamma}}\big)\big(1+\|u\|_{\mathcal{Y}^{0}}\big),\label{eq:fu_H-1}\\[2mm]
		\|f(u)-f(v)\|_{\mathcal{Y}^{-1}}&\leq C\Big(\|u\|_{\mathcal{Y}^{\gamma}},\|v\|_{\mathcal{Y}^{\gamma}} \Big)\|u-v\|_{\mathcal{Y}^{0}},\label{eq:fu-fv}\\[2mm]
		\|f^{\prime}(u)v\|_{\mathcal{Y}^{-1}}&\leq C\big(\|u\|_{\mathcal{Y}^{\gamma}}\big)\|v\|_{\mathcal{Y}^{0}},\label{eq:fuv}\\[2mm]
		\|f(u)-f(v)\|_{\mathcal{Y}^{\gamma-1}}&\leq C\Big(\|u\|_{\mathcal{Y}^{\gamma}},\|v\|_{\mathcal{Y}^{\gamma}} \Big)\|u-v\|^{\alpha}_{\mathcal{Y}^{0}},\label{eq:fu-fv_gamma-1}
	\end{align}
	where $C\big(\|u\|_{\mathcal{Y}^{\gamma}}\big)$ and $C\big(\|u\|_{\mathcal{Y}^{\gamma}},\|v\|_{\mathcal{Y}^{\gamma}}\big)$ denote constants that depend on $\|u\|_{\mathcal{Y}^{\gamma}}$ and both $\|u\|_{\mathcal{Y}^{\gamma}}$ and $\|v\|_{\mathcal{Y}^{\gamma}}$, respectively. Specifically, $\alpha\in(0,1]$ can be chosen as
	\begin{align}\label{eq:choose_alpha}
		\alpha=\left\{\begin{array}{l}
			\frac{1}{2},\qquad\qquad\qquad\qquad\qquad\qquad\quad\;\; \text{for}\quad d=2,\quad \gamma=1,\\[2mm]
			\min\left\{1,\frac{d-2\gamma}{2\gamma}\left(\frac{d+2(1-\gamma)}{d-2\gamma}-\xi_{\gamma}\right)\right\},\quad\, \text{for } d=2,\ 0<\gamma<1,\ \text{or}\ d=3,\ 0<\gamma\le 1.
		\end{array}\right.
	\end{align}
\end{lemma}

\begin{proof}[Proof (sketch)]
	We note that, for \(0\le \ell\le 1\), the embeddings
	\(\mathcal{Y}^{\ell} \hookrightarrow H^{\ell}(\Omega)\) and
	\(H^{-\ell}(\Omega) \hookrightarrow \mathcal{Y}^{-\ell}\) hold.
	Moreover, since \(0<\gamma\le 1\), the estimates follow from the growth assumption
	\eqref{eq:growth_condition}, the mean-value theorem, Hölder's inequality, and Sobolev embeddings.
	As the argument is standard but somewhat technical, we defer the proof to Appendix~\ref{sec:nonlinear}.
\end{proof}

\section{Proof of the convergence result}\label{sec:proof}

In this section, we prove the main convergence theorem. We begin by estimating the local error of the frequency-localized auxiliary scheme~\eqref{eq:semi-discrete}. Using the semigroup difference estimates from Lemma~\ref{lem:estimate_etL-etLh} and the nonlinear estimates from Lemma~\ref{lem:nonlinear_term_estimates}, we then derive error bounds for the fully discrete scheme~\eqref{eq:fully-discrete} in the \(\mathcal{Y}^{0}\times\mathcal{Y}^{-1}\) norm. We also establish uniform bounds for the numerical solution in \(\mathcal{Y}^{\gamma}\times\mathcal{Y}^{\gamma-1}\), which allow us to control the nonlinear contributions. Combining these ingredients, we close the induction argument and obtain the desired error estimate.
\subsection{Consistency error estimate for the auxiliary scheme}\label{subsec:semi-discrete}

Recalling the definitions of the remainder terms in \eqref{eq:R_1} and \eqref{eq:R2R3R4}, we establish the following result under the stated assumptions.
\begin{lemma}\label{lem:estimate_R1R2R3R4}
Under the conditions of Theorem~\ref{thm:main}, the following estimates hold:
	\begin{align}\label{eq:estimate_R1R2R3R4}
		\|R_1(t_n)\|_{0}\lesssim\tau^2,\quad \|R_2(t_n)\|_{0}+\|R_4(t_n)\|_{0}\lesssim \tau r^{-\gamma}, \quad \text{and}\quad \|R_3(t_n)\|_{0}\lesssim \tau^2 r^{1-\gamma}.
	\end{align}
\end{lemma}

\begin{proof}
Using the Taylor expansion, we have
\begin{align}\label{eq:FU-FeU}
	&F\big(U(t_n+s)\big)-F(e^{s L}U(t_n))=\int_0^1 F^{\prime}(\theta U(t_n+s)+(1-\theta)e^{sL}U(t_n))\d \theta\cdot\left(U(t_n+s)-e^{s L}U(t_n)\right)\notag\\
	&\quad=\int_0^1 F^{\prime}(\theta U(t_n+s)+(1-\theta)e^{sL}U(t_n))\d \theta\cdot\int_{0}^{s}e^{(s-\sigma)L}F\big(U(t_n+\sigma)\big)\d \sigma.
\end{align}
Therefore,
\begin{align*}
	&\|R_1(t_n)\|_{0}\lesssim \tau^2 \sup_{\theta,s,\sigma}\left\|F^{\prime}\big(\theta U(t_n+s)+(1-\theta)e^{sL}U(t_n)\big)e^{(s-\sigma)L}F\big(U(t_n+\sigma)\right\|_{0}.
\end{align*}
	Let us denote $U_{\theta}=(U_{\theta,1},U_{\theta,2})^{\top}:=\theta U(t_n+s)+(1-\theta)e^{sL}U(t_n)$, and
	\begin{align*}
		F^{\prime}(U_\theta)=\begin{pmatrix}
			0 & 0 \\ f^{\prime}(U_{\theta,1}) & 0
		\end{pmatrix}.
	\end{align*}
	Using the nonlinear estimate \eqref{eq:fuv}, we find
	\begin{align*}
		&\left\|F^{\prime}\big(\theta U(t_n+s)+(1-\theta)e^{sL}U(t_n)\big)e^{(s-\sigma)L}F\big(U(t_n+\sigma)\right\|_{0}\\
		&=\left\|f^{\prime}(U_{\theta,1})\Big(e^{(s-\sigma)L}F\big(U(t_n+\sigma)\Big)_1\right\|_{\mathcal{Y}^{-1}}\lesssim C(\|U_{\theta,1}\|_{\mathcal{Y}^{\gamma}})\cdot\left\|\Big(e^{(s-\sigma)L}F\big(U(t_n+\sigma)\Big)_1\right\|_{\mathcal{Y}^0}.
	\end{align*}
	where \(\big(e^{(s-\sigma)L}F(U(t_n+\sigma))\big)_1\) denotes the first component of the vector \(e^{(s-\sigma)L}F(U(t_n+\sigma))\).
	
	From the semigroup property \eqref{semigroup_prop_etL} and the estimate \eqref{eq:fu_H-1}, we further obtain
	\begin{align*}
		\left\|\Big(e^{(s-\sigma)L}F\big(U(t_n+\sigma)\Big)_1\right\|_{\mathcal{Y}^0}\lesssim \left\|f\big(u(t_n+\sigma)\right\|_{\mathcal{Y}^{-1}}\leq  C\Big( \|u(t_n+\sigma)\|_{\mathcal{Y}^{\gamma}}\Big)\cdot\big(1+\|u(t_n+\sigma)\|_{\mathcal{Y}^{0}}\big).
	\end{align*}
	
	Moreover, from the definition of \(U_{\theta}\) we have
	\begin{align*}
		\|U_{\theta,1}\|_{\mathcal{Y}^{\gamma}}&\leq \|U_\theta\|_{\gamma}\lesssim \|U(t_n+s)\|_{\gamma}+\|U(t_n)\|_{\gamma}.
	\end{align*}
	Combining the above estimates yields
	\begin{align}\label{eq:estimate_R1}
		\|R_1(t_n)\|_{0}\lesssim \tau^2 C\Big(\|U\|_{C([0,T]; \mathcal{Y}^{\gamma}\times \mathcal{Y}^{\gamma-1})}\Big)\lesssim \tau^2,
		\end{align}
	provided \(U \in C([0,T]; \mathcal{Y}^{\gamma}\times \mathcal{Y}^{\gamma-1})\).
	
	For \(R_2\) and \(R_4\), using the semigroup property \eqref{semigroup_prop_etL} we note that
	\begin{align*}
		\left\|F(e^{s L}U(t_n))-F(e^{s L}\Pi_r U(t_n))\right\|_0=\left\|f(\tilde{u}(t_n+s))-f(\Pi_r \tilde{u}(t_n+s))\right\|_{\mathcal{Y}^{-1}},
	\end{align*}
	 where \((\tilde{u}(t_n+s), \tilde{v}(t_n+s))^\top := e^{sL}(u(t_n), v(t_n))^\top\).
	 
	 The boundedness of \(e^{(\tau-s)L}\) and \(\varphi_1(\tau L)\), together with Bernstein-type inequalities and the nonlinear estimates \eqref{eq:fu-fv} and \eqref{eq:f_gamma-1}, gives
\begin{align*}
	&\|R_2(t_n)\|_{0}+\|R_4(t_n)\|_{0}\\
	&\lesssim \tau \sup_{s} \left\|f(\tilde{u}(t_n+s))-f(\Pi_{r}\tilde{u}(t_n+s))\right\|_{\mathcal{Y}^{-1}}+\tau\sup_{s} \left\|f(\Pi_{r}\tilde{u}(t_n+s))-\Pi_{r} f(\Pi_{r}\tilde{u}(t_n+s))\right\|_{\mathcal{Y}^{-1}}\\
	&\lesssim \tau \sup_{s}C\Big(\|\tilde{u}(t_n+s)\|_{\mathcal{Y}^{\gamma}}, \|\Pi_r\tilde{u}(t_n+s)\|_{\mathcal{Y}^{\gamma}}\Big)\cdot \|\tilde{u}(t_n+s)-\Pi_{r}\tilde{u}(t_n+s)\|_{\mathcal{Y}^{0}}\\
	&\quad +\tau r^{-\gamma} \sup_{s}\|f(\Pi_r \tilde{u}(t_n+s))\|_{\mathcal{Y}^{\gamma-1}}\\
	&\lesssim \tau r^{-\gamma} \sup_{s}C\Big(\|\tilde{u}(t_n+s)\|_{\mathcal{Y}^{\gamma}}, \|\Pi_r\tilde{u}(t_n+s)\|_{\mathcal{Y}^{\gamma}}\Big).
\end{align*}

	For \(R_3\), we compute
	\begin{align}\label{eq:d-eFe}
		\frac{\d}{\d \sigma}F(e^{\sigma L}\Pi_r U(t_n))&=\begin{pmatrix}
			0& 0\\
			f^{\prime}\big(\Pi_r\tilde{u}(t_n+\sigma) )& 0
		\end{pmatrix}\begin{pmatrix}
			0& 1\\
			-A & 0
		\end{pmatrix}\begin{pmatrix}
			\Pi_r\tilde{u}(t_n+\sigma)\\
			\Pi_r\tilde{v}(t_n+\sigma)
		\end{pmatrix}\notag\\
		&=\begin{pmatrix}
			0\\ f^{\prime}\big(\Pi_r\tilde{u}(t_n+\sigma)\big)\cdot\Pi_r\tilde{v}(t_n+\sigma)
		\end{pmatrix} ,
	\end{align}
	where, as before, \((\tilde{u}(t_n+\sigma), \tilde{v}(t_n+\sigma))^\top = e^{\sigma L}(u(t_n), v(t_n))^\top\).
	
	Applying \eqref{eq:d-eFe}, together with estimate \eqref{eq:fuv} and Bernstein’s inequality, we obtain
	\begin{align*}
		&\left\|\frac{\d}{\d \sigma}F(e^{\sigma L}\Pi_r U(t_n))\right\|_0\lesssim \left\|f^{\prime}\big(\Pi_r\tilde{u}(t_n+\sigma)\big)\cdot\Pi_r\tilde{v}(t_n+\sigma)\right\|_{\mathcal{Y}^{-1}}\\
		&\quad\leq C\Big(\|\Pi_r\tilde{u}(t_n+\sigma)\|_{\mathcal{Y}^{\gamma}}\Big) \|\Pi_r \tilde{v}(t_n+\sigma)\|_{\mathcal{Y}^0}\lesssim C\Big(\|\Pi_r\tilde{u}(t_n+\sigma)\|_{\mathcal{Y}^{\gamma}}\Big)\cdot r^{1-\gamma} \|U(t_n)\|_{\gamma}.
	\end{align*}
	Substituting this estimate into \eqref{eq:R3} completes the proof of the bound for \(R_3\) in \eqref{eq:estimate_R1R2R3R4}.
	
\end{proof}

\subsection{Error estimate for the fully discrete scheme}

Next, we employ the expansions of the exact and the numerical solution, and
denote \(\mathcal{L}^n = \sum_{i=1}^4 R_i(t_n)\). By iteratively applying the consistency error formula \eqref{eq:consistency_error}, we conclude that
\begin{align}\label{eq:iterate_U}
	U(t_{n+1})=e^{t_{n+1}L}U(0)+\sum_{j=0}^{n}\tau e^{(t_{n+1}-t_{j+1})L}\varphi_1(\tau L)\Pi_r F(\Pi_r U(t_j))+\sum_{j=0}^n e^{(t_{n+1}-t_{j+1})L} \mathcal{L}^j.
\end{align}
Similarly, for the fully discrete scheme \eqref{eq:fully-discrete}, we have
\begin{align}\label{eq:iterate_Uh}
	U^{n+1}_h=e^{t_{n+1} L_h}\P_hU^0+\sum_{j=0}^n\tau e^{(t_{n+1}-t_{j+1}) L_h}\varphi_1(\tau L_h)\P_h F(U^j_h).
\end{align}
Subtracting \eqref{eq:iterate_Uh} from \eqref{eq:iterate_U}, we obtain
\begin{align}\label{eq:U-Uh}
	&U(t_{n+1})-U^{n+1}_h=e^{t_{n+1}L}U(0)-e^{t_{n+1} L_h}\P_h U^0\notag\\
	&+\sum_{j=0}^{n}\tau \left(e^{(t_{n+1}-t_{j+1})L}\varphi_1(\tau L)\Pi_r F(\Pi_r U(t_j))-e^{(t_{n+1}-t_{j+1}) L_h}\varphi_1(\tau L_h)\P_h F(U^j_h)\right)\notag\\
	&+\sum_{j=0}^n e^{(t_{n+1}-t_{j+1})L} \mathcal{L}^j.
\end{align}

In the next lemma, we estimate the three terms on the right-hand side of \eqref{eq:U-Uh}.

\begin{lemma}\label{lem:estimate_todo}
Under the conditions of Theorem~\ref{thm:main}, and further assume that $1\leq r\leq h^{-1}$,
there exist constants \(C,C_0>0\), such that
\begin{align}\label{eq:error_estimate_Y0}
	\left\|U(t_{n+1})-U^{n+1}_h\right\|_{0}\leq & C_0 (r^{-\gamma}+h^{\ell_k}r^{\ell_k+1-\gamma})+\sum_{j=0}^{n}\tau \cdot C\left(\big\|U^j_h\big\|_{\gamma}\right) \|U(t_j)-U^j_h\|_0\notag\\
	&+\sum_{j=0}^{n}C_0 \tau\cdot  (\tau r^{1-\gamma}+r^{-\gamma}+h^{\ell_k}r^{\ell_k+1-\gamma})+ \sum_{j=0}^{n}C\left(\big\|U^j_h\big\|_{\gamma}\right)\tau \cdot r^{-\gamma},
\end{align}
where \(\ell_k\) is defined in \eqref{eq:def_lk}. The constant \(C_0\) is independent of \(\tau\), \(h\), \(r\), and the bound on the numerical solution, whereas \(C\) depends only on that bound.
\end{lemma}

\begin{proof}
We begin by estimating the last term on the right-hand side of \eqref{eq:U-Uh}. By the semigroup property~\eqref{semigroup_prop_etL} and the semidiscrete consistency estimate in Lemma~\ref{lem:estimate_R1R2R3R4}, we obtain
\begin{align}\label{eq:consistency_error_semi_discrete}
		\left\|e^{(t_{n+1}-t_{j+1})L} \mathcal{L}^j\right\|_0\lesssim \tau^2+\tau r^{-\gamma}+\tau^2 r^{1-\gamma}\lesssim \tau\left( r^{-\gamma}+\tau r^{1-\gamma}\right).
\end{align}

For the first term on the right-hand side of~\eqref{eq:U-Uh}, we decompose
\begin{align}\label{eq:decompose_U0}
	&e^{t_{n+1}L}U(0)-e^{t_{n+1} L_h}U^0_h\notag\\
	&=\underbrace{\big(e^{t_{n+1}L}U(0) - e^{t_{n+1}L}\Pi_r U(0)\big)}_{(I)}+ \underbrace{\big(e^{t_{n+1}L}\Pi_r U(0) - e^{t_{n+1}L_h}\P_h\Pi_r U(0)\big)}_{(II)}
	\notag\\
	&\quad+ \underbrace{\big(e^{t_{n+1}L_h}\P_h\Pi_r U(0) - e^{t_{n+1}L_h}\P_h U(0)\big)}_{(III)}.
\end{align}
For term~(I), by the definition of \(\Pi_r\) and~\eqref{semigroup_prop_etL},
\begin{align}\label{estimate_U_1}
	\left\|e^{t_{n+1}L}U(0)-e^{t_{n+1}L}\Pi_r U(0)\right\|_0\lesssim \left\|(\text{id}-\Pi_r) U(0)\right\|_0\lesssim r^{-\gamma} \|U(0)\|_{\gamma}.
\end{align}
For term~(III), using assumption~\eqref{assump_etLh} and the stability of \(\P_h\) (Assumption~\ref{assumption3} (a)), 
\begin{align}\label{estimate_U_2}
	&\|e^{t_{n+1}L_h}\P_h\Pi_r U(0)-e^{t_{n+1}L_h}\P_h U(0)\|_0\lesssim	\|\P_h\Pi_r U(0)-\P_h U(0)\|_0\notag\\
	&\quad\lesssim \|\P_h\left(\Pi_r u(0)-u(0)\right)\|_{\mathcal{Y}^0}+\|\P_h\left(\Pi_r v(0)-v(0)\right)\|_{\mathcal{Y}^{-1}}\notag\\
	&\quad\lesssim \|\Pi_r U(0)-U(0)\|_0\lesssim r^{-\gamma}\|U(0)\|_{\gamma}.
\end{align}
For term~(II), by Lemma~\ref{lem:estimate_etL-etLh},
\begin{align}\label{estimate_U_3}
	\left\|e^{t_{n+1}L}\Pi_r U(0)-e^{t_{n+1}L_h}\P_h\Pi_r U(0)\right\|_0\lesssim (h^{\gamma}+h^{\ell_k}r^{\ell_k+1-\gamma})\|U(0)\|_\gamma.
\end{align}

Combining \eqref{estimate_U_1}–\eqref{estimate_U_3}, and noting that $r\leq h^{-1}$, gives
\begin{align}\label{eq:estimate_U0}
	\left\|e^{t_{n+1}L}U(0)-e^{t_{n+1} L_h}U^0_h\right\|_0\lesssim \big(r^{-\gamma}+h^{\ell_k}r^{\ell_k+1-\gamma}\big)\|U(0)\|_{\gamma}.
\end{align}

Similarly, for the nonlinear contribution in \eqref{eq:U-Uh}, we write
\begin{align}\label{eq:decopose_F}
	&e^{(t_{n+1}-t_{j+1})L}\varphi_1(\tau L)\Pi_r F(\Pi_r U(t_j))-e^{(t_{n+1}-t_{j+1}) L_h}\varphi_1(\tau L_h)\P_h F(U^j_h)\notag\\
	&=\underbrace{\Big(e^{(t_{n+1}-t_{j+1})L}\varphi_1(\tau L)\Pi_r F(\Pi_r U(t_j))-e^{(t_{n+1}-t_{j+1})L_h}\varphi_1(\tau L_h)\P_h\Pi_r F(\Pi_r U(t_j))\Big)}_{(i)}\notag\\
	&\quad +\underbrace{\Big(e^{(t_{n+1}-t_{j+1})L_h}\varphi_1(\tau L_h)\P_h\Pi_r F(\Pi_r U(t_j))-e^{(t_{n+1}-t_{j+1})L_h}\varphi_1(\tau L_h)\P_h F(\Pi_r U(t_j))\Big)}_{(ii)}\notag\\
	&\quad +\underbrace{\Big(e^{(t_{n+1}-t_{j+1})L_h}\varphi_1(\tau L_h)\P_h F(\Pi_r U(t_j))-e^{(t_{n+1}-t_{j+1})L_h}\varphi_1(\tau L_h)\P_h F(U^j_h)\Big)}_{(iii)}.
\end{align}

For term~(ii), by assumption~\eqref{assump_etLh}, the Bernstein inequality~\eqref{lem:Bernstein}, and the stability of \(\P_h\), together with Lemma~\ref{lem:nonlinear_term_estimates}, 
\begin{align}\label{estimate_F_1}
	&\Big\|e^{(t_{n+1}-t_{j+1})L_h}\varphi_1(\tau L_h)\P_h\Pi_r F(\Pi_r U(t_j))-e^{(t_{n+1}-t_{j+1})L_h}\varphi_1(\tau L_h)\P_h F(\Pi_r U(t_j))\Big\|_0\notag\\
	&\lesssim \Big\|\P_h\Pi_r F(\Pi_r U(t_j))-\P_h F(\Pi_r U(t_j))\Big\|_0\notag\\
	&\lesssim \Big\|\Pi_r F(\Pi_r U(t_j))-F(\Pi_r U(t_j))\Big\|_0=\Big\|\Pi_{r}f(\Pi_{r} u(t_j))-f(\Pi_{r} u(t_j))\Big\|_{\mathcal{Y}^{-1}}\notag\\
	&\lesssim r^{-\gamma}\|f(\Pi_r u(t_j))\|_{\mathcal{Y}^{\gamma-1}}\leq r^{-\gamma}C\Big(\|u(t_j)\|_{\mathcal{Y}^{\gamma}}\Big).
\end{align}
For term~(iii),
\begin{align}\label{estimate_F_2}
	&\Big\|e^{(t_{n+1}-t_{j+1})L_h}\varphi_1(\tau L_h)\P_h F(\Pi_r U(t_j))-e^{(t_{n+1}-t_{j+1})L_h}\varphi_1(\tau L_h)\P_h F(U^j_h)\Big\|_0\notag\\
	&\lesssim \left\|F(\Pi_r U(t_j))-F(U^j_h)\right\|_0=\|f(\Pi_r u(t_j))-f(u^j_h)\|_{\mathcal{Y}^{-1}}\notag\\
	&\lesssim C\Big(\|u(t_j)\|_{\mathcal{Y}^{\gamma}}, \|u^j_h\|_{\mathcal{Y}^{\gamma}}\Big)\|\Pi_r u(t_j)-u^j_h\|_{\mathcal{Y}^0}\notag\\
	&\leq C\Big(\|u(t_j)\|_{\mathcal{Y}^{\gamma}}, \|u^j_h\|_{\mathcal{Y}^{\gamma}}\Big)\cdot\Big(\|\Pi_{>r}u(t_j)\|_{\mathcal{Y}^0}+\| u(t_j)-u^j_h\|_{\mathcal{Y}^0}\Big)\notag\\
	&\lesssim C\Big(\|u(t_j)\|_{\mathcal{Y}^{\gamma}}, \|u^j_h\|_{\mathcal{Y}^{\gamma}}\Big)\cdot\Big(r^{-\gamma}+\|U(t_j)-U^j_h\|_0\Big).
\end{align}

For term~(i), using the integral representation of \(\varphi_1\),
\begin{align*}
	&e^{(t_{n+1}-t_{j+1})L}\varphi_1(\tau L)\Pi_r F(\Pi_r U(t_j))-e^{(t_{n+1}-t_{j+1})L_h}\varphi_1(\tau L_h)\P_h\Pi_r F(\Pi_r U(t_j))\\
	&=\frac{1}{\tau}\int_0^{\tau}\left(e^{(t_{n+1}-t_j-s)L}\Pi_r F(\Pi_r U(t_j))-e^{(t_{n+1}-t_j-s)L_h}\P_h\Pi_r F(\Pi_r U(t_j)) \right)\d s,
\end{align*}
and applying Lemma~\ref{lem:estimate_etL-etLh} with \(V(0) = F(\Pi_r U(t_j))\), together with Lemma~\ref{lem:nonlinear_term_estimates}, yields
\begin{align}\label{estimate_F_3}
	&\Big\|e^{(t_{n+1}-t_{j+1})L}\varphi_1(\tau L)\Pi_r F(\Pi_r U(t_j))-e^{(t_{n+1}-t_{j+1})L_h}\varphi_1(\tau L_h)\P_h\Pi_r F(\Pi_r U(t_j))\Big\|_0\notag\\
	&\lesssim (h^{\gamma}+h^{\ell_k}r^{\ell_k+1-\gamma})\|F(\Pi_r U(t_j))\|_\gamma\notag\\
	&= (h^{\gamma}+h^{\ell_k}r^{\ell_k+1-\gamma})\|f(\Pi_r u(t_j))\|_{\mathcal{Y}^{\gamma-1}}\lesssim (r^{-\gamma}+h^{\ell_k}r^{\ell_k+1-\gamma})\cdot  C\Big(\|u(t_j)\|_{\mathcal{Y}^{\gamma}}\Big),
\end{align}
where, in the last inequality, we have used \(r\le h^{-1}\).

Collecting \eqref{estimate_F_1}–\eqref{estimate_F_3}, we obtain
\begin{align}\label{eq:estimate_F}
	&\left\|e^{(t_{n+1}-t_{j+1})L}\varphi_1(\tau L)\Pi_r F(\Pi_r U(t_j))-e^{(t_{n+1}-t_{j+1}) L_h}\varphi_1(\tau L_h)\P_h F(U^j_h)\right\|_0\notag\\
	& \leq  C\Big(\|u(t_j)\|_{\mathcal{Y}^{\gamma}}, \|u^j_h\|_{\mathcal{Y}^{\gamma}}\Big)\|U(t_j)-U^j_h\|_0+ C\Big(\|u(t_j)\|_{\mathcal{Y}^{\gamma}}, \|u^j_h\|_{\mathcal{Y}^{\gamma}}\Big)r^{-\gamma}\notag\\
	&\quad+(r^{-\gamma}+h^{\ell_k}r^{\ell_k+1-\gamma})C\Big(\|u(t_j)\|_{\mathcal{Y}^{\gamma}}\Big).
\end{align}
Substituting \eqref{eq:estimate_U0}, \eqref{eq:estimate_F}, and the local consistency bound~\eqref{eq:consistency_error_semi_discrete} into~\eqref{eq:U-Uh}, 
and using the regularity \(U\in C([0,T];\mathcal{Y}^{\gamma}\times\mathcal{Y}^{\gamma-1})\), 
we conclude that there exists a constant \(C_0>0\) such that
\eqref{eq:error_estimate_Y0} holds.
%
\end{proof}

\subsection{Estimate of the numerical solution in $\mathcal{Y}^{\gamma}\times\mathcal{Y}^{\gamma-1}$}
At this point, a direct application of Gronwall’s inequality is not yet possible, since the constant on the right-hand side of \eqref{eq:error_estimate_Y0} may depend on \(\|U_h^j\|_{\gamma}\), so such an argument requires a uniform bound for the numerical solution in \(\mathcal{Y}^{\gamma}\times \mathcal{Y}^{\gamma-1}\). Moreover, the estimate in \(\mathcal{Y}^0\times \mathcal{Y}^{-1}\) at the \((n+1)\)-th time step depends on the boundedness of the numerical solution at all previous time steps, that is, on \(\|U_h^j\|_{\gamma}\) for \(j=0,\dots,n\). We therefore estimate the error and the boundedness of the numerical solution simultaneously by induction. To this end, we take the \(\mathcal{Y}^{\gamma}\times\mathcal{Y}^{\gamma-1}\)-norm on both sides of \eqref{eq:fully-discrete}, which yields
\begin{align}\label{eq:U_H1}
	\|U^{n+1}_h\|_{\gamma}&\leq\left\| e^{t_{n+1} L_h}U^0_h\right\|_{\gamma}+\sum_{j=0}^n\tau \left\| e^{(t_{n+1}-t_{j+1}) L_h}\varphi_1(\tau L_h)\P_h F(U^j_h)\right\|_{\gamma}\notag\\
	&\lesssim \left\| U^0_h\right\|_{\gamma}+\sum_{j=0}^n\tau \left\| F(U^j_h)\right\|_{\gamma},
\end{align}
where we used the representation \(\varphi_1(\tau L_h)=\frac{1}{\tau}\int_0^\tau e^{sL_h}\,\mathrm{d}s\), together with the boundedness of \(e^{tL_h}\) and \(\P_h\) from Assumption~\ref{assumption3}(a),(b) and Remark~\ref{rem:assumption_corollary}.  For the second term on the right-hand side of \eqref{eq:U_H1}, 
we insert the term $F(U(t_j))$.  
By applying the estimates \eqref{eq:f_gamma-1} and \eqref{eq:fu-fv_gamma-1} in Lemma~\ref{lem:nonlinear_term_estimates}, we have
\begin{align}\label{eq:IhFhUh}
	\left\| F(U^j_h)\right\|_{\gamma}&\leq \left\|f(u(t_j))\right\|_{\mathcal{Y}^{\gamma-1}}+\left\|f(u(t_j))-f(u^j_h)\right\|_{\mathcal{Y}^{\gamma-1}}\notag\\
	&\leq  C\Big(\left\|u(t_j)\right\|_{\mathcal{Y}^{\gamma}}\Big)+C\Big(\left\|u(t_j)\right\|_{\mathcal{Y}^{\gamma}},\left\|u^j_h\right\|_{\mathcal{Y}^{\gamma}}\Big)\cdot \left\|u(t_j)-u^j_h\right\|^\alpha_{\mathcal{Y}^{0}}.
\end{align}
with $\alpha>0$ given by \eqref{eq:choose_alpha}.

Substituting this estimate into \eqref{eq:U_H1} and noting that $U\in C([0,T];\mathcal{Y}^{\gamma}\times\mathcal{Y}^{\gamma-1})$
we deduce
\begin{align}\label{eq:gronwall_2}
	&\|U^{n+1}_h\|_{\gamma}\leq C_0+\sum_{j=0}^n \tau C\Big(\|U^j_h\|_{\gamma}\Big)\cdot\left\|u(t_j)-u^j_h\right\|^\alpha_{\mathcal{Y}^{0}},
\end{align}
where $C_0$ denotes a constant independent of $\tau$, $h$, $r$, and the numerical solution bound.

\subsection{The proof of Theorem~\ref{thm:main}}
In \eqref{eq:error_estimate_Y0} and \eqref{eq:gronwall_2}, we have obtained estimates for
\(\|U(t_{n+1})-U_h^{n+1}\|_0\) and \(\|U_h^{n+1}\|_{\gamma}\).
However, these bounds rely on the corresponding boundedness assumptions at the previous time steps.
To close the argument, we therefore employ a bootstrap procedure, for which we first state the following technical lemma.
\begin{lemma}[Discrete bootstrap lemma]\label{lem:discrete_bootstrap}
	Let $(a^n)_{n\ge0}$ and $(b^n)_{n\ge0}$ be nonnegative sequences such that
	\begin{align}
		a^{n+1} &\le A + \sum_{j=0}^n \tau\, C(b^j)(\delta+a^j), \label{eq:db_a}\\
		b^{n+1} &\le B + \sum_{j=0}^n \tau\, C(b^j)(\delta+a^j)^\alpha, \label{eq:db_b}\\
		a^0&\le A\qquad b^0\le B,\label{eq:db_ab0}
	\end{align}
	for $n\ge0$, where $A,B,\delta\ge0$, \(\alpha>0\), \(\tau>0\), and
	\(C:[0,\infty)\to[0,\infty)\) is nondecreasing. 
	Given \(T>0\), set
	\[
	\widetilde C:=\sup_{0\le s\le B+1} C(s).
	\]
	If
	\begin{equation}\label{eq:db_smallness}
		\widetilde C T\Bigl(\delta + (A+\widetilde C T\,\delta)e^{T\widetilde C}\Bigr)^\alpha \le 1,
	\end{equation}
	then for all integers \(n\) with \(0\le n\le T/\tau\),
	\begin{equation}\label{eq:db_conclusion}
		b^n \le B+1,
		\qquad
		a^n \le (A+\widetilde C T\,\delta)e^{n\tau\widetilde C}.
	\end{equation}
\end{lemma}

\begin{proof}
	We use induction on \(n\). The claim is immediate for \(n=0\). Assume that
	\[
	b^j\le B+1,
	\qquad
	a^j\le (A+\widetilde C T\,\delta)e^{j\tau\widetilde C},
	\qquad 0\le j\le n.
	\]
	Then \(C(b^j)\le \widetilde C\) for \(0\le j\le n\), and hence
	\[
	b^{n+1}
	\le B+\widetilde C\sum_{j=0}^n \tau\,(\delta+a^j)^\alpha
	\le B+\widetilde C T\Bigl(\delta+(A+\widetilde C T\,\delta)e^{T\widetilde C}\Bigr)^\alpha
	\le B+1,
	\]
	by \eqref{eq:db_smallness}. Moreover,
	\[
	a^{n+1}
	\le A+\widetilde C\sum_{j=0}^n \tau\,(\delta+a^j)
	\le A+\widetilde C T\,\delta+\widetilde C\sum_{j=0}^n \tau\,a^j.
	\]
	By the discrete Gronwall inequality,
	\[
	a^{n+1}\le (A+\widetilde C T\,\delta)e^{(n+1)\tau\widetilde C}.
	\]
	This closes the induction and proves \eqref{eq:db_conclusion}.
\end{proof}
Going back to the estimate \eqref{eq:error_estimate_Y0}, to balance the high-frequency cut-off error and the finite element approximation error, we choose \(r^{-\gamma} = h^{\ell_k} r^{\ell_k+1-\gamma}\), which leads to
\begin{align*}
	r=h^{-\frac{\ell_k}{\ell_k+1}},\quad \text{so that} \quad r^{-\gamma}=h^{\ell_k}r^{\ell_k+1-\gamma}=h^{\frac{\ell_k}{\ell_k+1}\gamma}.
\end{align*}
Set
\begin{align*}
	a^n=\left\|U(t_{n})-U^{n}_h\right\|_{0},\quad \text{and}\quad b^n=\|U^n_h\|_{\gamma}.
\end{align*}
Then, the error estimate
\eqref{eq:error_estimate_Y0} together with the bound \eqref{eq:gronwall_2} implies that
\begin{align}
	a^{n+1}&\leq C_0( \tau h^{\frac{\ell_k}{\ell_k+1}(\gamma-1)} +h^{\frac{\ell_k}{\ell_k+1}\gamma})+\sum_{j=0}^n \tau\cdot C(b^j)\cdot\left(h^{\frac{\ell_k}{\ell_k+1}\gamma}+a^j\right),\label{eq:ineq_a}\\
	b^{n+1}&\leq C_0+\sum_{j=0}^n \tau\cdot C(b^j)\cdot (a^j)^{\alpha},\label{eq:ineq_b}
\end{align}
Moreover, by the projection error estimate in Lemma~\ref{lem:projection} and the boundedness of \(\P_h\) from Assumption~\ref{assumption3} and Remark~\ref{rem:assumption_corollary}, we have
\begin{align}\label{eq:a0b0}
	a^0=\left\|U(0)-\P_h U(0)\right\|_{0}\leq C_0h^{\gamma}
	\quad \text{and}\quad b^0=\|\P_h U(0)\|_{\gamma}\leq C_0.
\end{align}
Here and below, \(C_0\) denotes a generic constant independent of \(\tau\), \(h\), \(r\), and the a priori bound for the numerical solution.
Thus $a^n$ and $b^n$ satisfy the inequalities \eqref{eq:db_a}--\eqref{eq:db_ab0} with 
\begin{align*}
	A=C_0(\tau h^{\frac{\ell_k}{\ell_k+1}(\gamma-1)} +h^{\frac{\ell_k}{\ell_k+1}\gamma}),\quad B=C_0,\quad \delta=h^{\frac{\ell_k}{\ell_k+1}\gamma}.
\end{align*}
where the initial condition \eqref{eq:db_ab0} is satisfied because \(h^\gamma \le h^{\frac{\ell_k}{\ell_k+1}\gamma}\) for $h\in(0,1)$.
Using the step-size condition from Theorem~\ref{thm:main}, namely,
\(\tau h^{\frac{\ell_k}{\ell_k+1}(\gamma-1)} \lesssim h^\mu\)
for an arbitrarily small but fixed \(\mu>0\), we obtain
\begin{align*}
	A\leq C_0(h^{\mu} +h^{\frac{\ell_k}{\ell_k+1}\gamma}).
\end{align*}
Since \(\mu>0\) and \(\frac{\ell_k}{\ell_k+1}\gamma>0\), there exists \(h_0>0\) such that
\begin{align}\label{eq:tau0h0}
	\widetilde{C}T\left[h_0^{\frac{\ell_k}{\ell_k+1}\gamma}
	+\Bigl(C_0\bigl(h_0^{\mu}+h_0^{\frac{\ell_k}{\ell_k+1}\gamma}\bigr)
	+\widetilde{C}Th_0^{\frac{\ell_k}{\ell_k+1}\gamma}\Bigr)e^{T \widetilde{C}}\right]^{\alpha}\le 1.
\end{align}
Hence, for \(0<h<h_0\) and \(\tau\lesssim h^{\frac{\ell_k}{\ell_k+1}(1-\gamma)+\mu}\), the smallness condition \eqref{eq:db_smallness} holds.
By applying Lemma~\ref{lem:discrete_bootstrap} and substituting the definitions of \(A\) and \(\delta\), we conclude that
\begin{align*}
	\sup_{0\leq n\leq T/\tau}\|U(t_n)-U^n_h\|_{0}&\leq\left(C_0(\tau \,h^{\frac{\ell_k}{\ell_k+1}(\gamma-1)} +h^{\frac{\ell_k}{\ell_k+1}\gamma})+\widetilde C Th^{\frac{\ell_k}{\ell_k+1}\gamma}\right)e^{T\widetilde{C}}\\
	&= C_0 e^{T\widetilde{C}}\cdot \tau \,h^{\frac{\ell_k}{\ell_k+1}(\gamma-1)} +\big(C_0+\widetilde{C}T\big)e^{T\widetilde{C}} \cdot h^{\frac{\ell_k}{\ell_k+1}\gamma}.
\end{align*}
This completes the proof of Theorem~\ref{thm:main}.

\section{Numerical experiments}\label{sec:numerical_experiments}


In this section, we present two numerical experiments on two different domains, using quasi-uniform, shape-regular simplicial partitions.
First, we consider the unit cube, for which Assumption~\ref{assump:mesh} is satisfied for all \(k\in\mathbb{N}\), see Remark~\ref{rem:mesh_ass}. 
Second, we consider an equilateral triangle. In this case, \(\theta=\frac{\pi}{3}\) in \eqref{eq:embedded_Yl}, and hence
\[
\mathcal{Y}^{\ell_\Omega}\hookrightarrow H^{\ell_\Omega}(\Omega)
\qquad\text{with}\qquad
\ell_\Omega = 4-\varepsilon,
\]
for an arbitrarily small \(\varepsilon>0\). Therefore, Assumption~\ref{assump:mesh} is fully justified for \(k=1,2\), while for \(k=3\) it fails only by an arbitrarily small margin on quasi-uniform meshes. 

We perform numerical experiments for polynomial degrees \(k=1,2,3\). In the triangular-domain case, the convergence rate for cubic elements is therefore expected to be reduced by an arbitrarily small amount compared with the ideal rate. Nevertheless, since this loss is negligible at the theoretical level, we still expect to observe essentially the predicted convergence rates in all experiments.

\subsection{Implementation}

We first discuss the implementation of the algorithm 
in Python. We first generate the meshes in
\texttt{Gmsh} 
\cite{geuzaine2009gmsh}
and use 
\texttt{ADIOS4DOLFINx}
\cite{ADIOS4DOLFINx2024}
for writing and reading of the meshes and functions. 
We then use the \texttt{DOLFINx} environment \cite{DOLFINx23} to assemble the mass and stiffness matrices $\massmatrix$ and $\stiffmatrix$ defined by 
\begin{equation}
	(\massmatrix)_{i,j} = \int_\Omega \varphi_i \varphi_j \, dx,
	\qquad
	(\stiffmatrix)_{i,j} = \int_\Omega \nabla \varphi_i \cdot  \nabla \varphi_j \, dx ,
\end{equation}
and store them as objects in the \textup{PETSc for Python} library \cite{petsc_2,petsc_1}, which allows for many options in the solution of the upcoming linear systems. 
We denote by
  $(\solutionvec^{n},\solutiondtvec^n)$
  the coefficient vectors of $( u_h^n, v_h^n ) \approx ( u(t_n), u'(t_n) )$,
 %
   and note that the exponential integrator in \eqref{eq:fully-discrete} is equivalent to the exact solution $(\solutionvec,\solutionvec')$  at time $t_{n+1} = t_n + \tau$ of
\begin{equation} \label{eq:sys_for_rat_kry}
 \massmatrix \solutionvec''(t) = - \stiffmatrix \solutionvec(t)  + \loadvectorn, 
	\quad 
	t \in [t_n, t_n +\tau],
	\qquad
	\solutionvec(t_n) = \solutionvec^n,  
	\quad 
	\solutionvec'(t_n) = \solutiondtvec^n, 
\end{equation}
where $\loadvectorn$ denotes the load vector computed from  $\P_h f(u_h(t_n))$,
and the new approximations are defined as $\solutionvec^{n+1} = \solutionvec(t_{n+1})$
and 
$\solutiondtvec^{n+1} = \solutionvec'(t_{n+1})$.
The solution of \eqref{eq:sys_for_rat_kry} is computed using a rational Krylov approximation as suggested in 
\cite{GriH08} and \cite{HocPSTW15}.
For the error computation we use the
$\massmatrix$-inner product for the $L^2(\Omega)$-norm and the $(\stiffmatrix + \massmatrix)^{-1}$-inner product for the $H^{-1}(\Omega)$-norm. Note that along the lines of Lemma~\ref{lemma:equivalence}, one can show that the inner products induced by 
$(\stiffmatrix + \massmatrix)^{-1}$ and $\stiffmatrix^{-1}$ are equivalent.
The code corresponding to the experiments in this section is made publicly available at
\begin{equation*}
	\text{\url{https://github.com/BenjaminDoerich/ExpoFem-LowRegWave}}.
\end{equation*}

\subsection{Indicator function of a square in the unit square}

\begin{figure}
	\centering
	\begin{subfigure}{0.45\textwidth}
	\includegraphics[scale=0.25]{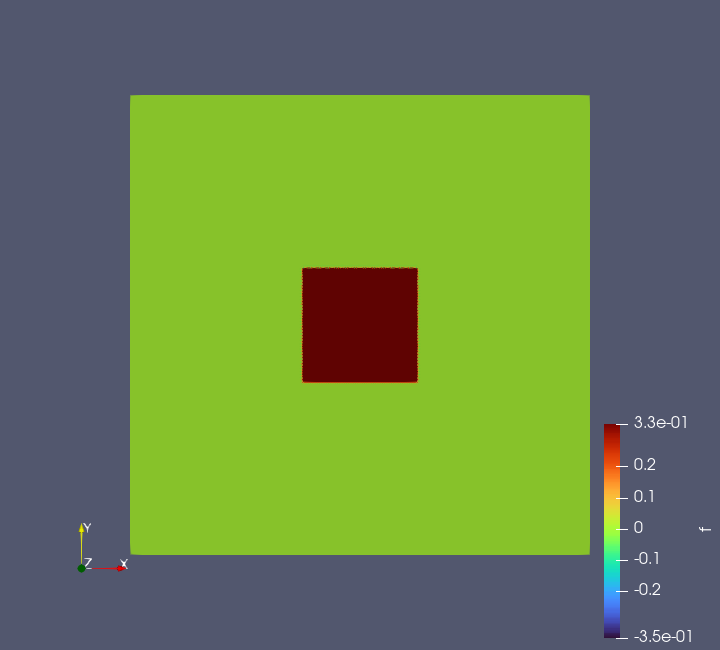}
\end{subfigure}%
	\begin{subfigure}{0.45\textwidth}
	\centering
	\includegraphics[scale=0.25]{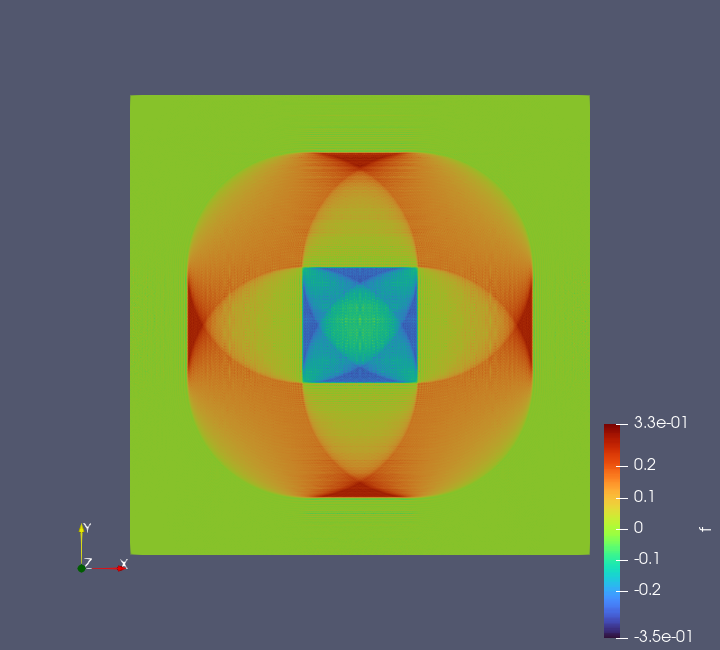}
\end{subfigure}
	\caption{Left: Initial datum defined in \eqref{eq:init_cube_in_cube}. Right: Numerical approximation obtained from the reference solution $u$ at end time $T=0.25$. Both plots are obtained from the simulation with $k=3$.}
	\label{fig:unit_cube_init_and_sol}
\end{figure}

In this subsection, we follow Section~7.2 in \cite{CLLY2024} and consider the nonlinear wave equation defined on $\Omega=[0,1]^2$,
time interval $[0,0.25]$, and $f(u)=4 \sin(u)$ for the initial states $v^{0} = 0$ and
\begin{equation} \label{eq:init_cube_in_cube}
u^{0}(x)=\left\{
	\begin{array}{ll}
		{\displaystyle 0.5},\quad & {\displaystyle x\in \big [0.375,0.625\big]^{2}},\\[1mm]
		{\displaystyle 0},\quad &\text{else. }
	\end{array}
	\right. 
\end{equation}
We provide a plot of the initial state $u^{0}$ in Figure~\ref{fig:unit_cube_init_and_sol}.
Here, we have $u^0 \in H^{1/2 - \varepsilon}(\Omega)$ for any $\varepsilon>0$, and we are in the case of 
Theorem~\ref{thm:main} with $\gamma < \frac12$. 
For the experiments, we choose the scaling for $\ell_k$ defined in \eqref{eq:def_lk}
\begin{equation} \label{eq:relation_h_tau}
h \approx \frac{1}{T} \taufac  \tau^{\frac{\ell_k+1}{\ell_k}} 
\end{equation} 
with some factor $\taufac$ to balance out the error contributions, cf. \eqref{eq:main-1} and Remark~\ref{rem:CFL}, and to remain in a regime of a moderate number of degrees of freedom.  This factor is the same for the numerical approximations as well as the reference solution, computed with some smaller time step size $\tauref$ and finer mesh parameter $\hspaceref$, see Table~\ref{tab:exp1_data} for the precise data. The error $E(T)$ is always computed in the $L^2(\Omega) \times H^{-1}(\Omega)$-norm at the final time $T = 0.25$ and scaled by the $L^2(\Omega) \times H^{-1}(\Omega)$-norm of the reference solution.
%
\begin{table}[h!]
\centering
\begin{tabular}{||c c c c c||} 
 \hline
 $k$ & $\taufac$ & $\tauref$ & $\hspaceref$ &  $\text{dofs}_{\text{ref}}$ \\ [0.5ex] 
 \hline\hline
 1 & 10.8 & $\sim  1\cdot 10^{-3}$ & $\sim 1\cdot 10^{-3}$ & 986,035 \\
 \hline
 2 & 7.2 & $\sim  1\cdot 10^{-3}$  & $\sim 2\cdot 10^{-3}$  &  1,116,329 \\
 \hline
 3 & 7.2 & $\sim  1\cdot 10^{-3}$ & $\sim 3\cdot 10^{-3}$ & 1,526,890 \\  
 \hline
 \end{tabular}
\caption{Choice of the parameters for the experiment corresponding to the initial datum defined in \eqref{eq:init_cube_in_cube}.}
 \label{tab:exp1_data}
\end{table}

We perform two different experiments for this domain. First, we study the convergence under the coupling \eqref{eq:relation_h_tau}, and we expect by Theorem~\ref{thm:main} convergence of order $\tau^{1/2-\varepsilon}$ for any $\varepsilon>0$. This is depicted in Figure~\ref{fig:exp_1_2x2grid} in the upper left plot, where the dashed lines indicate order $\frac12$. The upper right plot shows the relation between the number of degrees of freedom and the error and intends to show the superior scaling for higher order polynomials, since fewer degrees of freedom are needed for $k=2,3$ to achieve the same error as linear elements at least for finer resolutions. 

Second, we compute solutions with a fixed $\tau_0 = 4 \tauref$ and decreasing values of $h$ in order to visualize the predicted order 
$h^{\frac{\ell_k}{\ell_k+1}\frac12}$. This is shown in the lower left part of Figure~\ref{fig:exp_1_2x2grid}.
From our theory, we expect convergence with the orders
$\frac{1}{3}, \frac{2}{5}, \frac{5}{12}$, which is again indicated with the dashed lines. As above also the scaling in the number of degrees of freedom shows a gain in the use of higher order polynomials.
\begin{figure}
\textbf{Case of $H^{0.5-\varepsilon}(\Omega)$-solution}\\[0.5em]
  \centering
  \begin{subfigure}{0.48\textwidth}
    \centering
	\begingroup%
	\StrCount{num_exp/square_in_square_exp/square_in_square_exp_joint_tau}{/}[\matches]%
	\StrBefore[\matches]{num_exp/square_in_square_exp/square_in_square_exp_joint_tau}{/}[\datapath]%
\providecommand{\logLogSlopeTriangle}[5]
{
	
	\pgfplotsextra
	{
		\pgfkeysgetvalue{/pgfplots/xmin}{\xmin}
		\pgfkeysgetvalue{/pgfplots/xmax}{\xmax}
		\pgfkeysgetvalue{/pgfplots/ymin}{\ymin}
		\pgfkeysgetvalue{/pgfplots/ymax}{\ymax}
		
		\pgfmathsetmacro{\xArel}{#1}
		\pgfmathsetmacro{\yArel}{#3}
		\pgfmathsetmacro{\xBrel}{#1-#2}
		\pgfmathsetmacro{\yBrel}{\yArel}
		\pgfmathsetmacro{\xCrel}{\xArel}
		
		\pgfmathsetmacro{\lnxB}{\xmin*(1-(#1-#2))+\xmax*(#1-#2)} 
		\pgfmathsetmacro{\lnxA}{\xmin*(1-#1)+\xmax*#1} 
		\pgfmathsetmacro{\lnyA}{\ymin*(1-#3)+\ymax*#3} 
		\pgfmathsetmacro{\lnyC}{\lnyA+#4*(\lnxA-\lnxB)}
		\pgfmathsetmacro{\yCrel}{\lnyC-\ymin)/(\ymax-\ymin)} 
		
		\coordinate (A) at (rel axis cs:\xArel,\yArel);
		\coordinate (B) at (rel axis cs:\xBrel,\yBrel);
		\coordinate (C) at (rel axis cs:\xCrel,\yCrel);
		
	
	\draw[#5]   (A)-- node[midway, auto, font=\scriptsize] {1}
	(B)-- 
	(C)-- node[midway, auto, font=\scriptsize] {#4}
	cycle;
}
}

    \begin{tikzpicture}
      \begin{axis}[
          width=0.925\linewidth,
          height=0.9\linewidth,
          grid=major,
          grid style={dashed,gray!30},
          xlabel=$\tau$,
          ylabel={$E(T=0.25)$},
          ylabel style={yshift= 0.5em},
          xmode=log,
          ymode=log,
          ymax=2,
          legend style={at={(0.7,0.21)},anchor=north}, 
          x tick label style={rotate=0,anchor=north}, 
          legend columns=-1,
          legend to name=Legendforall,
        ]
		\addplot+[mark=*, color=blue]
         table[x=taus,y=error_tot_1, col sep=comma]
         {\datapath/data/ref_sol_test_num_9_h_0.0026109660574412533_tau_0.0009765625_all_pol_deg_joint_conv_relative_time.txt};
        \addlegendentry{\(k=1\)}

         \addplot+[mark=*, color=red]
         table[x=taus,y=error_tot_2,col sep=comma] 
        {\datapath/data/ref_sol_test_num_9_h_0.0026109660574412533_tau_0.0009765625_all_pol_deg_joint_conv_relative_time.txt};
        \addlegendentry{\(k=2\)}

        \addplot+[mark=*, color=cyan,  mark options={draw=cyan,fill=cyan}]
        table[x=taus,y=error_tot_3,col sep=comma] 
        {\datapath/data/ref_sol_test_num_9_h_0.0026109660574412533_tau_0.0009765625_all_pol_deg_joint_conv_relative_time.txt};
        \addlegendentry{\(k=3\)}

       \addplot+[mark=None, dashed, color=blue]
         table[x=taus,y=line_over_curve_tau_1,col sep=comma] 
        {\datapath/data/ref_sol_test_num_9_h_0.0026109660574412533_tau_0.0009765625_all_pol_deg_joint_conv_relative_time.txt};

        \addplot+[mark=None, dashed, color=red]
         table[x=taus,y=line_over_curve_tau_2,col sep=comma] 
        {\datapath/data/ref_sol_test_num_9_h_0.0026109660574412533_tau_0.0009765625_all_pol_deg_joint_conv_relative_time.txt};

        \addplot+[mark=None, dashed, color=cyan]
         table[x=taus,y=line_over_curve_tau_3,col sep=comma] 
        {\datapath/data/ref_sol_test_num_9_h_0.0026109660574412533_tau_0.0009765625_all_pol_deg_joint_conv_relative_time.txt};
        \logLogSlopeTriangle{0.35}{0.2}{0.14}{0.5}{black, thick};


      \end{axis}
    \end{tikzpicture}
%
	\endgroup%

    \label{fig:exp_1_joint-tau}
  \end{subfigure}%
  \hfill
  \begin{subfigure}{0.48\textwidth}
    \centering
	\begingroup%
	\StrCount{num_exp/square_in_square_exp/square_in_square_exp_joint_dofs}{/}[\matches]%
	\StrBefore[\matches]{num_exp/square_in_square_exp/square_in_square_exp_joint_dofs}{/}[\datapath]%
	\includestandalone[]{num_exp/square_in_square_exp/square_in_square_exp_joint_dofs}%
	\endgroup%

    \label{fig:exp_1_joint-dofs}
  \end{subfigure}
  \par\bigskip
  \begin{subfigure}{0.48\textwidth}
    \centering
	\begingroup%
	\StrCount{num_exp/square_in_square_exp/square_in_square_exp_tau}{/}[\matches]%
	\StrBefore[\matches]{num_exp/square_in_square_exp/square_in_square_exp_tau}{/}[\datapath]%
	\includestandalone[]{num_exp/square_in_square_exp/square_in_square_exp_tau}%
	\endgroup%

    \label{fig:exp_1_onlyh-tau}
  \end{subfigure}%
  \hfill
  \begin{subfigure}{0.48\textwidth}
    \centering
	\begingroup%
	\StrCount{num_exp/square_in_square_exp/square_in_square_exp_dofs}{/}[\matches]%
	\StrBefore[\matches]{num_exp/square_in_square_exp/square_in_square_exp_dofs}{/}[\datapath]%
	\includestandalone[]{num_exp/square_in_square_exp/square_in_square_exp_dofs}%
	\endgroup%

    \label{fig:exp_1_onlyh-dofs}
  \end{subfigure}
  \par\bigskip
  \ref{Legendforall}
  \caption{Error plots for the initial data \eqref{eq:init_cube_in_cube}. Top row: Equilibrated error $E(T = 0.25)$ via \eqref{eq:relation_h_tau} for $k=1,2,3$. 
 In the left plot, the dashed lines indicate order~$\tau^\frac12$.
  Bottom row: Convergence in $h$ for fixed $\tau$. The dashed lines indicate the expected order of convergence $h^{\frac{\ell_k}{\ell_k+1}\frac12}$ from Theorem~\ref{thm:main}.}
  \label{fig:exp_1_2x2grid}
\end{figure}

\subsection{Scaled half sphere in an equilateral triangle}

\begin{figure}
	\centering
	\begin{subfigure}{0.33\textwidth}
	\includegraphics[scale=0.2]{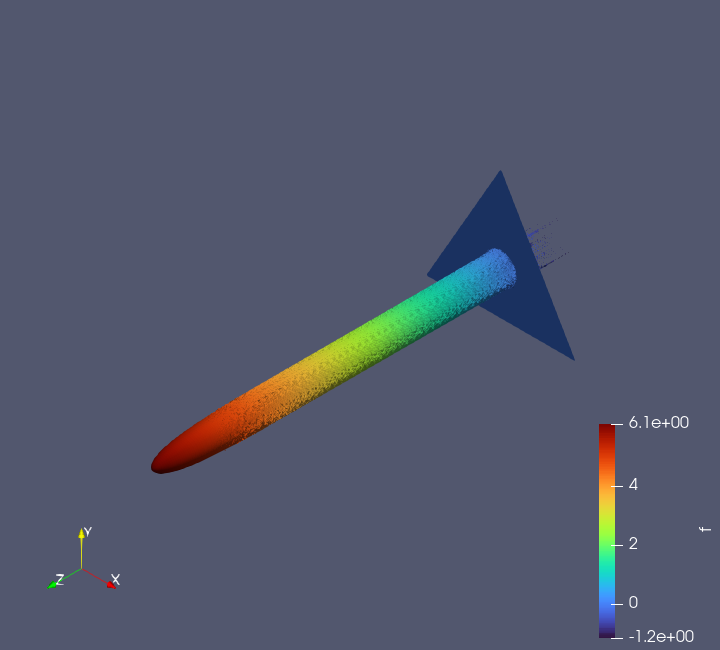}
\end{subfigure}%
	\centering
	\begin{subfigure}{0.33\textwidth}
	\includegraphics[scale=0.2]{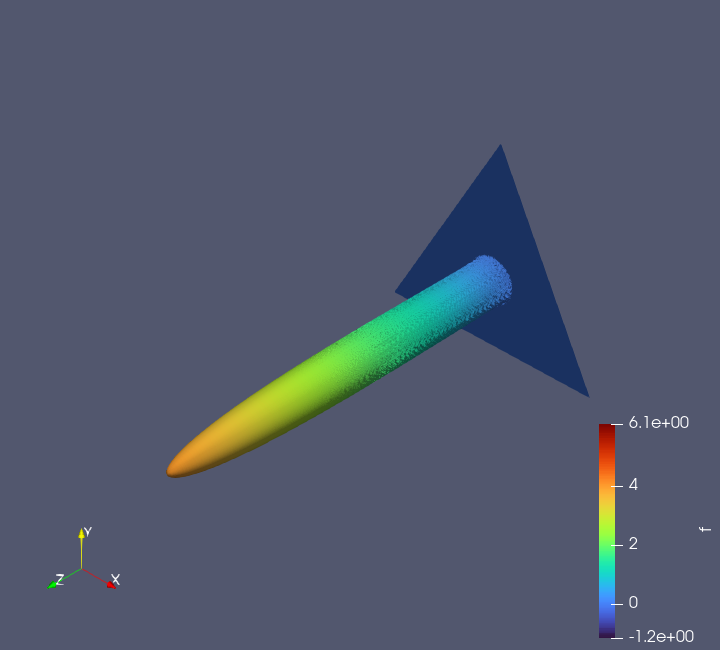}
\end{subfigure}%
	\centering
	\begin{subfigure}{0.33\textwidth}
	\includegraphics[scale=0.2]{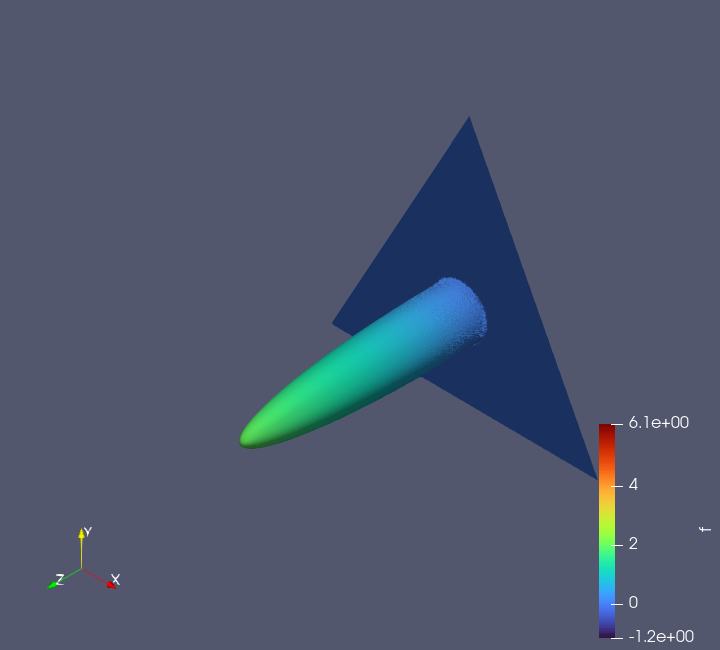}
\end{subfigure}%
	\caption{Initial datum defined in \eqref{eq:init_half_sphere} 
	for the parameters $\alpha = \frac{1}{10}, \frac14, \frac12$ (from left to right)
	computed with the $L^2$-projection and cubic elements.}
	\label{fig:triangle_inits}
\end{figure}

For the second experiment, we focus on a domain given by an equilateral triangle with center at the origin and corners lying on a circle of radius 
$1.2$.
%
We study  the convergence on a scale of regularity parameters $\alpha \in (0,\frac12)$. i.e., we consider the family of initial values, cf. Figure~\ref{fig:triangle_inits}, 
\begin{equation} \label{eq:init_half_sphere}
u_\alpha^{0}(x)=\left\{
	\begin{array}{ll}
		{\displaystyle \numexpheight (r^2 - |x|^2)^{\alpha}},\quad &\text{for } {\displaystyle |x| \leq r},\\[1mm]
		{\displaystyle 0},\quad &\text{elsewhere, }
	\end{array}
	\right.
\end{equation}
for the radius $r = \frac14$ and an amplitude $\numexpheight = 8$. 
Note that these initial values cover the range
 $u_\alpha^{0} \in H^{1/2 + \alpha - \varepsilon}(\Omega)$,
see Corollary~\ref{cor:reg_of_u_alpha} for a proof. As nonlinear term we use $f(u)= u^3$,
and thus by Theorem~\ref{thm:main} we expect convergence of order $\frac12+ \alpha- \varepsilon$. 
In the following, we present the results for $\alpha \in \{\frac{1}{10}, \frac14, \frac12\}$
on the time interval  $[0,0.3]$.

\begin{figure}
	\centering
	\begin{subfigure}{0.45\textwidth}
	\includegraphics[scale=0.3]{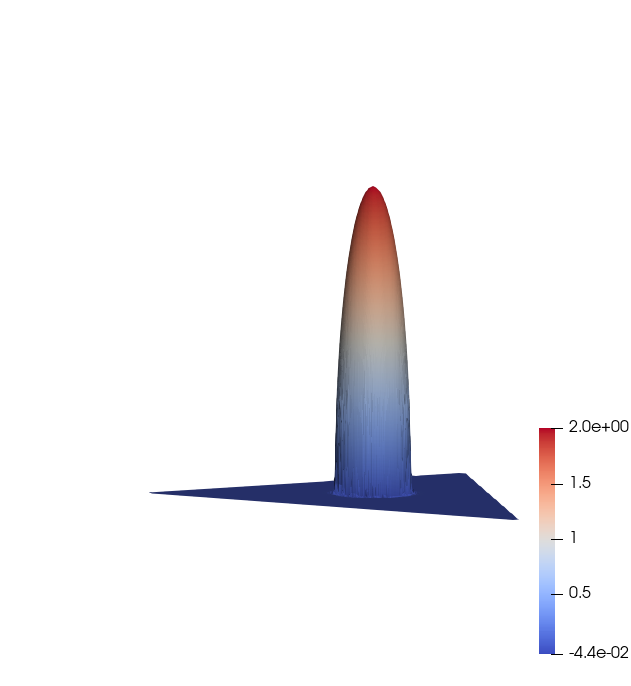}
\end{subfigure}%
	\centering
	\begin{subfigure}{0.45\textwidth}
	\includegraphics[scale=0.3]{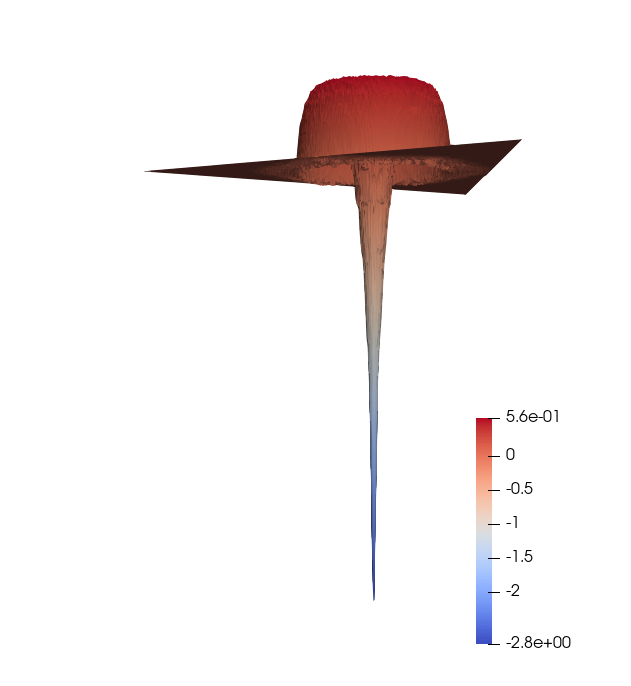}
\end{subfigure}%
	\caption{Plot of the initial datum \eqref{eq:init_half_sphere} 
	for $\alpha = \frac12$ (left) as well as the numerical approximation at time $t \approx 0.25$ computed on a coarse mesh with quadratic elements.}
	\label{fig:triangle_H1_loglog}
\end{figure}

 \begin{remark}\upshape

		For \(\alpha=\tfrac12\), the initial data satisfy $(u^0,v^0)\in H^{1-\varepsilon}(\Omega)\times H^{-\varepsilon}(\Omega)$ for all $\varepsilon>0$.
		In the two-dimensional case, the cubic nonlinearity \(f(u)=u^3\) satisfies the growth condition \eqref{eq:growth_condition}. For sufficiently small initial data, one expects local solutions on a time interval \([0,T_{\max}]\), where \(T_{\max}>0\) depends on the size of the data, for instance on the amplitude parameter \(\numexpheight\), and may become large when \(\numexpheight\) is small.
		
		Moreover, by finite speed of propagation, for points \((x,t)\) such that \(B(x,t)\subset\Omega\), the solution agrees with the corresponding full-space solution and hence admits the standard two-dimensional Duhamel representation formula; see, for example, \cite{evans2022partial}:
		\begin{align}
			u(x,t)
			=
			\frac12\fint_{B(x,t)}
			\frac{t\,u_\alpha^0(y)+t\,\nabla u_\alpha^0(y)\cdot (y-x)}
			{\sqrt{t^2-|y-x|^2}}
			\,\d y
			+
			\frac12\int_0^t \fint_{B(x,t-s)}
			\frac{(t-s)^2u(s,y)^3}
			{\sqrt{(t-s)^2-|y-x|^2}}
			\,\d y\,\d s .
		\end{align}
		Here \(\fint_{B(x,t)}\) denotes the average integral over the ball \(B(x,t)\).
		
		Using the expression of \(u_\alpha^0\) in \eqref{eq:init_half_sphere}, one formally finds that the linear contribution at the origin behaves like
		\begin{align*}
			\left|\fint_{B(0,t)}
			\frac{t\,\nabla u_\alpha^0(y)\cdot y}{\sqrt{t^2-|y|^2}}\,\d y\right|
			\sim
			\numexpheight\,|\ln|t-r||,
			\qquad t\to r.
		\end{align*}
		This indicates the formation of a logarithmic singularity in the linear part at \(t=r\). 
		Accordingly, one expects the full solution to lose boundedness in \(L^\infty\) near \(t=r\), since the nonlinear term does not appear to provide a mechanism for cancelling the singular behaviour of the linear contribution. At the same time, the interaction with the nonlinear term may further amplify and complicate the singular profile. From a physical and analytical point of view, this reflects the combined effect of wave superposition (or energy focusing) and nonlinearity. This illustrates that a solution may remain in the low-regularity space \(H^{1-\varepsilon}(\Omega)\times H^{-\varepsilon}(\Omega)\) while exhibiting singular pointwise behavior, which further motivates the use of low-regularity spaces in both the analysis and the numerical approximation.
		
		This formal prediction is consistent with our numerical observations: in the right panel of Figure~\ref{fig:triangle_H1_loglog}, as \(t\) approaches \(r\), the solution develops a very sharp, high-amplitude spike near the focusing point. This numerically observed peak strongly suggests the onset of the singular behavior described above.
	\end{remark}


We proceed as in the previous experiment and first consider the convergence under the coupling \eqref{eq:relation_h_tau} and reference solution with data from Table~\ref{tab:exp234_data}.
\begin{table}[h!]
\centering
\begin{tabular}{||c c c c c||} 
 \hline
 $k$ & $\taufac$ & $\tauref$ & $\hspaceref$ &  $\text{dofs}_{\text{ref}}$ \\ [0.5ex] 
 \hline\hline
 1 & 15 & $\sim  1\cdot 10^{-3}$ & $\sim 2\cdot 10^{-3}$ &  623,286 \\ 
 \hline
 2 & 10 & $\sim  1\cdot 10^{-3}$  & $\sim 4\cdot 10^{-3}$  & 907,878 \\
 \hline
 3 & 10 & $\sim  1\cdot 10^{-3}$ & $\sim 5\cdot 10^{-3}$ & 1,304,920 \\
 \hline
 \end{tabular}
\caption{Choice of the parameters for the experiment corresponding to the initial datum defined in \eqref{eq:init_half_sphere}.}
 \label{tab:exp234_data}
\end{table}

We observe in the upper left plots  
of Figures~\ref{fig:exp_2_2x2grid}, \ref{fig:exp_3_2x2grid}, and \ref{fig:exp_4_2x2grid} that under the coupling condition \eqref{eq:relation_h_tau}  
 the convergence order roughly behaves as $\tau^\gamma$
with  $\gamma = \frac35, \frac34, 1$, for $\alpha = \frac{1}{10}, \frac14, \frac12$, respectively (indicated by the dashed lines), as predicted. In the lower left plots
we studied again the convergence in $h$ and included dashed lines of order $h^{\frac{\ell_{k}}{\ell_{k}+1}\gamma}$.
%

\begin{figure}
\textbf{Case of $H^{0.6-\varepsilon}(\Omega)$ solution}\\[0.5em]
  \centering
  \begin{subfigure}{0.48\textwidth}
    \centering
	\begingroup%
	\StrCount{num_exp/cirlce_in_triangle_exp_060/cirlce_in_triangle_exp_60_joint_tau}{/}[\matches]%
	\StrBefore[\matches]{num_exp/cirlce_in_triangle_exp_060/cirlce_in_triangle_exp_60_joint_tau}{/}[\datapath]%
	\includestandalone[]{num_exp/cirlce_in_triangle_exp_060/cirlce_in_triangle_exp_60_joint_tau}%
	\endgroup%

    \label{fig:exp_2_joint-tau}
  \end{subfigure}%
  \hfill
  \begin{subfigure}{0.48\textwidth}
    \centering
	\begingroup%
	\StrCount{num_exp/cirlce_in_triangle_exp_060/cirlce_in_triangle_exp_60_joint_dofs}{/}[\matches]%
	\StrBefore[\matches]{num_exp/cirlce_in_triangle_exp_060/cirlce_in_triangle_exp_60_joint_dofs}{/}[\datapath]%
\providecommand{\logLogSlopeTriangle}[5]
{

    \pgfplotsextra
    {
        \pgfkeysgetvalue{/pgfplots/xmin}{\xmin}
        \pgfkeysgetvalue{/pgfplots/xmax}{\xmax}
        \pgfkeysgetvalue{/pgfplots/ymin}{\ymin}
        \pgfkeysgetvalue{/pgfplots/ymax}{\ymax}

        \pgfmathsetmacro{\xArel}{#1}
        \pgfmathsetmacro{\yArel}{#3}
        \pgfmathsetmacro{\xBrel}{#1-#2}
        \pgfmathsetmacro{\yBrel}{\yArel}
        \pgfmathsetmacro{\xCrel}{\xArel}

        \pgfmathsetmacro{\lnxB}{\xmin*(1-(#1-#2))+\xmax*(#1-#2)} 
        \pgfmathsetmacro{\lnxA}{\xmin*(1-#1)+\xmax*#1} 
        \pgfmathsetmacro{\lnyA}{\ymin*(1-#3)+\ymax*#3} 
        \pgfmathsetmacro{\lnyC}{\lnyA+#4*(\lnxA-\lnxB)}
        \pgfmathsetmacro{\yCrel}{\lnyC-\ymin)/(\ymax-\ymin)} 

        \coordinate (A) at (rel axis cs:\xArel,\yArel);
        \coordinate (B) at (rel axis cs:\xBrel,\yBrel);
        \coordinate (C) at (rel axis cs:\xCrel,\yCrel);

        \draw[#5]   (A)-- node[pos=0.5,anchor=north] {1}
                    (B)-- 
                    (C)-- node[pos=0.5,anchor=west] {#4}
                    cycle;
    }
}

    \begin{tikzpicture}
      \begin{axis}[
          width=0.925\linewidth,
          height=0.9\linewidth,
          grid=major,
          grid style={dashed,gray!30},
          xlabel= dofs,
          xmode=log,
          ymode=log,
          ymax=1.1,
          legend style={at={(0.7,0.21)},anchor=north}, 
          legend columns = 2
          x tick label style={rotate=0,anchor=north}, 
          legend columns=-1,
          legend to name=Legendforall
        ]
        \addplot+[mark=*, color=blue]
        table[x=dofs_1,y=error_tot_1,col sep=comma] 
        {\datapath/data/ref_sol_test_num_13_reg_060_tri_h_0.004830917874396135_tau_0.001171875_all_pol_deg_joint_conv_relative_time.txt};
        \addlegendentry{\(k=1\)}

        
        \addplot+[mark=*, color=red]
        table[x=dofs_2,y=error_tot_2,col sep=comma] 
        {\datapath/data/ref_sol_test_num_13_reg_060_tri_h_0.004830917874396135_tau_0.001171875_all_pol_deg_joint_conv_relative_time.txt};
        \addlegendentry{\(k=2\)}
		        
        
        \addplot+[mark=*, color=cyan,  mark options={draw=cyan,fill=cyan}]
      	table[x=dofs_3,y=error_tot_3,col sep=comma] 
		{\datapath/data/ref_sol_test_num_13_reg_060_tri_h_0.004830917874396135_tau_0.001171875_all_pol_deg_joint_conv_relative_time.txt};
		\addlegendentry{\(k=3\)}




      \end{axis}
    \end{tikzpicture}
%
	\endgroup%

    \label{fig:exp_2_joint-dofs}
  \end{subfigure}
  \par\bigskip
  \begin{subfigure}{0.48\textwidth}
    \centering
	\begingroup%
	\StrCount{num_exp/cirlce_in_triangle_exp_060/cirlce_in_triangle_exp_60_tau}{/}[\matches]%
	\StrBefore[\matches]{num_exp/cirlce_in_triangle_exp_060/cirlce_in_triangle_exp_60_tau}{/}[\datapath]%
	\includestandalone[]{num_exp/cirlce_in_triangle_exp_060/cirlce_in_triangle_exp_60_tau}%
	\endgroup%

    \label{fig:exp_2_onlyh-tau}
  \end{subfigure}%
  \hfill
  \begin{subfigure}{0.48\textwidth}
    \centering
	\begingroup%
	\StrCount{num_exp/cirlce_in_triangle_exp_060/cirlce_in_triangle_exp_60_dofs}{/}[\matches]%
	\StrBefore[\matches]{num_exp/cirlce_in_triangle_exp_060/cirlce_in_triangle_exp_60_dofs}{/}[\datapath]%
\providecommand{\logLogSlopeTriangle}[5]
{

    \pgfplotsextra
    {
        \pgfkeysgetvalue{/pgfplots/xmin}{\xmin}
        \pgfkeysgetvalue{/pgfplots/xmax}{\xmax}
        \pgfkeysgetvalue{/pgfplots/ymin}{\ymin}
        \pgfkeysgetvalue{/pgfplots/ymax}{\ymax}

        \pgfmathsetmacro{\xArel}{#1}
        \pgfmathsetmacro{\yArel}{#3}
        \pgfmathsetmacro{\xBrel}{#1-#2}
        \pgfmathsetmacro{\yBrel}{\yArel}
        \pgfmathsetmacro{\xCrel}{\xArel}

        \pgfmathsetmacro{\lnxB}{\xmin*(1-(#1-#2))+\xmax*(#1-#2)} 
        \pgfmathsetmacro{\lnxA}{\xmin*(1-#1)+\xmax*#1} 
        \pgfmathsetmacro{\lnyA}{\ymin*(1-#3)+\ymax*#3} 
        \pgfmathsetmacro{\lnyC}{\lnyA+#4*(\lnxA-\lnxB)}
        \pgfmathsetmacro{\yCrel}{\lnyC-\ymin)/(\ymax-\ymin)} 

        \coordinate (A) at (rel axis cs:\xArel,\yArel);
        \coordinate (B) at (rel axis cs:\xBrel,\yBrel);
        \coordinate (C) at (rel axis cs:\xCrel,\yCrel);

        \draw[#5]   (A)-- node[pos=0.5,anchor=north] {1}
                    (B)-- 
                    (C)-- node[pos=0.5,anchor=west] {#4}
                    cycle;
    }
}

    \begin{tikzpicture}
      \begin{axis}[
         width=0.925\linewidth,
          height=0.9\linewidth,       
          grid=major,
          grid style={dashed,gray!30},
          xlabel=dofs,
          xmode=log,
          ymode=log,
          ymax=1.1,
          legend style={at={(0.7,0.21)},anchor=north}, 
          legend columns = 2
          x tick label style={rotate=0,anchor=north}, 
          legend columns=-1,
          legend to name=Legendforall
        ]
        \addplot+[mark=*, color=blue]
        table[x=dofs_1,y=error_tot_1,col sep=comma] 
        {\datapath/data/ref_sol_test_num_13_reg_060_tri_h_0.004830917874396135_tau_0.001171875_all_pol_deg_h_conv_relative_time.txt};
        \addlegendentry{\(k=1\)}

        
        \addplot+[mark=*, color=red]
         table[x=dofs_2,y=error_tot_2,col sep=comma] 
        {\datapath/data/ref_sol_test_num_13_reg_060_tri_h_0.004830917874396135_tau_0.001171875_all_pol_deg_h_conv_relative_time.txt};
        \addlegendentry{\(k=2\)}

        
        \addplot+[mark=*, color=cyan,  mark options={draw=cyan,fill=cyan}]
        table[x=dofs_3,y=error_tot_3,col sep=comma] 
		{\datapath/data/ref_sol_test_num_13_reg_060_tri_h_0.004830917874396135_tau_0.001171875_all_pol_deg_h_conv_relative_time.txt};
		\addlegendentry{\(k=3\)}




      \end{axis}
    \end{tikzpicture}
%
	\endgroup%

    \label{fig:exp_2_onlyh-dofs}
  \end{subfigure}
  \par\bigskip
  \ref{Legendforall}
  \caption{Error plots for the initial data \eqref{eq:init_half_sphere} with $\alpha = 0.1$. Top row: Equilibrated error $E(T = 0.3)$ via \eqref{eq:relation_h_tau} for $k=1,2,3$. 
  In the left plot, the dashed lines indicate order~$\tau^\frac35$.
  Bottom row: Convergence in $h$ for fixed $\tau$. The dashed lines indicate the expected order of convergence $h^{\frac{\ell_k}{\ell_k+1}\frac35}$ from Theorem~\ref{thm:main}.
  }
  \label{fig:exp_2_2x2grid}
\end{figure}

\begin{figure}
\textbf{Case of $H^{0.75-\varepsilon}(\Omega)$ solution}\\[0.5em]
  \centering
  \begin{subfigure}{0.48\textwidth}
    \centering
	\begingroup%
	\StrCount{num_exp/cirlce_in_triangle_exp_075/cirlce_in_triangle_exp_75_joint_tau}{/}[\matches]%
	\StrBefore[\matches]{num_exp/cirlce_in_triangle_exp_075/cirlce_in_triangle_exp_75_joint_tau}{/}[\datapath]%
	\includestandalone[]{num_exp/cirlce_in_triangle_exp_075/cirlce_in_triangle_exp_75_joint_tau}%
	\endgroup%

    \label{fig:exp_3_joint-tau}
  \end{subfigure}%
  \hfill
  \begin{subfigure}{0.48\textwidth}
    \centering
	\begingroup%
	\StrCount{num_exp/cirlce_in_triangle_exp_075/cirlce_in_triangle_exp_75_joint_dofs}{/}[\matches]%
	\StrBefore[\matches]{num_exp/cirlce_in_triangle_exp_075/cirlce_in_triangle_exp_75_joint_dofs}{/}[\datapath]%
\providecommand{\logLogSlopeTriangle}[5]
{

    \pgfplotsextra
    {
        \pgfkeysgetvalue{/pgfplots/xmin}{\xmin}
        \pgfkeysgetvalue{/pgfplots/xmax}{\xmax}
        \pgfkeysgetvalue{/pgfplots/ymin}{\ymin}
        \pgfkeysgetvalue{/pgfplots/ymax}{\ymax}

        \pgfmathsetmacro{\xArel}{#1}
        \pgfmathsetmacro{\yArel}{#3}
        \pgfmathsetmacro{\xBrel}{#1-#2}
        \pgfmathsetmacro{\yBrel}{\yArel}
        \pgfmathsetmacro{\xCrel}{\xArel}

        \pgfmathsetmacro{\lnxB}{\xmin*(1-(#1-#2))+\xmax*(#1-#2)} 
        \pgfmathsetmacro{\lnxA}{\xmin*(1-#1)+\xmax*#1} 
        \pgfmathsetmacro{\lnyA}{\ymin*(1-#3)+\ymax*#3} 
        \pgfmathsetmacro{\lnyC}{\lnyA+#4*(\lnxA-\lnxB)}
        \pgfmathsetmacro{\yCrel}{\lnyC-\ymin)/(\ymax-\ymin)} 

        \coordinate (A) at (rel axis cs:\xArel,\yArel);
        \coordinate (B) at (rel axis cs:\xBrel,\yBrel);
        \coordinate (C) at (rel axis cs:\xCrel,\yCrel);

        \draw[#5]   (A)-- node[pos=0.5,anchor=north] {1}
                    (B)-- 
                    (C)-- node[pos=0.5,anchor=west] {#4}
                    cycle;
    }
}

    \begin{tikzpicture}
      \begin{axis}[
          width=0.925\linewidth,
          height=0.9\linewidth,
          grid=major,
          grid style={dashed,gray!30},
          xlabel=dofs,
          xmode=log,
          ymode=log,
          ymax=1.1,
          legend style={at={(0.7,0.21)},anchor=north}, 
          x tick label style={rotate=0,anchor=north}, 
          legend columns=-1,
          legend to name=Legendforall,
        ]
       \addplot+[mark=*, color=blue]
        table[x=dofs_1,y=error_tot_1,col sep=comma] 
        {\datapath/data/ref_sol_test_num_13_reg_075_tri_h_0.004830917874396135_tau_0.001171875_all_pol_deg_joint_conv_relative_time.txt};
        \addlegendentry{\(k=1\)}

        
        \addplot+[mark=*, color=red]
        table[x=dofs_2,y=error_tot_2,col sep=comma] 
        {\datapath/data/ref_sol_test_num_13_reg_075_tri_h_0.004830917874396135_tau_0.001171875_all_pol_deg_joint_conv_relative_time.txt};
        \addlegendentry{\(k=2\)}

        
        \addplot+[mark=*, color=cyan,  mark options={draw=cyan,fill=cyan}]
		table[x=dofs_3,y=error_tot_3,col sep=comma] 
		{\datapath/data/ref_sol_test_num_13_reg_075_tri_h_0.004830917874396135_tau_0.001171875_all_pol_deg_joint_conv_relative_time.txt};
		\addlegendentry{\(k=3\)}




      \end{axis}
    \end{tikzpicture}
%
	\endgroup%

    \label{fig:exp_3_joint-dofs}
  \end{subfigure}
  \par\bigskip
  \begin{subfigure}{0.48\textwidth}
    \centering
	\begingroup%
	\StrCount{num_exp/cirlce_in_triangle_exp_075/cirlce_in_triangle_exp_75_tau}{/}[\matches]%
	\StrBefore[\matches]{num_exp/cirlce_in_triangle_exp_075/cirlce_in_triangle_exp_75_tau}{/}[\datapath]%
\providecommand{\logLogSlopeTriangle}[5]
{
	
	\pgfplotsextra
	{
		\pgfkeysgetvalue{/pgfplots/xmin}{\xmin}
		\pgfkeysgetvalue{/pgfplots/xmax}{\xmax}
		\pgfkeysgetvalue{/pgfplots/ymin}{\ymin}
		\pgfkeysgetvalue{/pgfplots/ymax}{\ymax}
		
		\pgfmathsetmacro{\xArel}{#1}
		\pgfmathsetmacro{\yArel}{#3}
		\pgfmathsetmacro{\xBrel}{#1-#2}
		\pgfmathsetmacro{\yBrel}{\yArel}
		\pgfmathsetmacro{\xCrel}{\xArel}
		
		\pgfmathsetmacro{\lnxB}{\xmin*(1-(#1-#2))+\xmax*(#1-#2)} 
		\pgfmathsetmacro{\lnxA}{\xmin*(1-#1)+\xmax*#1} 
		\pgfmathsetmacro{\lnyA}{\ymin*(1-#3)+\ymax*#3} 
		\pgfmathsetmacro{\lnyC}{\lnyA+#4*(\lnxA-\lnxB)}
		\pgfmathsetmacro{\yCrel}{\lnyC-\ymin)/(\ymax-\ymin)} 
		
		\coordinate (A) at (rel axis cs:\xArel,\yArel);
		\coordinate (B) at (rel axis cs:\xBrel,\yBrel);
		\coordinate (C) at (rel axis cs:\xCrel,\yCrel);
		
	
	\draw[#5]   (A)-- node[midway, auto, font=\scriptsize] {1}
	(B)-- 
	(C)-- node[midway, auto, font=\scriptsize] {#4}
	cycle;
}
}

\providecommand{\logLogSlopeTriangleFlip}[5]
{

\pgfplotsextra
{
	\pgfkeysgetvalue{/pgfplots/xmin}{\xmin}
	\pgfkeysgetvalue{/pgfplots/xmax}{\xmax}
	\pgfkeysgetvalue{/pgfplots/ymin}{\ymin}
	\pgfkeysgetvalue{/pgfplots/ymax}{\ymax}
	
	\pgfmathsetmacro{\xArel}{#1}
	\pgfmathsetmacro{\yArel}{#3}
	\pgfmathsetmacro{\xBrel}{#1-#2}
	\pgfmathsetmacro{\yBrel}{\yArel}
	\pgfmathsetmacro{\xCrel}{\xArel}
	
	\pgfmathsetmacro{\lnxB}{\xmin*(1-(#1-#2))+\xmax*(#1-#2)} 
	\pgfmathsetmacro{\lnxA}{\xmin*(1-#1)+\xmax*#1} 
	\pgfmathsetmacro{\lnyA}{\ymin*(1-#3)+\ymax*#3} 
	\pgfmathsetmacro{\lnyC}{\lnyA+#4*(\lnxA-\lnxB)}
	\pgfmathsetmacro{\yCrel}{\lnyC-\ymin)/(\ymax-\ymin)} 
	
	\coordinate (A) at (rel axis cs:\xBrel,\yCrel); 
	\coordinate (B) at (rel axis cs:\xBrel,\yBrel); 
	\coordinate (C) at (rel axis cs:\xCrel,\yCrel); 
	
	\draw[#5]
	(A) -- node[midway, auto, swap, font=\scriptsize] {$1$} (B)
	-- (C)
	-- node[midway, auto , swap,font=\scriptsize] {#4}
	cycle;
}
}

    \begin{tikzpicture}
      \begin{axis}[
          width=0.925\linewidth,
          height=0.9\linewidth,
          grid=major,
          grid style={dashed,gray!30},
          xlabel=$h$,
          ylabel={$E(T=0.3)$},
          ylabel style={yshift= 0.5em},
          xmode=log,
          ymode=log,
          ymax=1.1,
          legend style={at={(0.7,0.21)},anchor=north}, 
          legend columns = 2
          x tick label style={rotate=0,anchor=north}, 
          legend columns=-1,
          legend to name=Legendforall
        ]
        \addplot+[mark=*, color=blue]
        table[x=has_1,y=error_tot_1,col sep=comma] 
        {\datapath/data/ref_sol_test_num_13_reg_075_tri_h_0.004830917874396135_tau_0.001171875_all_pol_deg_h_conv_relative_time.txt};
        \addlegendentry{\(k=1\)}

        
        \addplot+[mark=*, color=red]
        table[x=has_2,y=error_tot_2,col sep=comma] 
        {\datapath/data/ref_sol_test_num_13_reg_075_tri_h_0.004830917874396135_tau_0.001171875_all_pol_deg_h_conv_relative_time.txt};
        \addlegendentry{\(k=2\)}

        
        \addplot+[mark=*, color=cyan,  mark options={draw=cyan,fill=cyan}]
		table[x=has_3,y=error_tot_3,col sep=comma] 
		{\datapath/data/ref_sol_test_num_13_reg_075_tri_h_0.004830917874396135_tau_0.001171875_all_pol_deg_h_conv_relative_time.txt};
		\addlegendentry{\(k=3\)}

        \addplot+[mark=None, dashed, color=blue]
        table[x=has_1,y=line_over_curve_has_1,col sep=comma] 
        {\datapath/data/ref_sol_test_num_13_reg_075_tri_h_0.004830917874396135_tau_0.001171875_all_pol_deg_h_conv_relative_time.txt};
		\logLogSlopeTriangleFlip{0.28}{0.15}{0.36}{0.53}{blue, thick};

		\addplot+[mark=None, dashed, color=red]
		table[x=has_2,
 		y expr=\thisrow{line_over_curve_has_2}*0.95,
		col sep=comma] 
		{\datapath/data/ref_sol_test_num_13_reg_075_tri_h_0.004830917874396135_tau_0.001171875_all_pol_deg_h_conv_relative_time.txt};
		\logLogSlopeTriangleFlip{0.28}{0.15}{0.15}{0.64}{red, thick};

 	     \addplot+[mark=None, dashed, color=cyan]
		table[x=has_3,y=line_over_curve_has_3,col sep=comma] 
		{\datapath/data/ref_sol_test_num_13_reg_075_tri_h_0.004830917874396135_tau_0.001171875_all_pol_deg_h_conv_relative_time.txt};
		\logLogSlopeTriangle{0.42}{0.15}{0.1}{0.67}{cyan, thick};
		
      \end{axis}
    \end{tikzpicture}
%
	\endgroup%

    \label{fig:exp_3_onlyh-tau}
  \end{subfigure}%
  \hfill
  \begin{subfigure}{0.48\textwidth}
    \centering
	\begingroup%
	\StrCount{num_exp/cirlce_in_triangle_exp_075/cirlce_in_triangle_exp_75_dofs}{/}[\matches]%
	\StrBefore[\matches]{num_exp/cirlce_in_triangle_exp_075/cirlce_in_triangle_exp_75_dofs}{/}[\datapath]%
\providecommand{\logLogSlopeTriangle}[5]
{

    \pgfplotsextra
    {
        \pgfkeysgetvalue{/pgfplots/xmin}{\xmin}
        \pgfkeysgetvalue{/pgfplots/xmax}{\xmax}
        \pgfkeysgetvalue{/pgfplots/ymin}{\ymin}
        \pgfkeysgetvalue{/pgfplots/ymax}{\ymax}

        \pgfmathsetmacro{\xArel}{#1}
        \pgfmathsetmacro{\yArel}{#3}
        \pgfmathsetmacro{\xBrel}{#1-#2}
        \pgfmathsetmacro{\yBrel}{\yArel}
        \pgfmathsetmacro{\xCrel}{\xArel}

        \pgfmathsetmacro{\lnxB}{\xmin*(1-(#1-#2))+\xmax*(#1-#2)} 
        \pgfmathsetmacro{\lnxA}{\xmin*(1-#1)+\xmax*#1} 
        \pgfmathsetmacro{\lnyA}{\ymin*(1-#3)+\ymax*#3} 
        \pgfmathsetmacro{\lnyC}{\lnyA+#4*(\lnxA-\lnxB)}
        \pgfmathsetmacro{\yCrel}{\lnyC-\ymin)/(\ymax-\ymin)} 

        \coordinate (A) at (rel axis cs:\xArel,\yArel);
        \coordinate (B) at (rel axis cs:\xBrel,\yBrel);
        \coordinate (C) at (rel axis cs:\xCrel,\yCrel);

        \draw[#5]   (A)-- node[pos=0.5,anchor=north] {1}
                    (B)-- 
                    (C)-- node[pos=0.5,anchor=west] {#4}
                    cycle;
    }
}

    \begin{tikzpicture}
      \begin{axis}[
          width=0.925\linewidth,
          height=0.9\linewidth,          grid=major,
          grid style={dashed,gray!30},
          xlabel=dofs,
          xmode=log,
          ymode=log,
          ymax=1.1,
          legend style={at={(0.7,0.21)},anchor=north}, 
          legend columns = 2
          x tick label style={rotate=0,anchor=north}, 
          legend columns=-1,
          legend to name=Legendforall
        ]
        \addplot+[mark=*, color=blue]
        table[x=dofs_1,y=error_tot_1,col sep=comma] 
        {\datapath/data/ref_sol_test_num_13_reg_075_tri_h_0.004830917874396135_tau_0.001171875_all_pol_deg_h_conv_relative_time.txt};
        \addlegendentry{\(k=1\)}

        
        \addplot+[mark=*, color=red]
        table[x=dofs_2,y=error_tot_2,col sep=comma] 
        {\datapath/data/ref_sol_test_num_13_reg_075_tri_h_0.004830917874396135_tau_0.001171875_all_pol_deg_h_conv_relative_time.txt};
        \addlegendentry{\(k=2\)}

        
        \addplot+[mark=*, color=cyan,  mark options={draw=cyan,fill=cyan}]
		table[x=dofs_3,y=error_tot_3,col sep=comma] 
		{\datapath/data/ref_sol_test_num_13_reg_075_tri_h_0.004830917874396135_tau_0.001171875_all_pol_deg_h_conv_relative_time.txt};
		\addlegendentry{\(k=3\)}




      \end{axis}
    \end{tikzpicture}
%
	\endgroup%

    \label{fig:exp_3_onlyh-dofs}
  \end{subfigure}
  \par\bigskip
  \ref{Legendforall}
  \caption{Error plots for the initial data \eqref{eq:init_half_sphere} with $\alpha = 0.25$. Top row: Equilibrated error $E(T = 0.3)$ via \eqref{eq:relation_h_tau} for $k=1,2,3$. 
  In the left plot, the dashed lines indicate order~$\tau^\frac34$.
  Bottom row: Convergence in $h$ for fixed $\tau$. The dashed lines indicate the expected order of convergence $h^{\frac{\ell_k}{\ell_k+1}\frac34}$ from Theorem~\ref{thm:main}.
  }
 \label{fig:exp_3_2x2grid}
\end{figure}

\begin{figure}
\textbf{Case of $H^{1-\varepsilon}(\Omega)$ solution}\\[0.5em]
  \centering
  \begin{subfigure}{0.48\textwidth}
    \centering
	\begingroup%
	\StrCount{num_exp/cirlce_in_triangle_exp_100/cirlce_in_triangle_exp_100_joint_tau}{/}[\matches]%
	\StrBefore[\matches]{num_exp/cirlce_in_triangle_exp_100/cirlce_in_triangle_exp_100_joint_tau}{/}[\datapath]%
	\includestandalone[]{num_exp/cirlce_in_triangle_exp_100/cirlce_in_triangle_exp_100_joint_tau}%
	\endgroup%

    \label{fig:exp_4_joint-tau}
  \end{subfigure}%
  \hfill
  \begin{subfigure}{0.48\textwidth}
    \centering
	\begingroup%
	\StrCount{num_exp/cirlce_in_triangle_exp_100/cirlce_in_triangle_exp_100_joint_dofs}{/}[\matches]%
	\StrBefore[\matches]{num_exp/cirlce_in_triangle_exp_100/cirlce_in_triangle_exp_100_joint_dofs}{/}[\datapath]%
\providecommand{\logLogSlopeTriangle}[5]
{

    \pgfplotsextra
    {
        \pgfkeysgetvalue{/pgfplots/xmin}{\xmin}
        \pgfkeysgetvalue{/pgfplots/xmax}{\xmax}
        \pgfkeysgetvalue{/pgfplots/ymin}{\ymin}
        \pgfkeysgetvalue{/pgfplots/ymax}{\ymax}

        \pgfmathsetmacro{\xArel}{#1}
        \pgfmathsetmacro{\yArel}{#3}
        \pgfmathsetmacro{\xBrel}{#1-#2}
        \pgfmathsetmacro{\yBrel}{\yArel}
        \pgfmathsetmacro{\xCrel}{\xArel}

        \pgfmathsetmacro{\lnxB}{\xmin*(1-(#1-#2))+\xmax*(#1-#2)} 
        \pgfmathsetmacro{\lnxA}{\xmin*(1-#1)+\xmax*#1} 
        \pgfmathsetmacro{\lnyA}{\ymin*(1-#3)+\ymax*#3} 
        \pgfmathsetmacro{\lnyC}{\lnyA+#4*(\lnxA-\lnxB)}
        \pgfmathsetmacro{\yCrel}{\lnyC-\ymin)/(\ymax-\ymin)} 

        \coordinate (A) at (rel axis cs:\xArel,\yArel);
        \coordinate (B) at (rel axis cs:\xBrel,\yBrel);
        \coordinate (C) at (rel axis cs:\xCrel,\yCrel);

        \draw[#5]   (A)-- node[pos=0.5,anchor=north] {1}
                    (B)-- 
                    (C)-- node[pos=0.5,anchor=west] {#4}
                    cycle;
    }
}

    \begin{tikzpicture}
      \begin{axis}[
          width=0.925\linewidth,
          height=0.9\linewidth,
          grid=major,
          grid style={dashed,gray!30},
          xlabel=dofs,
          xmode=log,
          ymode=log,
          ymax=1.1,
          legend style={at={(0.7,0.21)},anchor=north}, 
          legend columns = 2
          x tick label style={rotate=0,anchor=north}, 
          legend columns=-1,
          legend to name=Legendforall
        ]
        \addplot+[mark=*, color=blue]
        table[x=dofs_1,y=error_tot_1,col sep=comma] 
        {\datapath/data/ref_sol_test_num_13_reg_100_tri_h_0.004830917874396135_tau_0.001171875_all_pol_deg_joint_conv_relative_time.txt};
        \addlegendentry{\(k=1\)}


        \addplot+[mark=*, color=red]
        table[x=dofs_2,y=error_tot_2,col sep=comma] 
        {\datapath/data/ref_sol_test_num_13_reg_100_tri_h_0.004830917874396135_tau_0.001171875_all_pol_deg_joint_conv_relative_time.txt};
        \addlegendentry{\(k=2\)}

        
        \addplot+[mark=*, color=cyan,  mark options={draw=cyan,fill=cyan}]
		table[x=dofs_3,y=error_tot_3,col sep=comma] 
		{\datapath/data/ref_sol_test_num_13_reg_100_tri_h_0.004830917874396135_tau_0.001171875_all_pol_deg_joint_conv_relative_time.txt};
		\addlegendentry{\(k=3\)}




      \end{axis}
    \end{tikzpicture}
%
	\endgroup%

    \label{fig:exp_4_joint-dofs}
  \end{subfigure}
  \par\bigskip
  \begin{subfigure}{0.48\textwidth}
    \centering
	\begingroup%
	\StrCount{num_exp/cirlce_in_triangle_exp_100/cirlce_in_triangle_exp_100_tau}{/}[\matches]%
	\StrBefore[\matches]{num_exp/cirlce_in_triangle_exp_100/cirlce_in_triangle_exp_100_tau}{/}[\datapath]%
	\includestandalone[]{num_exp/cirlce_in_triangle_exp_100/cirlce_in_triangle_exp_100_tau}%
	\endgroup%

    \label{fig:exp_4_onlyh-tau}
  \end{subfigure}%
  \hfill
  \begin{subfigure}{0.48\textwidth}
    \centering
	\begingroup%
	\StrCount{num_exp/cirlce_in_triangle_exp_100/cirlce_in_triangle_exp_100_dofs}{/}[\matches]%
	\StrBefore[\matches]{num_exp/cirlce_in_triangle_exp_100/cirlce_in_triangle_exp_100_dofs}{/}[\datapath]%
\providecommand{\logLogSlopeTriangle}[5]
{

    \pgfplotsextra
    {
        \pgfkeysgetvalue{/pgfplots/xmin}{\xmin}
        \pgfkeysgetvalue{/pgfplots/xmax}{\xmax}
        \pgfkeysgetvalue{/pgfplots/ymin}{\ymin}
        \pgfkeysgetvalue{/pgfplots/ymax}{\ymax}

        \pgfmathsetmacro{\xArel}{#1}
        \pgfmathsetmacro{\yArel}{#3}
        \pgfmathsetmacro{\xBrel}{#1-#2}
        \pgfmathsetmacro{\yBrel}{\yArel}
        \pgfmathsetmacro{\xCrel}{\xArel}

        \pgfmathsetmacro{\lnxB}{\xmin*(1-(#1-#2))+\xmax*(#1-#2)} 
        \pgfmathsetmacro{\lnxA}{\xmin*(1-#1)+\xmax*#1} 
        \pgfmathsetmacro{\lnyA}{\ymin*(1-#3)+\ymax*#3} 
        \pgfmathsetmacro{\lnyC}{\lnyA+#4*(\lnxA-\lnxB)}
        \pgfmathsetmacro{\yCrel}{\lnyC-\ymin)/(\ymax-\ymin)} 

        \coordinate (A) at (rel axis cs:\xArel,\yArel);
        \coordinate (B) at (rel axis cs:\xBrel,\yBrel);
        \coordinate (C) at (rel axis cs:\xCrel,\yCrel);

        \draw[#5]   (A)-- node[pos=0.5,anchor=north] {1}
                    (B)-- 
                    (C)-- node[pos=0.5,anchor=west] {#4}
                    cycle;
    }
}

    \begin{tikzpicture}
      \begin{axis}[
          width=0.925\linewidth,
          height=0.9\linewidth,      
          grid=major,
          grid style={dashed,gray!30},
          xlabel=dofs,
          xmode=log,
          ymode=log,
          ymax=1.1,
          legend style={at={(0.7,0.21)},anchor=north}, 
          legend columns = 2
          x tick label style={rotate=0,anchor=north}, 
          legend columns=-1,
          legend to name=Legendforall
        ]
        \addplot+[mark=*, color=blue]
        table[x=dofs_1,y=error_tot_1,col sep=comma] 
        {\datapath/data/ref_sol_test_num_13_reg_100_tri_h_0.004830917874396135_tau_0.001171875_all_pol_deg_h_conv_relative_time.txt};
        \addlegendentry{\(k=1\)}

        
        \addplot+[mark=*, color=red]
        table[x=dofs_2,y=error_tot_2,col sep=comma] 
        {\datapath/data/ref_sol_test_num_13_reg_100_tri_h_0.004830917874396135_tau_0.001171875_all_pol_deg_h_conv_relative_time.txt};
        \addlegendentry{\(k=2\)}

        
        \addplot+[mark=*, color=cyan,  mark options={draw=cyan,fill=cyan}]
		table[x=dofs_3,y=error_tot_3,col sep=comma] 
		{\datapath/data/ref_sol_test_num_13_reg_100_tri_h_0.004830917874396135_tau_0.001171875_all_pol_deg_h_conv_relative_time.txt};
		\addlegendentry{\(k=3\)}


		

      \end{axis}
    \end{tikzpicture}
%
	\endgroup%

    \label{fig:exp_4_onlyh-dofs}
  \end{subfigure}
  \par\bigskip
  \ref{Legendforall}
\caption{Error plots for the initial data \eqref{eq:init_half_sphere} with $\alpha = 0.5$. Top row: Equilibrated error $E(T = 0.3)$ via \eqref{eq:relation_h_tau} for $k=1,2,3$. 
In the left plot, the dashed lines indicate order~$\tau$.
Bottom row: Convergence in $h$ for fixed $\tau$. The dashed lines indicate the expected order of convergence $h^{\frac{\ell_k}{\ell_k+1}}$ from Theorem~\ref{thm:main}.
}
  \label{fig:exp_4_2x2grid}
\end{figure}

\subsection*{Acknowledgments}

We
would like to thank Tim Buchholz for his valuable input on the \texttt{DOLFINx} 
and \texttt{PETSc} implementation of the finite element solver.

\bibliographystyle{abbrv} 
\bibliography{literature}

\appendix

\section{Proofs of Lemma~\ref{lem:estimate_etL-etLh}}\label{sec:proof_or_etL-etLh}
To prove the estimate of the difference between two semigroups $e^{tL}$ and $e^{tL_h}$ in the weak norm $\mathcal{Y}^0\times \mathcal{Y}^{-1}$, we need first show the following equivalence norm result:
\begin{lemma}[Equivalence between the discrete and continuous $\mathcal{Y}^{-1}$ norms]\label{lemma:equivalence}
	Under the Assumption~\ref{assumption3}(a), there exist constants $c,C>0$, independent of $h$, such that
	\[
	c\,\|w_h\|_{\mathcal{Y}^{-1}}
	\le
	\|A_h^{-1/2}w_h\|_{\mathcal{Y}^{0}}
	\le
	C\,\|w_h\|_{\mathcal{Y}^{-1}}
	\qquad\forall\, w_h\in X^k_h.
	\]
\end{lemma}

\begin{proof}
	For $w_h\in X_h^k$, define
	\[
	\|w_h\|_{-1,h}
	:=
	\sup_{\substack{v_h\in X_h^k\\ \|v_h\|_{\mathcal{Y}^1}=1}}
	\langle w_h,v_h\rangle_H.
	\]
	Since $X_h^k\subset H_0^1(\Omega)=\mathcal{Y}^1$, Poincar\'e's inequality implies that
	\[
	\|w_h\|_{-1,h}
	\sim
	\sup_{0\neq v_h\in X_h^k}
	\frac{\langle w_h,v_h\rangle_H}{\|\nabla v_h\|_{\mathcal{Y}^0}}.
	\]
	
	Let $z_h:=A_h^{-1}w_h$. Then, by the definition of $A_h$,
	\[
	\langle w_h,v_h\rangle_H
	=
	\langle A_h z_h,v_h\rangle_H
	=
	\langle \nabla z_h,\nabla v_h\rangle_H
	\qquad\forall\,v_h\in X_h^k.
	\]
	Therefore,
	\[
	\|w_h\|_{-1,h}
	\sim
	\sup_{0\neq v_h\in X_h^k}
	\frac{\langle \nabla z_h,\nabla v_h\rangle_H}{\|\nabla v_h\|_{\mathcal{Y}^0}}
	=
	\|\nabla z_h\|_{\mathcal{Y}^0},
	\]
	where the last identity follows from the Cauchy--Schwarz inequality and the choice
	$v_h=z_h$ if $z_h\neq 0$.
	
	On the other hand,
	\[
	\|A_h^{-1/2}w_h\|_{\mathcal{Y}^0}^2
	=
	\langle A_h^{-1}w_h,w_h\rangle_H
	=
	\langle z_h,A_h z_h\rangle_H
	=
	\|\nabla z_h\|_{\mathcal{Y}^0}^2.
	\]
	Hence
	\begin{equation}\label{eq:Ah-half-discrete-Hminus1}
		\|A_h^{-1/2}w_h\|_{\mathcal{Y}^0}
		\sim
		\|w_h\|_{-1,h}.
	\end{equation}
	
	Since $X_h^k\subset H_0^1(\Omega)=\mathcal{Y}^1$, we immediately have
	\begin{align}\label{eq:estimate-minus1h}
		\|w_h\|_{-1,h}\le \|w_h\|_{\mathcal{Y}^{-1}}.
	\end{align}
	Conversely, for any $v\in \mathcal{Y}^1$,
	\[
	\langle w_h,v\rangle_H
	=
	\langle w_h,\P_h v\rangle_H
	\le
	\|w_h\|_{-1,h}\,\|\P_h v\|_{\mathcal{Y}^1}
	\le
	C\|w_h\|_{-1,h}\,\|v\|_{\mathcal{Y}^1},
	\]
	where we used that $\P_h$ is the $H$-orthogonal projection onto $X_h^k$ and is
	bounded in $\mathcal{Y}^1$ by Assumption~\ref{assumption3}(a).
	
	Taking the supremum over all $v\in \mathcal{Y}^1$ with $\|v\|_{\mathcal{Y}^1}=1$ yields
	\[
	\|w_h\|_{\mathcal{Y}^{-1}}
	\le
	C\|w_h\|_{-1,h}.
	\]
	Combining this with \eqref{eq:Ah-half-discrete-Hminus1} and \eqref{eq:estimate-minus1h} proves the result.
\end{proof}

\begin{proof}[Proof of \eqref{eq:improved_etL-etLh_Lipschitz} for \(k\ge 2\)]
	We denote
	\begin{align*}
		e^{tL}\Pi_r V(0) &=: V_r(t) = \bigl(v_r(t),\partial_t v_r(t)\bigr),\\
		e^{tL_h}\P_h\Pi_r V(0) &=: V_{h,r}(t) = \bigl(v_{h,r}(t),\partial_t v_{h,r}(t)\bigr).
	\end{align*}
	By the definitions of \(e^{tL}\) and \(e^{tL_h}\), the functions \(v_r(t)\) and \(v_{h,r}(t)\) satisfy
	\begin{align}
		\langle \partial_{tt} v_r(t), w_h\rangle_H + \langle A v_r(t), w_h\rangle_H &= 0, \label{eq:continuous_LW}\\
		\langle \partial_{tt} v_{h,r}(t), w_h\rangle_H + \langle A_h v_{h,r}(t), w_h\rangle_H &= 0, \label{eq:discrete_LW}
	\end{align}
	for all \(w_h\in X_h^k\).
	
	We next decompose the error \(v_r-v_{h,r}\) as
	\begin{equation}\label{eq:vr-vhr}
		v_r(t)-v_{h,r}(t)=v_r(t)-\R_h v_r(t)+\theta_h(t),
		\qquad
		\theta_h(t):=\R_h v_r(t)-v_{h,r}(t).
	\end{equation}
	Hence,
	\begin{align}
		&\left\|e^{tL}\Pi_r V(0)-e^{tL_h}\P_h\Pi_r V(0)\right\|_0
		\lesssim \|v_r(t)-v_{h,r}(t)\|_{\mathcal{Y}^0}
		+ \|\partial_t v_r(t)-\partial_t v_{h,r}(t)\|_{\mathcal{Y}^{-1}} \notag\\
		&\quad\le \|\theta_h(t)\|_{\mathcal{Y}^0}
		+ \|\partial_t\theta_h(t)\|_{\mathcal{Y}^{-1}}
		+ \|v_r(t)-\R_h v_r(t)\|_{\mathcal{Y}^0}
		+ \|\partial_t v_r(t)-\R_h \partial_t v_r(t)\|_{\mathcal{Y}^{-1}}.
		\label{eq:decompose_V}
	\end{align}
	
	Since
	\[
	\langle \partial_{tt} v_r(t), w_h\rangle_H=\langle \P_h\partial_{tt} v_r(t), w_h\rangle_H,
	\qquad
	\langle A v_r(t), w_h\rangle_H=\langle A_h\R_h v_r(t), w_h\rangle_H,
	\]
	subtracting \eqref{eq:discrete_LW} from \eqref{eq:continuous_LW}, and using the fact that \(\partial_{tt}\) commutes with \(\P_h\), yields
	\[
	\langle \partial_{tt}\theta_h(t), w_h\rangle_H
	+ \langle A_h\theta_h(t), w_h\rangle_H
	= -\langle \partial_{tt}\rho_h(t), w_h\rangle_H,
	\qquad \forall\, w_h\in X_h^k,
	\]
	where
	\begin{align}\label{eq:def_rho}
		\rho_h(t):=\P_h v_r(t)-\R_h v_r(t).
	\end{align}
	Since \(\theta_h(t)\in X_h^k\) and \(A_h\) is invertible on \(X_h^k\), we choose $w_h=A_h^{-1}\partial_t\theta_h(t).$
	Using also that \(A_h\) commutes with \(\partial_t\), we obtain
	\begin{align*}
		\frac{\d}{\d t}\left(
		\|A_h^{-1/2}\partial_t\theta_h(t)\|_{\mathcal{Y}^0}^2
		+\|\theta_h(t)\|_{\mathcal{Y}^0}^2
		\right)
		&\le
		\|A_h^{-1/2}\partial_{tt}\rho_h(t)\|_{\mathcal{Y}^0}
		\cdot
		\|A_h^{-1/2}\partial_t\theta_h(t)\|_{\mathcal{Y}^0}\\
		&\le
		\frac12\|A_h^{-1/2}\partial_{tt}\rho_h(t)\|_{\mathcal{Y}^0}^2
		+\frac12\|A_h^{-1/2}\partial_t\theta_h(t)\|_{\mathcal{Y}^0}^2.
	\end{align*}
	An application of Gronwall's inequality therefore gives
	\[
	\|A_h^{-1/2}\partial_t\theta_h(t)\|_{\mathcal{Y}^0}^2
	+\|\theta_h(t)\|_{\mathcal{Y}^0}^2
	\lesssim
	\|A_h^{-1/2}\partial_t\theta_h(0)\|_{\mathcal{Y}^0}^2
	+\|\theta_h(0)\|_{\mathcal{Y}^0}^2
	+
	\sup_{s\in[0,t]}
	\|A_h^{-1/2}\partial_{tt}\rho_h(s)\|_{\mathcal{Y}^0}^2.
	\]
	By the equivalence between the discrete and continuous \(\mathcal{Y}^{-1}\)-norms, we further obtain
	\begin{equation}\label{eq:theta_estimate}
		\|\partial_t\theta_h(t)\|_{\mathcal{Y}^{-1}}^2
		+\|\theta_h(t)\|_{\mathcal{Y}^0}^2
		\lesssim
		\|\partial_t\theta_h(0)\|_{\mathcal{Y}^{-1}}^2
		+\|\theta_h(0)\|_{\mathcal{Y}^0}^2
		+
		\sup_{s\in[0,t]}
		\|\partial_{tt}\rho_h(s)\|_{\mathcal{Y}^{-1}}^2.
	\end{equation}
	
	We next estimate the last term on the right-hand side of \eqref{eq:theta_estimate}. By the definition of \(\rho_h\) in \eqref{eq:def_rho} and the estimates \eqref{eq:estimate_Ph_Hminus1} and \eqref{eq:estimate_Rh_Hminus1}, we deduce, for \(k\geq 2\),
	\begin{align}\label{eq:dual_rho_1}
		\|\partial_{tt}\rho_h(s)\|_{\mathcal{Y}^{-1}}
		&\leq \|\partial_{tt}v_r(s)-\P_h \partial_{tt}v_r(s)\|_{\mathcal{Y}^{-1}}
		+\|\partial_{tt}v_r(s)-\R_h \partial_{tt}v_r(s)\|_{\mathcal{Y}^{-1}}\notag\\
		&\lesssim h^{k+2}\|\partial_{tt}v_r(s)\|_{\mathcal Y^{k+1}}.
	\end{align}
	Moreover, since \(\partial_{tt}v_r(s)=-A v_r(s)\), we have
	\[
	\|\partial_{tt}v_r(s)\|_{\mathcal{Y}^{k+1}}
	\lesssim
	\|v_r(s)\|_{\mathcal{Y}^{k+3}}.
	\]
	Combining this with \eqref{eq:dual_rho_1}, and using the boundedness of \(e^{tL}\) together with the spectral localization of \(\Pi_r\), we deduce
	\begin{align}\label{eq:estimate_rho_h}
		\|\partial_{tt}\rho_h(s)\|_{H^{-1}(\Omega)}
		\lesssim
		h^{k+2}\|v_r(s)\|_{\mathcal{Y}^{k+3}}
		\lesssim
		h^{k+2}r^{k+3-\gamma}
		\|\Pi_r V(0)\|_{\mathcal{Y}^{\gamma}\times\mathcal{Y}^{\gamma-1}}.
	\end{align}
	
	We now estimate the first two terms on the right-hand side of \eqref{eq:theta_estimate}. Since
	\[
	\theta_h(0)=\R_h v_r(0)-\P_h v_r(0),
	\qquad
	\partial_t\theta_h(0)=\R_h\partial_t v_r(0)-\P_h\partial_t v_r(0),
	\]
	by \eqref{eq:estimate_PhRh_L2}, \eqref{eq:estimate_Ph_Hminus1}, and \eqref{eq:estimate_Rh_Hminus1}, we have
	\begin{align}\label{eq:theta_0_dt_theta_0}
		\|\theta_h(0)\|_{\mathcal{Y}^0}+\|\partial_t\theta_h(0)\|_{\mathcal{Y}^{-1}}
		\leq {}&\|v_r(0)-\R_h v_r(0)\|_{\mathcal{Y}^0}
		+\|v_r(0)-\P_h v_r(0)\|_{\mathcal{Y}^0}\notag\\
		&+\|\partial_t v_r(0)-\R_h\partial_t v_r(0)\|_{\mathcal{Y}^{-1}}
		+\|\partial_t v_r(0)-\P_h\partial_t v_r(0)\|_{\mathcal{Y}^{-1}}\notag\\
		\lesssim{}& h^{k+1}\|v_r(0)\|_{\mathcal{Y}^{k+1}}
		+h^{k+2}\|\partial_t v_r(0)\|_{\mathcal{Y}^{k+1}}.
	\end{align}
	Noting that \((v_r(0),\partial_t v_r(0))=\Pi_r V(0)\), we obtain from the Bernstein-type inequality in Lemma~\ref{lem:Bernstein} that
	\begin{align*}
		\|v_r(0)\|_{\mathcal{Y}^{k+1}}
		\lesssim
		r^{k+1-\gamma}\|\Pi_r V(0)\|_{\mathcal{Y}^{\gamma}\times\mathcal{Y}^{\gamma-1}},
		\qquad
		\|\partial_t v_r(0)\|_{\mathcal{Y}^{k+1}}
		\lesssim
		r^{k+2-\gamma}\|\Pi_r V(0)\|_{\mathcal{Y}^{\gamma}\times\mathcal{Y}^{\gamma-1}}.
	\end{align*}
	Therefore, combining this with \eqref{eq:theta_0_dt_theta_0}, and using the constraint  \(h\leq r^{-1}\), we obtain
	\begin{align}\label{eq:estimate_theta0}
		\|\theta_h(0)\|_{\mathcal{Y}^0}+\|\partial_t\theta_h(0)\|_{\mathcal{Y}^{-1}}
		&\lesssim
		\left(h^{k+1}r^{k+1-\gamma}+h^{k+2}r^{k+2-\gamma}\right)
		\|\Pi_r V(0)\|_{\mathcal{Y}^{\gamma}\times\mathcal{Y}^{\gamma-1}} \notag\\
		&\lesssim
		h^{\gamma}\,
		\|\Pi_r V(0)\|_{\mathcal{Y}^{\gamma}\times\mathcal{Y}^{\gamma-1}}.
	\end{align}
	Combining \eqref{eq:theta_estimate}, \eqref{eq:estimate_rho_h}, and \eqref{eq:estimate_theta0}, we deduce
	\begin{align}\label{eq:estimate_theta_h}
		\|\partial_t\theta_h(t)\|_{\mathcal{Y}^{-1}}
		+\|\theta_h(t)\|_{\mathcal{Y}^0}
		\lesssim
		\big(h^{k+2}r^{k+3-\gamma}+h^{\gamma}\big)
		\|\Pi_r V(0)\|_{\mathcal{Y}^{\gamma}\times\mathcal{Y}^{\gamma-1}}.
	\end{align}
	Recalling the decomposition of the error \(v_r-v_{h,r}\) in \eqref{eq:vr-vhr}, we have
	\begin{align}\label{eq:estimate_etL}
		\left\|e^{t L} \Pi_r V(0) - e^{t L_h} \P_h \Pi_r V(0)\right\|_0
		\lesssim{}& \|v_r(t)-\R_h v_r(t)\|_{\mathcal{Y}^0}
		+\|\partial_{t}v_r(t)-\R_h \partial_{t}v_r(t)\|_{\mathcal{Y}^{-1}}
		\notag\\
		&+\|\theta_h(t)\|_{\mathcal{Y}^0}
		+\|\partial_t\theta_h(t)\|_{\mathcal{Y}^{-1}}.
	\end{align}
	The first two terms can be estimated in the same way as above by using \eqref{eq:estimate_PhRh_L2}, \eqref{eq:estimate_Rh_Hminus1}, and the boundedness of \(e^{tL}\):
	\begin{align}\label{eq:vr-Rhvr}
		&\|v_r(t)-\R_h v_r(t)\|_{\mathcal{Y}^0}
		+\|\partial_{t}v_r(t)-\R_h \partial_{t}v_r(t)\|_{\mathcal{Y}^{-1}}\notag\\
		&\quad\lesssim
		h^{k+1}\|v_r(t)\|_{\mathcal{Y}^{k+1}}
		+h^{k+2}\|\partial_{t}v_r(t)\|_{\mathcal{Y}^{k+1}}\notag\\
		&\quad\lesssim
		h^{\gamma}\,
		\|e^{tL}\Pi_r V(0)\|_{\mathcal{Y}^{\gamma}\times\mathcal{Y}^{\gamma-1}}
		\lesssim
		h^{\gamma}\|\Pi_r V(0)\|_{\mathcal{Y}^{\gamma}\times\mathcal{Y}^{\gamma-1}}.
	\end{align}
	Together with \eqref{eq:estimate_theta_h} and \eqref{eq:estimate_etL}, this implies
	\[
	\left\|e^{tL}\Pi_r V(0)-e^{tL_h}\P_h\Pi_r V(0)\right\|_0
	\lesssim
	\left(
	h^{k+2}r^{k+3-\gamma}
	+
	h^{\gamma}
	\right)
	\|\Pi_r V(0)\|_{\mathcal{Y}^{\gamma}\times\mathcal{Y}^{\gamma-1}},
	\]
	which proves \eqref{eq:improved_etL-etLh_Lipschitz} for \(k\ge2\).
\end{proof}

\begin{proof}[Proof of \eqref{eq:improved_etL-etLh_Lipschitz} for \(k=1\)]
	When \(k=1\), the \(\mathcal{Y}^{-1}\)-estimate for the Ritz projection is no longer optimal. Therefore, unlike the case \(k\geq 2\), we need to reformulate the argument so as to rely on \(\mathcal{Y}^{0}\)-estimates for the Ritz projection. To this end, we follow the argument in \cite{baker1976error}. We continue to use the notation
	\begin{align*}
		e^{tL}\Pi_r V(0) &=: V_r(t)=\bigl(v_r(t),\partial_t v_r(t)\bigr),\\
		e^{tL_h}\P_h\Pi_r V(0) &=: V_{h,r}(t)=\bigl(v_{h,r}(t),\partial_t v_{h,r}(t)\bigr),
	\end{align*}
	and the decomposition
	\begin{equation}\label{eq:k1_split}
		v_r(t)-v_{h,r}(t)=v_r(t)-\R_h v_r(t)+\theta_h(t),
		\qquad
		\theta_h(t):=\R_h v_r(t)-v_{h,r}(t).
	\end{equation}
	Then
	\begin{align}
		\left\|e^{tL}\Pi_rV(0)-e^{tL_h}\P_h\Pi_rV(0)\right\|_0
		&\lesssim
		\|\theta_h(t)\|_{\mathcal{Y}^0}
		+\|\partial_t\theta_h(t)\|_{\mathcal{Y}^{-1}} \notag\\
		&\quad
		+\|v_r(t)-\R_h v_r(t)\|_{\mathcal{Y}^0}
		+\|\partial_{t}v_r(t)-\R_h \partial_{t}v_r(t)\|_{\mathcal{Y}^{-1}}.
		\label{eq:k1_decomp}
	\end{align}
	
	As in the case \(k\ge2\), the error term \(\theta_h\) satisfies
	\begin{equation}\label{eq:k1_theta_eq}
		\langle \partial_{tt}\theta_h(t),w_h\rangle_H
		+\langle A_h\theta_h(t),w_h\rangle_H
		=
		-\langle \partial_{tt}\rho_h(t),w_h\rangle_H,
		\qquad \forall\, w_h\in X_h^1,
	\end{equation}
	where \(\rho_h(t):=\P_h v_r(t)-\R_h v_r(t)\).
	Define
	\[
	\Theta_h(t):=\int_0^t \theta_h(s)\,ds.
	\]
	Integrating \eqref{eq:k1_theta_eq} from \(0\) to \(t\), we obtain
	\begin{equation}\label{eq:k1_theta_integrated}
		\langle \partial_t\theta_h(t),w_h\rangle_H
		+\langle A_h\Theta_h(t),w_h\rangle_H
		=
		-\langle \partial_t\rho_h(t)-\partial_t\rho_h(0),w_h\rangle_H
		+\langle \partial_t\theta_h(0),w_h\rangle_H
	\end{equation}
	for all \(w_h\in X_h^1\).
	
	Choosing \(w_h=\theta_h(t)=\partial_t\Theta_h(t)\), we obtain
	\begin{align*}
		\frac12\frac{\d}{\d t}\|\theta_h(t)\|_{\mathcal{Y}^0}^2
		+\frac12\frac{\d}{\d t}a(\Theta_h(t),\Theta_h(t))
		&\le
		\Bigl(
		\|\partial_t\rho_h(t)-\partial_t\rho_h(0)\|_{\mathcal{Y}^0}
		+\|\partial_t\theta_h(0)\|_{\mathcal{Y}^0}
		\Bigr)
		\|\theta_h(t)\|_{\mathcal{Y}^0}.
	\end{align*}
	Hence, by Gronwall's inequality,
	\begin{align}\label{eq:k1_theta_energy}
		&\|\theta_h(t)\|_{\mathcal{Y}^0}+\|\nabla\Theta_h(t)\|_{\mathcal{Y}^0}\notag\\
		&\quad\lesssim
		\|\partial_t\rho_h(0)\|_{\mathcal{Y}^0}
		+\|\partial_t\theta_h(0)\|_{\mathcal{Y}^0}
		+\|\theta_h(0)\|_{\mathcal{Y}^0}
		+\sup_{s\in[0,t]}
		\|\partial_t\rho_h(s)\|_{\mathcal{Y}^0},
	\end{align}
	where we used that \(\Theta_h(0)=0\).
	
	We now estimate the terms on the right-hand side of \eqref{eq:k1_theta_energy}. By the error estimate \eqref{eq:estimate_PhRh_L2} in Lemma~\ref{lem:PhRh},
	\begin{equation}\label{eq:k1_rhot_L2}
		\|\partial_t\rho_h(t)\|_{\mathcal{Y}^0}
		\leq 
		\|(\mathrm{id}-\P_h)\partial_t v_r(t)\|_{\mathcal{Y}^0}
		+\|(\mathrm{id}-\R_h)\partial_t v_r(t)\|_{\mathcal{Y}^0}
		\lesssim
		h^2\|\partial_t v_r(t)\|_{\mathcal Y^2}.
	\end{equation}
	Similarly, since \(\partial_t\theta_h(0)=\R_h\partial_t v_r(0)-\P_h\partial_t v_r(0)\), we have
	\begin{equation}\label{eq:k1_theta0}
		\|\partial_t\theta_h(0)\|_{\mathcal{Y}^0}
		\le
		\|(\mathrm{id}-\R_h)\partial_t v_r(0)\|_{\mathcal{Y}^0}
		+
		\|(\mathrm{id}-\P_h)\partial_t v_r(0)\|_{\mathcal{Y}^0}
		\lesssim
		h^2\|\partial_t v_r(0)\|_{\mathcal Y^2},
	\end{equation}
	and, since \(\theta_h(0)=\R_hv_r(0)-\P_hv_r(0)\),
	\begin{align}\label{eq:k1_theta0_1}
		\|\theta_h(0)\|_{L^2(\Omega)}
		\lesssim
		h^2\|v_r(0)\|_{\mathcal Y^2}.
	\end{align}
	Using the boundedness of \(e^{tL}\) and the spectral localization of \(\Pi_r\), we have
	\[
	\|\partial_t v_r(t)\|_{\mathcal Y^2}
	\lesssim
	r^{3-\gamma}\|\Pi_rV(0)\|_{\mathcal Y^\gamma\times\mathcal Y^{\gamma-1}},
	\qquad
	\|v_r(t)\|_{\mathcal Y^2}
	\lesssim
	r^{2-\gamma}\|\Pi_rV(0)\|_{\mathcal Y^\gamma\times\mathcal Y^{\gamma-1}}.
	\]
	Substituting these bounds into \eqref{eq:k1_rhot_L2}, \eqref{eq:k1_theta0}, and \eqref{eq:k1_theta0_1}, and noting that \(r\geq 1\), we arrive at
	\[
	\sup_{s\in[0,t]}
	\|\partial_t\rho_h(s)\|_{\mathcal{Y}^0}
	+\|\partial_t\rho_h(0)\|_{\mathcal{Y}^0}
	+\|\partial_t\theta_h(0)\|_{\mathcal{Y}^0}
	+\|\theta_h(0)\|_{\mathcal{Y}^0}
	\lesssim
	h^2 r^{3-\gamma}
	\|\Pi_rV(0)\|_{\mathcal Y^\gamma\times\mathcal Y^{\gamma-1}}.
	\]
	Hence, by \eqref{eq:k1_theta_energy},
	\begin{equation}\label{eq:k1_theta_final}
		\|\theta_h(t)\|_{\mathcal{Y}^0}
		\lesssim
		h^2 r^{3-\gamma}
		\|\Pi_rV(0)\|_{\mathcal Y^\gamma\times\mathcal Y^{\gamma-1}}.
	\end{equation}
	
	Next, we estimate \(\partial_t\theta_h(t)\) in \(\mathcal{Y}^{-1}\). By \eqref{eq:k1_theta_integrated} and the norm equivalence in Lemma~\ref{lemma:equivalence},
	\begin{align}\label{eq:k1_theta_t_Hm1}
		\|\partial_t\theta_h(t)\|_{\mathcal{Y}^{-1}}
		&\lesssim \sup_{ w_h\in X_h^1\atop \|w_h\|_{\mathcal{Y}^1}=1}\langle \partial_t\theta_h(t), w_h\rangle_{H}\notag\\
		&=\sup_{w_h\in X_h^1\atop \|w_h\|_{\mathcal{Y}^1}=1}\Big(-\langle A_h\Theta_h(t),w_h\rangle_H
		-\langle \partial_t\rho_h(t)-\partial_t\rho_h(0),w_h\rangle_H
		+\langle \partial_t\theta_h(0),w_h\rangle_H\Big)\notag\\
		&\lesssim
		\|\nabla\Theta_h(t)\|_{\mathcal{Y}^0}
		+\|\partial_t\rho_h(t)\|_{\mathcal{Y}^0}
		+\|\partial_t\rho_h(0)\|_{\mathcal{Y}^0}
		+\|\partial_t\theta_h(0)\|_{\mathcal{Y}^0}.
	\end{align}
	Combining \eqref{eq:k1_theta_energy}, \eqref{eq:k1_rhot_L2}, \eqref{eq:k1_theta0}, and \eqref{eq:k1_theta_t_Hm1}, we infer
	\begin{equation}\label{eq:k1_theta_final_pre}
		\|\theta_h(t)\|_{\mathcal{Y}^{0}}
		+\|\partial_t\theta_h(t)\|_{\mathcal{Y}^{-1}}
		\lesssim
		h^2 r^{3-\gamma}
		\|\Pi_rV(0)\|_{\mathcal Y^\gamma\times\mathcal Y^{\gamma-1}}.
	\end{equation}
	
	Finally, applying the error estimate \eqref{eq:estimate_PhRh_L2}, we obtain
	\begin{align*}
		&\|v_r(t)-\R_h v_r(t)\|_{\mathcal{Y}^0}
		+\|\partial_{t}v_r(t)-\R_h \partial_{t}v_r(t)\|_{\mathcal{Y}^{-1}}\\
		&\quad\leq \|v_r(t)-\R_h v_r(t)\|_{\mathcal{Y}^0}
		+\|\partial_{t}v_r(t)-\R_h \partial_{t}v_r(t)\|_{\mathcal{Y}^{0}}
		\lesssim
		h^2\big(\|v_r(t)\|_{\mathcal Y^2}
		+\|\partial_t v_r(t)\|_{\mathcal Y^2}\big).
	\end{align*}
	Using again the boundedness of \(e^{tL}\) and the spectral localization of \(\Pi_r\), we obtain
	\[
	\|v_r(t)\|_{\mathcal Y^2}
	\lesssim
	r^{2-\gamma}\|\Pi_rV(0)\|_{\mathcal Y^\gamma\times\mathcal Y^{\gamma-1}},
	\qquad
	\|\partial_t v_r(t)\|_{\mathcal Y^2}
	\lesssim
	r^{3-\gamma}\|\Pi_rV(0)\|_{\mathcal Y^\gamma\times\mathcal Y^{\gamma-1}}.
	\]
	Therefore,
	\begin{align}
		\|v_r(t)-\R_h v_r(t)\|_{\mathcal{Y}^0}
		+\|\partial_{t}v_r(t)-\R_h \partial_{t}v_r(t)\|_{\mathcal{Y}^{-1}}
		&\lesssim
		\left(h^2r^{2-\gamma}+h^2r^{3-\gamma}\right)
		\|\Pi_rV(0)\|_{\mathcal Y^\gamma\times\mathcal Y^{\gamma-1}} \notag\\
		&\lesssim
		h^2r^{3-\gamma}
		\|\Pi_rV(0)\|_{\mathcal Y^\gamma\times\mathcal Y^{\gamma-1}},
		\label{eq:k1_rho_final}
	\end{align}
	where we used \(r\geq 1\).
	
	Finally, combining \eqref{eq:k1_decomp}, \eqref{eq:k1_theta_final_pre}, and \eqref{eq:k1_rho_final}, we conclude that
	\[
	\left\|e^{t L} \Pi_r V(0) - e^{t L_h} \P_h \Pi_r V(0)\right\|_0
	\lesssim h^2 r^{3-\gamma}
	\|\Pi_rV(0)\|_{\mathcal Y^\gamma\times\mathcal Y^{\gamma-1}},
	\]
	which proves \eqref{eq:improved_etL-etLh_Lipschitz} for \(k=1\).
\end{proof}

\section{Proofs of the estimates for the nonlinear terms}\label{sec:nonlinear}
{\sc Proof of \eqref{eq:f_gamma-1}}: Since \(\gamma\in(0,1]\), it follows from \eqref{eq:relation_Yl} and the Sobolev embedding \(H_0^{1-\gamma}\hookrightarrow L^{\frac{2d}{d-2(1-\gamma)}}\) that \(L^{\frac{2d}{d+2(1-\gamma)}}\hookrightarrow (H_0^{1-\gamma})^*=H^{\gamma-1}\). Therefore,
\begin{align*}
	\|f(u)\|_{\mathcal{Y}^{\gamma-1}}\lesssim\|f(u)\|_{H^{\gamma-1}}\lesssim\|f(u)\|_{L^{\frac{2d}{d+2(1-\gamma)}}}
	\lesssim \|1+|u|^{\xi_{\gamma}}\|_{L^{\frac{2d}{d+2(1-\gamma)}}}\lesssim 1+\|u\|_{L^{\frac{2d\xi_{\gamma}}{d+2(1-\gamma)}}}^{\xi_{\gamma}}.
\end{align*}
Moreover, by the Sobolev embedding \(\mathcal{Y}^\gamma \hookrightarrow H^\gamma(\Omega)\hookrightarrow L^p(\Omega)\), where
\begin{align}\label{eq:p}
	p\left\{\begin{array}{l}
		<\infty,\;\qquad\text{for }d=2,\; \gamma=1,\\[2mm]
		=\frac{2d}{d-2\gamma},\quad\;\text{for } d=2,\ 0<\gamma<1,\ \text{or}\quad d=3,\ 0<\gamma\le 1,
	\end{array}\right.
\end{align}
and by the assumption on \(\xi_{\gamma}\) in \eqref{eq:xi_gamma}, it is straightforward to verify that
\begin{align*}
	\frac{2d}{d-2\gamma}>\frac{2d\xi_{\gamma}}{d+2(1-\gamma)}.
\end{align*}
Therefore,
\begin{align*}
	\|f(u)\|_{\mathcal{Y}^{\gamma-1}}\lesssim 1+\|u\|_{L^{\frac{2d\xi_{\gamma}}{d+2(1-\gamma)}}}^{\xi_{\gamma}}\lesssim 1+\|u\|_{H^\gamma}^{\xi_{\gamma}}\lesssim 1+\|u\|_{\mathcal{Y}^{\gamma}}^{\xi_{\gamma}},
\end{align*}
which proves \eqref{eq:f_gamma-1}.

 {\sc Proof of \eqref{eq:fu-fv_gamma}}:
We only prove the cases \(d=2,\ 0<\gamma<1\) and \(d=3,\ 0<\gamma\le 1\), since the case \(d=2,\ \gamma=1\) is analogous. Note that
\begin{align*}
	f(u)-f(v)=(u-v)\int_{0}^{1}f^{\prime}(\theta u +(1-\theta)v)\,d \theta.
\end{align*}
By \eqref{eq:growth_condition},
\begin{align}\label{eq:taylor_fu-fv}
	|f(u)-f(v)|\lesssim \Big(1+|u|^{\xi_{\gamma}-1}+|v|^{\xi_{\gamma}-1}\Big)|u-v|.
\end{align}
As before, \(L^{\frac{2d}{d+2(1-\gamma)}} \hookrightarrow (H_0^{1-\gamma})^*=H^{\gamma-1}\), and hence
\begin{align*}
	\|f(u)-f(v)\|_{\mathcal{Y}^{\gamma-1}}
	&\lesssim \|f(u)-f(v)\|_{L^{\frac{2d}{d+2(1-\gamma)}}(\Omega)}\\
	&\lesssim \Big(1+\||u|^{\xi_{\gamma}-1}\|_{L^q(\Omega)}+\||v|^{\xi_{\gamma}-1}\|_{L^q(\Omega)}\Big)\|u-v\|_{L^p(\Omega)},
\end{align*}
where \(p=\frac{2d}{d-2\gamma}\) and
\begin{align*}
	\frac{1}{q}=\frac{d+2(1-\gamma)}{2d}-\frac{1}{p}=\frac{1}{d},
\end{align*}
so that \(q=d\). By \eqref{eq:xi_gamma}, \(q(\xi_{\gamma}-1)<p\). Therefore, by the Sobolev embedding \(\mathcal{Y}^{\gamma}\hookrightarrow H^\gamma(\Omega)\hookrightarrow L^p(\Omega)\),
\begin{align*}
	\|f(u)-f(v)\|_{\mathcal{Y}^{\gamma-1}}
	&\lesssim \Big(1+\|u\|_{L^p(\Omega)}^{\xi_{\gamma}-1}+\|v\|_{L^p(\Omega)}^{\xi_{\gamma}-1}\Big)\|u-v\|_{L^p(\Omega)}\\
	&\lesssim \Big(1+\|u\|_{\mathcal{Y}^{\gamma}}^{\xi_{\gamma}-1}+\|v\|_{\mathcal{Y}^{\gamma}}^{\xi_{\gamma}-1}\Big)\|u-v\|_{\mathcal{Y}^{\gamma}},
\end{align*}
which proves \eqref{eq:fu-fv_gamma}.

{\sc Proof of \eqref{eq:fu_H-1}--\eqref{eq:fuv}}: Since the inequalities \eqref{eq:fu_H-1}--\eqref{eq:fuv} can be proved in a similar way, 
we present only one representative case, namely \eqref{eq:fu-fv}. 

Since \(\mathcal{Y}^{-1} = H^{-1}(\Omega)\), the Sobolev embeddings yield 
\(L^{p_1}(\Omega)\hookrightarrow H^{-1}(\Omega)\), where
\begin{align}\label{eq:p1p2}
	p_{1}=\left\{\begin{array}{l}
		1+,\qquad\;\text{for}\quad d=2,\\[2mm]
		\frac{2d}{d+2},\qquad \text{for}\quad d= 3.
	\end{array}\right.
\end{align}
Thus, by \eqref{eq:taylor_fu-fv}:  
\begin{align*}
	\|f(u)-f(v)\|_{\mathcal{Y}^{-1}}&\lesssim \|f(u)-f(v)\|_{L^{p_{1}}}\\
	&\lesssim \Big(1+\left\||u|^{\xi_{\gamma}-1}\right\|_{L^{q_1}(\Omega)}+\left\||v|^{\xi_{\gamma}-1}\right\|_{L^{q_1}(\Omega)}\Big)\cdot\|u-v\|_{L^2(\Omega)},
\end{align*}
where H\"older's inequality is used with
\begin{align}\label{eq:q_1}
	\frac{1}{q_1}=\frac{1}{p_1}-\frac{1}{2},\quad \text{that is}\quad q_1=\frac{2p_1}{2-p_1}.
\end{align}
Hence,
\begin{align}\label{eq:estimate_fu-fv}
	\|f(u)-f(v)\|_{\mathcal{Y}^{-1}}\lesssim& \big(1+\big\|u\big\|^{\xi_{\gamma}-1}_{L^{q_1(\xi_{\gamma}-1)}}+\big\|v\big\|^{\xi_{\gamma}-1}_{L^{q_1(\xi_{\gamma}-1)}}\Big)\cdot\|u-v\|_{L^2}.
\end{align}
When \(d=3\), combining \eqref{eq:p1p2} and \eqref{eq:q_1} gives
\begin{align*}
	q_1(\xi_{\gamma}-1)=d(\xi_{\gamma}-1)<\frac{2d}{d-2\gamma}.
\end{align*}
Recalling the Sobolev embedding  \(\mathcal{Y}^\gamma \hookrightarrow H^\gamma(\Omega) \hookrightarrow L^p(\Omega)\) 
with \(p\) given in \eqref{eq:p}, we obtain
\begin{align*}
	\|f(u)-f(v)\|_{\mathcal{Y}^{-1}}\lesssim \big(1+\big\|u\big\|^{\xi_{\gamma}-1}_{H^\gamma}+\big\|v\big\|^{\xi_{\gamma}-1}_{H^\gamma}\Big)\cdot\|u-v\|_{L^2}\lesssim C\Big(\|u\|_{\mathcal{Y}^{\gamma}},\|v\|_{\mathcal{Y}^{\gamma}} \Big)\|u-v\|_{\mathcal{Y}^{0}}.
\end{align*}
This proves \eqref{eq:fu-fv} for \(d=3\); the case \(d=2\) can be verified analogously.

{\sc Proof of \eqref{eq:fu-fv_gamma-1}:} From the Sobolev embedding \(H^{1-\gamma}_0 \hookrightarrow L^{\frac{2d}{d-2(1-\gamma)}}\) and the boundedness estimate in \eqref{eq:taylor_fu-fv}, we have
\begin{align}\label{eq:estimate_fu-fv_1}
	\|f(u)-f(v)\|_{\mathcal{Y}^{\gamma-1}}&\lesssim \|f(u)-f(v)\|_{H^{\gamma-1}}\lesssim\|f(u)-f(v)\|_{L^{\frac{2d}{d+2(1-\gamma)}}}\notag\\
	&\lesssim \left\|\Big(1+|u|^{\xi_{\gamma}-1}+|v|^{\xi_{\gamma}-1}\Big)\cdot|u-v|\right\|_{L^{\frac{2d}{d+2(1-\gamma)}}}\notag\\
	&\leq \left\|\Big(1+|u|^{\xi_{\gamma}-1}+|v|^{\xi_{\gamma}-1}\Big)\cdot|u-v|^{1-\alpha}\right\|_{L^{\theta_1}}\||u-v|^{\alpha}\|_{L^{\theta_2}},
\end{align}
for $\alpha\in(0,1]$, where H\"older's inequality is applied in the last step with exponents satisfying
\begin{align}\label{eq:theta1theta2}
	\frac{1}{\theta_1}+\frac{1}{\theta_2}=\frac{d+2(1-\gamma)}{2d}\quad \text{and} \quad \theta_2=\frac{2}{\alpha}.
\end{align}
Since
\begin{align*}
	\Big(1+|u|^{\xi_{\gamma}-1}+|v|^{\xi_{\gamma}-1}\Big)\cdot|u-v|^{1-\alpha}\lesssim 1+|u|^{\xi_\gamma-\alpha}+|v|^{\xi_\gamma-\alpha},  
\end{align*}
we deduce from \eqref{eq:estimate_fu-fv_1} and \eqref{eq:theta1theta2} that
\begin{align}\label{eq:estimate_fu-fv_2}
	\|f(u)-f(v)\|_{\mathcal{Y}^{\gamma-1}}\lesssim \left(1+\big\|u\big\|^{\xi_{\gamma}-\alpha}_{L^{(\xi_{\gamma}-\alpha)\theta_1}}+\big\|v\big\|^{\xi_{\gamma}-\alpha}_{L^{(\xi_{\gamma}-\alpha)\theta_1}}\right)\|u-v\|_{L^2}^{\alpha}.
\end{align}
We first consider the case \((d,\gamma)\ne(2,1)\). Set
\[
\eta_\gamma :=
\frac{d-2\gamma}{2\gamma}
\left(
\frac{d+2(1-\gamma)}{d-2\gamma}-\xi_\gamma
\right)>0 ,
\]
where the positivity follows from the growth condition on \(\xi_\gamma\).
If \(0<\eta_\gamma\le1\), we choose \(\alpha=\eta_\gamma\). A direct
calculation using \eqref{eq:theta1theta2} gives
\[
(\xi_\gamma-\alpha)\theta_1
=
\frac{2d}{d-2\gamma}.
\]
Hence, by the Sobolev embedding \(\mathcal Y^\gamma
\hookrightarrow
H^\gamma(\Omega)
\hookrightarrow
L^{\frac{2d}{d-2\gamma}}(\Omega)\), we deduce from \eqref{eq:estimate_fu-fv_2} that
\[
\|f(u)-f(v)\|_{\mathcal Y^{\gamma-1}}
\lesssim
C\bigl(\|u\|_{\mathcal Y^\gamma},\|v\|_{\mathcal Y^\gamma}\bigr)
\|u-v\|_{\mathcal Y^0}^{\alpha}.
\]

It remains, still with \((d,\gamma)\ne(2,1)\), to consider the case
\(\eta_\gamma>1\). We choose \(\alpha=1\). Then
\(
\theta_1=\frac{2d}{2-2\gamma},
\)
with the usual interpretation \(\theta_1=\infty\) if \(\gamma=1\).
Moreover, the condition \(\eta_\gamma>1\) implies
\[
\xi_\gamma<
\frac{d+2-4\gamma}{d-2\gamma},
\qquad
(\xi_\gamma-1)\theta_1
<
\frac{2d}{d-2\gamma}.
\]
Thus the same argument based on \eqref{eq:estimate_fu-fv_2} and the
Sobolev embedding above yields
\[
\|f(u)-f(v)\|_{\mathcal Y^{\gamma-1}}
\lesssim
C\bigl(\|u\|_{\mathcal Y^\gamma},\|v\|_{\mathcal Y^\gamma}\bigr)
\|u-v\|_{\mathcal Y^0}.
\]

Finally, let \(d=2\) and \(\gamma=1\). We choose \(\alpha=\frac12\).
Then \(\theta_1<\infty\), and, since \(\xi_\gamma<\infty\),
\[
(\xi_\gamma-\alpha)\theta_1<\infty.
\]
Using the two-dimensional Sobolev embedding \(H^1_0(\Omega)\hookrightarrow L^p(\Omega)\) for every \(p<\infty\), and applying \eqref{eq:estimate_fu-fv_2}, we again obtain
\[
\|f(u)-f(v)\|_{\mathcal Y^{0}}
\lesssim
C\bigl(\|u\|_{\mathcal Y^1},\|v\|_{\mathcal Y^1}\bigr)
\|u-v\|_{\mathcal Y^0}^{1/2}.
\]
Combining the above cases proves \eqref{eq:fu-fv_gamma-1}.

\section{A regularity result}

\begin{lemma}
	Let $0<a<b$ and consider the annulus around 0 denoted by $A_{a,b} = B_b(0) \setminus \overline{B_a}(0) \subset \mathbb{R}^2$. Further, let 
	$f \in W^{1,p}(A_{a,b})$
	 be radially symmetric.
	  Then, it holds
	\begin{equation}
		\| f \|_{H^{s}(A_{a,b})} \lesssim \| f \|_{W^{1,p}(A_{a,b})}
	\end{equation}
	for $s < \frac32 - \frac{1}{p}$ and $p \in (1,2)$.
\end{lemma}

\begin{proof}
	We only prove this for smooth functions $f$. To each $f$, we associate a profile $\varphi$ such that $f(x) = \varphi(|x|)$ for all $x \in A_{a,b}$. A straightforward calculation, then shows that there are constants only depending on $a$ and $b$ such that
	%
	\begin{subequations} \label{eq:norm_eq_radial}
	\begin{align}
		c_0 \|f\|_{L^p(A_{a,b})} &\leq  \|\varphi\|_{L^p(a,b)}  \leq C_0   \|f\|_{L^p(A_{a,b})} \label{eq:norm_eq_radial_1}
		\\
		c_1 \|\nabla f\|_{L^p(A_{a,b})} &\leq  \|\varphi' \|_{L^p(a,b)} 
		 \leq C_1   \|\nabla f\|_{L^p(A_{a,b})}. \label{eq:norm_eq_radial_2}
		\end{align}
	\end{subequations}
	%
	An interpolation argument for $p=2$ 
	in the first estimates 
	in \eqref{eq:norm_eq_radial_1} and \eqref{eq:norm_eq_radial_2}
	then shows that there
	is
	a constant which only depends on $a$ and $b$ such that for $s \in (0,1)$
	\begin{align}
	c_s \|f\|_{H^s(A_{a,b})} 
	&\leq  \|\varphi\|_{H^s(a,b)}  .
\end{align}
	Thus, the claim follows from the one dimensional Sobolev embedding
	$W^{1,p}(a,b) \hookrightarrow H^s(a,b)$ for $s < \frac32 - \frac{1}{p}$ and the
	 second estimates 
	 in \eqref{eq:norm_eq_radial_1} and \eqref{eq:norm_eq_radial_2}.
\end{proof}

\begin{ciao} \label{cor:reg_of_u_alpha}
	The functions defined in \eqref{eq:init_half_sphere} defined for $\alpha \in (0,\frac12]$ is in $H^s(\Omega)$ for $s < \frac12 + \alpha$.
\end{ciao}

\begin{proof}
	We first consider a radial partition of unity $\psi_1+\psi_2=1$ 
	in a neighborhood of the circle $\{|x|=\frac14\}$, where 
	$\psi_1$ is compactly supported in an annulus $A_{a,b}$ with $0<a<\frac14<b<1$, and is equal to one in a neighborhood of $\{|x|=\frac14\}$.
	Then, we write 
	\begin{equation}
		u_\alpha^0 = \psi_1 u_\alpha^0 + \psi_2 u_\alpha^0 = w_1 + w_2,
	\end{equation}
	and directly see that $w_2$ is a smooth function. For $w_2$ one can easily compute that 
	$w_2 \in W^{1,p}(\Omega)$ for $p < \frac{1}{1-\alpha}$. 
	Further, since $w_1$ is compactly supported inside the annulus, we conclude by the above lemma that
	\begin{equation}
		\| w_1 \|_{H^s(\Omega)} \lesssim \| w_1 \|_{H^s(A_{a,b})} \lesssim
		\| w_1 \|_{W^{1,p}(A_{a,b})} =   \| w_1 \|_{W^{1,p}(\Omega)}.
	\end{equation}
	Finally, as $s < \frac32 - \frac{1}{p} < \frac12 + \alpha$, we obtain the claim. 
\end{proof}

\end{document}